\documentclass[acmsmall,nonacm,screen]{acmart}
\AtBeginDocument{%
  \providecommand\BibTeX{{%
    \normalfont B\kern-0.5em{\scshape i\kern-0.25em b}\kern-0.8em\TeX}}}

\acmJournal{TOMS}
\acmVolume{0}
\acmNumber{0}
\acmArticle{0}
\acmMonth{1}

\usepackage{amsmath}
\usepackage{array}
\usepackage{graphicx}
\usepackage{epsfig}
\usepackage{float}
\usepackage{enumitem}  
\usepackage{stmaryrd}
\usepackage{standalone}
\usepackage{xcolor}
\usepackage{listings}
\usepackage{bm}
\usepackage{cleveref}
\usepackage{mathtools}
\usepackage{wrapfig}
\usepackage{algorithm}
\usepackage{algorithmic}
\usepackage{pgfplots}

\definecolor{codegreen}{rgb}{0,0.6,0}
\definecolor{codegray}{rgb}{0.5,0.5,0.5}
\definecolor{codepurple}{rgb}{0.58,0,0.82}
\definecolor{backcolor}{rgb}{0.95,0.95,0.92}
\definecolor{codeblue}{rgb}{0.05,0.3,0.7}
\definecolor{codedarkblue}{rgb}{0.3,0.1,0.7}
\definecolor{codered}{rgb}{0.54,0.03,0.06}
\definecolor{codedarkmagenta}{rgb}{0.99, 0.3, 0.4}

\lstdefinestyle{mystylepython}{
    language=Python,
    frame=single,
    backgroundcolor=\color{backcolor},   
    commentstyle=\color{codegreen},
    keywordstyle={\color{magenta}},
    keywordstyle = [2]{\color{codeblue}},
    keywordstyle = [3]{\color{codedarkblue}},
    keywordstyle = [4]{\color{codepurple}},
    keywordstyle = [5]{\color{codedarkmagenta}},
    numberstyle=\tiny\color{codegray},
    stringstyle=\color{codered},
    basicstyle=\ttfamily\footnotesize,
    breakatwhitespace=false,         
    breaklines=true,                 
    captionpos=b,                    
    keepspaces=true,
    numbersep=5pt,                  
    showspaces=false,                
    showstringspaces=false,
    showtabs=false,                  
    tabsize=2,
    morekeywords={map, range, },
    morekeywords=[2]{dx, ds, dx_, dS, div, dot, grad, inner, UnitSquareMesh, UnitCubeMesh, FunctionSpace, TrialFunction, TestFunction, MeshFunction, CompiledSubDomain, jump, assemble, sym, Vector, MPI, as_backend_type, Measure, Constant, FacetArea, CellDiameter, Expression, avg, sqrt},
    morekeywords=[3]{OuterNormal, Trace, block_form, block_mat, block_vec, ii_assemble, ii_convert, HXDiv, Jacobi, RA, block_diag_mat, MinRes, EmbeddedMesh, AMG, block_base, Precond, block_to_haz, ConjGrad, SOR, Circle, Average, PETSc_to_dCSRmat, LU},
    morekeywords=[4]{AMG_param, create_precond_ra, param_amg_set_dict, apply_precond, create_dvector, dvec_create_p, fenics_metric_amg_solver, create_precond_amg},
    morekeywords=[5]{UA_AMG, NL_AMLI_CYCLE, SMOOTHER_GS, DIRECT, VMB}
}

\lstdefinestyle{mystyleC}{
    language=C,
    frame=single,
    backgroundcolor=\color{backcolor},   
    commentstyle=\color{codegreen},
    keywordstyle={\color{magenta}},
    keywordstyle = [2]{\color{codepurple}},
        keywordstyle = [3]{\color{codedarkmagenta}},
    numberstyle=\tiny\color{codegray},
    stringstyle=\color{codered},
    basicstyle=\ttfamily\footnotesize,
    breakatwhitespace=false,         
    breaklines=true,                 
    captionpos=b,                    
    keepspaces=true,
    numbersep=5pt,                  
    showspaces=false,                
    showstringspaces=false,
    showtabs=false,                  
    tabsize=2,
    morekeywords=[2]{create_precond_ra, param_amg_set_dict, apply_precond, create_dvector, dvec_create_p, fenics_metric_amg_solver, create_precond_amg, precond_amg, param_prec_to_amg, array_cp, dvec_set, mgcycle, precond_ra_fenics, dcsr_pcg, array_ax, array_axpy, fenics_metric_amg_solver, param_input_init, param_input, linear_solver_bdcsr_krylov_metric_amg, linear_solver_bdcsr_krylov, precond_bdcsr_metric_amg_symmetric, smoother_dcsr_Schwarz_forward, mgcycle_bdcsr, smoother_dcsr_Schwarz_backward, param_amg_init},
    morekeywords=[3]{AMG_data, AMG_param, dCSRmat, dvector, block_dCSRmat, input_param,  Schwarz_param, Schwarz_data, AMG_data_bdcsr, REAL, INT, precond_data, UA_AMG, NL_AMLI_CYCLE, SMOOTHER_GS, DIRECT, VMB, precond}
}

\newcommand{\half}{{\frac{1}{2}}}
\newcommand{\RR}{{\mathbb{R}}}
\newcommand{\CC}{{\mathbb{C}}}
\renewcommand{\div}{{\operatorname{div}}}
\newcommand{\diag}{{\operatorname{diag}}}

\newcommand{\vertiii}[1]{{\vert\kern-0.25ex\vert\kern-0.25ex\vert #1 
    \vert\kern-0.25ex\vert\kern-0.25ex\vert}}

\begin{document}

\title{HAZniCS -- Software Components for Multiphysics Problems}

\author{Ana Budi\v{s}a}
\email{ana@simula.no}
\affiliation{%
  \institution{Simula Research Laboratory}
  \streetaddress{P.O. Box 134}
  \postcode{1325}
  \city{Lysaker}
  \country{Norway}
}

\author{Xiaozhe Hu}
\affiliation{%
  \institution{Department of Mathematics, Tufts University}
  \streetaddress{503 Boston Avenue}
  \city{Medford}
  \postcode{02155}
  \state{Massachusetts}
  \country{USA}
  }
\email{xiaozhe.hu@tufts.edu}

\author{Miroslav Kuchta}
\affiliation{%
  \institution{Simula Research Laboratory}
  \streetaddress{P.O. Box 134}
  \postcode{1325}
  \city{Lysaker}
  \country{Norway}
  }
\email{miroslav@simula.no}

\author{Kent--Andr\'{e} Mardal}
\affiliation{%
 \institution{Department of Mathematics, University of Oslo}
 \streetaddress{P.O. Box 1053, Blindern}
 \postcode{0316}
 \city{Oslo}
 \country{Norway}}
\email{kent-and@math.uio.no}

\author{Ludmil T. Zikatanov}
\affiliation{%
  \institution{Department of Mathematics, Penn State}
  \streetaddress{239 McAllister Building}
  \city{University Park}
  \postcode{16802}
  \state{Pennsylvania}
  \country{USA}}
\email{ludmil@psu.edu}

\renewcommand{\shortauthors}{A. Budi\v{s}a, X. Hu, M. Kuchta, K.-A. Mardal and L.~T.~Zikatanov}

\begin{abstract}
  We introduce the software toolbox HAZniCS for solving interface-coupled multiphysics problems. 
  HAZniCS is a suite of modules that combines the well-known FEniCS
  framework for finite element discretization 
  with solver and graph library HAZmath. The focus of the paper is on the design and implementation of a pool of robust and efficient solver algorithms which tackle issues related to the complex interfacial 
  coupling of the physical problems often encountered in applications in brain biomechanics. The robustness and efficiency of the numerical algorithms and methods is shown in several numerical examples, namely the Darcy-Stokes equations that model flow of cerebrospinal fluid in the human brain and the mixed-dimensional model of electrodiffusion in the brain tissue.
\end{abstract}

%


\maketitle

\section{Introduction} \label{sec:introduction}

The present paper aims to introduce a novel collection of tools for interface coupled multiphysics problems modeled by partial differential equations (PDEs). The interface is a main driver of the processes in a way that strategies relying on decoupled single-physics problems typically suffer from slow convergence. Furthermore, we target multiphysics problems with geometrically complex interfaces and slow dynamics -- promoting monolithic solvers. Specifically, we exploit fractional operators and low-order interface perturbations as preconditioning techniques.  

Fractional operators appear naturally on interfaces in multiphysics problems. One common approach has been using Poincar\'{e}-Steklov operators for fluid-structure interaction problems~\cite{deparis2006fluid, quarteroni1991theory, agoshkov1988poincare}, which exploits Dirichlet-to-Neumann mappings. As the Poincar\'{e}-Steklov operator takes functions in the fractional Sobolev space $ H^{1/2} $ to functions in its dual $H^{-1/2}$, it is equivalent to a fractional Laplacian operator $ (-\Delta)^{1/2} $. However, the Poincar\'{e}-Steklov operator is not sufficient for parameter-dependent problems as it is sensitive to problem parameters, and often many sub-iterations are required. More sophisticated techniques that include problem parameters such as Robin-to-Dirichlet, -Neumann, or -Robin maps have been explored ~\cite{badia2009robin}, but the approach still requires tuning. We remark that the Poincar\'{e}-Steklov operator involves the extension to a domain in a higher dimension and is, as such, computationally expensive. However, the computational complexity is usually the same as the involved single-physics problems. 

As an alternative or generalization of Poincar\'{e}-Steklov operators, several recent papers~\cite{hornkjol2021, boon2021robust, Kuchta2021, holter2021robust} have considered multiphysics problems and derived order-optimal and parameter-robust algorithms. They are obtained by exploiting fractional Laplacians (or sums thereof) and metric terms on the interface. Here, the fractional Laplacians arise naturally due to trace operators appearing in the coupling conditions, e.g., conservation of mass, that connect the unknowns of the different single-physics problems. We note that the fractional operators arise both when Lagrange multipliers are used to prescribe the interface conditions, e.g., \cite{layton2002coupling, holter2020}, and when they are avoided \cite{hornkjol2021, boon2021robust}. The metric terms then arise because interface conditions, such as the balance of forces, are often expressed in terms of differences of a quantity (e.g., displacement) across the interface rather than the quantity itself~\cite{boon2021robust, DAngelo2008, Kuchta2021}. Recently, fast solution algorithms for fractional Laplacians have been proposed based on multilevel approaches~\cite{Baerland2019,baerland2019auxiliary,fuhrer2022multilevel, bramble2000computational,ZhaoHuCaiKarniadakis2017a} and rational approximations~\cite{harizanov2018optimal, harizanov2020numerical, harizanov2022rational}. Here, we explore the latter for sums of fractional Laplacians. In addition to rational approximations, we will consider multilevel algorithms that work robustly in the presence of strong metric terms at interfaces. That is multilevel algorithms with a space decomposition aware of the metric kernel.  

The software tools we developed aim to solve computational mesoscale multiphysics problems. By computational mesoscale, in this context, we refer to problems in the range of a few hundred thousand to tens of millions of degrees of freedom. These problems do not require parallel computing, but they may benefit significantly from advanced algorithms. The collection of tools presented in this paper are FEniCS ~\cite{fenics_book} add-ons for block assembly~\cite{fenicsii} and block preconditioning~\cite{mardal2012block} combined with a flexible algebraic multigrid (AMG) toolbox, implemented in C, called HAZmath~\cite{hazmath}. Hence, we have named the tool collection HAZniCS. One of the reasons for developing HAZniCS is precisely the mentioned flexibility and variety of the implementation of the AMG method in HAZmath. It allows us to easily modify available linear solvers and preconditioners or create new model-specific solvers for the multiphysics problems at hand. Additionally, with HAZniCS, we provide another wide range of efficient computational methods for solving PDEs with FEniCS, but also a bridge to Python for HAZmath to be used with other PDE simulation tools. Further in the paper, we highlight with a series of code snippets the implementation of several solvers, namely the aggregation-based and metric-perturbed AMG methods and the rational approximation method.

Moreover, we consider a series of examples of multiphysics problems mainly related to biomechanical processes. Namely, we include: (1) a simple three-dimensional example of a elliptic problem on a regular domain, (2) Darcy-Stokes equations describing the interaction of the viscous flow of cerebrospinal fluid flow surrounding the brain and interacting with the porous media flow of interstitial fluid inside the brain, and (3) the mixed-dimensional equations representing electric signal propagation in neurons and
the surrounding matter.

The outline of the current paper is as follows: in \Cref{sec:examples} we introduce the multiphysics models together with the necessary mathematical concepts and numerical methods. \Cref{sec:implementation} focuses on the implementation of those methods and the interface between the software components. In \Cref{sec:results} we present the solver capabilities of our software to simulate relevant biomechanical phenomena. Finally, we draw concluding remarks in \Cref{sec:conclusion}.

\section{Examples} \label{sec:examples}

The following three examples illustrate different single- and multiphysics PDE models, as well as the relevant mathematical and computational concepts. More specifically, the examples provide an overview of iterative methods and preconditioning techniques for interface-coupled problems that lead, e.g., to the utilization of sums of fractional operators weighted by material parameters. Additionally, we include several code snippets in each example that highlight the most important features of the implementation, while the full codes can be found in \cite{github_haznics}.

\subsection{Linear elliptic problem} \label{subsec:example_poisson}

We start with a linear elliptic problem on a three-dimensional (3$d$) regular domain. This example will serve as a baseline for the
solvers in HAZniCS. Our goal is to demonstrate that our solver performance is comparable to other established software. Additionally, the solution methods that we use here will be incorporated and adapted to the multiphysics problems in the later examples.

Let $ \Omega = [0, 1]^3 $ be the unit cube and let $ \partial \Omega $ denote its boundary. Given external force $ f :\Omega \to \RR $ and the boundary data $ g : \partial \Omega \to \RR $, we set to find the solution $ u : \Omega \to \RR $ that satisfies
\begin{subequations}\label{eq:poisson}
	\begin{align}
		- \Delta u + u & = f & \text{ in } \Omega, \label{eq:poisson_a} \\
		\dfrac{\partial u }{\partial \bm{n}} & = g & \text{ on } \partial \Omega. \label{eq:poisson_b}
	\end{align}
\end{subequations}
To solve \eqref{eq:poisson} computationally, we relate to \eqref{eq:poisson} the variational formulation and the discrete problem using the finite element method (FEM). First, let $ L^2 = L^2(\Omega)$ be the space of square-integrable functions on $ \Omega $ and $ H^s = H^s(\Omega) $ the Sobolev spaces with $s$ derivatives in $ L^2 $. The corresponding inner products and norms for any function space $ X $ are denoted with $ (\cdot, \cdot)_X $ and $\| \cdot \|_X $, respectively. Furthermore, we let $ \langle \cdot, \cdot \rangle $ denote a duality pairing between $X'$, the dual of $X$ and $X$.

Now, let $ V_h \subset H^1(\Omega) $ be a finite element space on triangulation of $\Omega$, e.g., of continuous piecewise linear functions ($\mathbb{P}_1$) . A discrete variational formulation of \eqref{eq:poisson} states to find $ u \in V_h $ such that
\begin{equation}
	\label{eq:weak_poisson}
	a(u, v) = \langle f, v \rangle \qquad \forall v \in V_h,
\end{equation} 
with $ a(u, v) = (u, v)_{L^2(\Omega)} + (\nabla u, \nabla v)_{L^2(\Omega)} $ and $ \langle f, v\rangle  = (f, v)_{L^2(\Omega)} + (g, v)_{L^2(\partial \Omega)} $. 

Furthermore, it is important for the solvers to obtain matrix-vector representation of the above system. Let the discrete operator $ A : V_h \to V_h' $ satisfy 
\begin{equation} \label{eq:poisson_discrete_operator}
    \langle A u, v \rangle = a(u, v), \quad u, v \in V_h
\end{equation}
with $ V_h'$ denoting the dual space of $ V_h $. Its actual implementation can be derived as follows. Let $ \psi_i, i = 1, 2, \dots, m $ be the finite element basis functions of $ V_h $. Define matrix $ \mathsf{A} \in \RR^{m \times m} $ and vectors $ \mathsf{f} \in \RR^m $ as
\begin{equation}\label{eq:matrix_rhs}
	(\mathsf{A})_{ij} = \langle A \psi_j, \psi_i \rangle , \quad \text{ and } \quad \mathsf{f}_i = \langle f, \psi_i\rangle, \quad \text{ for } i,j = 1, 2, \dots, m.
\end{equation}
Consequently, we get the discrete system of equations related to \eqref{eq:poisson}, i.e. we aim to solve for $ \mathsf{u} \in \RR^m $ the algebraic system
\begin{equation} \label{eq:poisson_matrix_system}
	\mathsf{A} \mathsf{u} = \mathsf{f}.
\end{equation}
We remark that $\mathsf{u}$ and $\mathsf{f}$ above are both vectors in $\RR^m$, but that $\mathsf{u}$ is in the so-called \emph{nodal} representation, i.e. $u=\sum_i \mathsf{u}_i\psi_i$, while $\mathsf{f}$ is in the \emph{dual} representation~\cite{bramble2019multigrid, mardal_winther_2011}. As such, the matrix	$\mathsf{A}$ maps the nodal representations of  $\RR^m$ to its dual representation.

Since $ \mathsf{A} $ is symmetric positive definite (SPD), we solve \eqref{eq:poisson_matrix_system} with the Conjugate Gradient (CG) method. It is well known that the number of iterations of a Krylov iterative scheme can be bounded in terms of the \emph{condition number} of the system, that is $ \kappa(\mathsf{A}) = \vertiii{\mathsf{A}} \vertiii{\mathsf{A}^{-1}} $ for some matrix norm $ \vertiii{\cdot} $. Therefore, to efficiently solve the problem \eqref{eq:poisson_matrix_system} we want as few iterations as possible and order optimal scalability of the solver with regards to the number of degrees of freedom. To that aim, we introduce a \emph{preconditioner} $ \mathsf{B} $ such that
\begin{equation}\label{eq:poisson_condition_number}
	\kappa(\mathsf{B} \mathsf{A}) = \vertiii{ \mathsf{B} \mathsf{A} } \vertiii{ (\mathsf{B} \mathsf{A})^{-1}} \approx \mathcal{O}(1),
\end{equation}
that is, $ \kappa(\mathsf{B} \mathsf{A}) $ stays bounded from above independently of discretization and other problem parameters. It is important to make sure that $ \mathsf{B} $ maps dual representations of vectors to nodal representations as the preconditioner $ \mathsf{B} $ is an approximation of the inverse of $ \mathsf{A} $.

The previous result is also true in the general case for symmetric operators on Hilbert spaces. Specifically, for $ A $ in \eqref{eq:poisson_discrete_operator} we find an operator $ B : V_h' \to V_h $ such that $ \kappa(BA) $ is uniformly bounded, where the matrix norm is replaced with the operator norm in the space of continuous linear operators defined on $ V_h $. In context of operator preconditioning \cite{mardal_winther_2011}, a common choice for the preconditioner operator $B$ is the \emph{Riesz mapping}, that is
\begin{equation} \label{eq:riesz_map}
	(B f, v)_{V_h} = \langle f, v \rangle, \quad \forall f \in V_h', v \in V_h. 
\end{equation}
The Riesz map guarantees a uniform bound on $ \kappa(B A) $ when $ A $ is a bounded operator that satisfies the inf-sup conditions \cite{babuska_aziz, babuska1971} independent of system parameters, such as in the case of the operator in \eqref{eq:poisson_discrete_operator}. Moreover, we can use any other preconditioner that gives a uniform bound on the condition number. If we find a \emph{spectrally equivalent} operator $ B_{SE} $ such that for parameter-independent constants $ c_1, c_2 > 0 $ it satisfies 
\begin{equation} \label{eq:norm_equivalent}
	c_1 \| v \|^2_{A} \leq \| v \|^2_{B_{SE}^{-1}} \leq c_2 \| v \|^2_{A},
\end{equation}
with $ \| v \|^2_A = \langle A v, v \rangle  $, then we retain a uniform bound on the condition number $ \kappa(B_{SE} A) \leq \frac{c_2}{c_1} \kappa(BA) $. This is relevant when an application of $ B $ on a function in $ V_h'$ is infeasible or inefficient. We want to replace it with a method that applies a spectrally equivalent operation. In the rest of the paper, we will note $ \kappa(B A) $ as the condition number for both operators ($A,\,B$) and their matrix representations ($ \mathsf{A}, \, \mathsf{B} $), clarifying along the way if ambiguity occurs.

In our case, it is well-known that multilevel methods, such as AMG, provide spectrally equivalent and order optimal algorithms for the inverse of discretizations of $ I - \Delta $. Thus, we define the preconditioner for \eqref{eq:poisson_matrix_system} as $ \mathsf{B} = \operatorname{AMG}(\mathsf{A}) $.

The implementation of the elliptic problem in FEniCS follows straightforwardly from the variational formulation \eqref{eq:weak_poisson} and is one of the basic examples of FEniCS software, see \Cref{lst:poisson}.
\begin{lstlisting}[style=mystylepython, caption={Implementation of the linear elliptic problem \eqref{eq:poisson}. Complete code can be
      found in scripts \texttt{HAZniCS-examples/demo\_elliptic*.py}
    }, captionpos=b, label={lst:poisson}]
	from block.iterative import ConjGrad
	from block.algebraic.hazmath import AMG
	from dolfin import *
	
	mesh = UnitCubeMesh(32, 32, 32)
	V = FunctionSpace(mesh, "CG", 1)
	u, v = TrialFunction(V), TestFunction(V)
	f = Expression("sin(pi*x[0])", degree=4)
	
	a = inner(u, v) * dx + inner(grad(u), grad(v)) * dx
	L = inner(f, v) * dx
	A = assemble(a)
	b = assemble(L)
	
	B = AMG(A, parameters={"max_levels": 10, "AMG_type": 1})
	Ainv = ConjGrad(A, precond=B, tolerance=1e-10)
	x = Ainv * b   # Solve for the coefficient vector of foo in V
\end{lstlisting}
For the preconditioner, we utilize the AMG method implemented in HAZmath, available through our HAZniCS library. We describe the AMG method and its implementation in more detail in \Cref{subsec:amg} and showcase the performance of HAZmath AMG as compared with HYPRE \cite{hypre} AMG in \Cref{subsec:poisson_results}.

\subsection{Modeling brain clearance during sleep with Darcy-Stokes equations} \label{subsec:example_ds}
\begin{figure}
  \centering
  \includegraphics[width=0.7\textwidth]{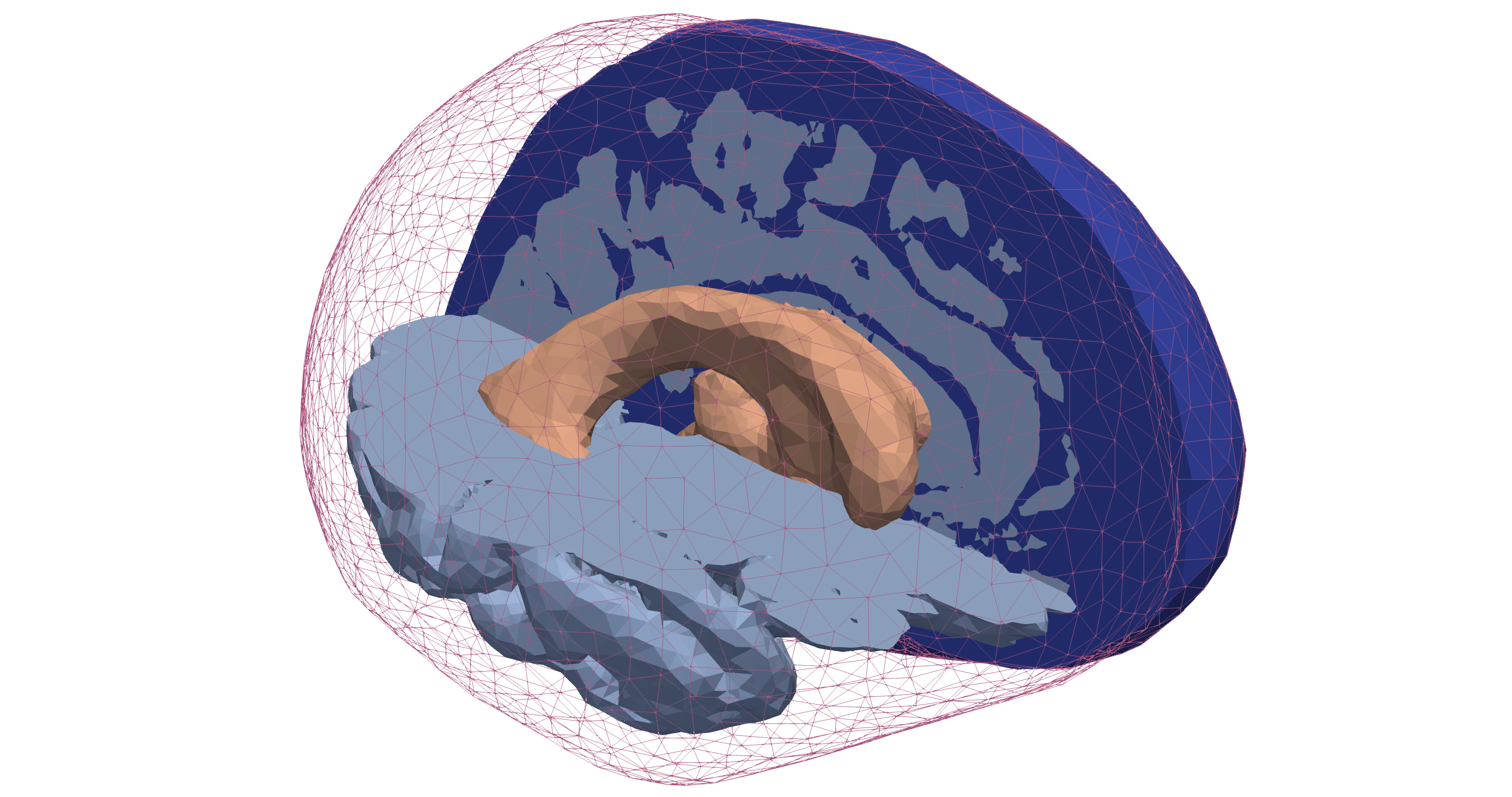}
  \includegraphics[width=0.2\textwidth]{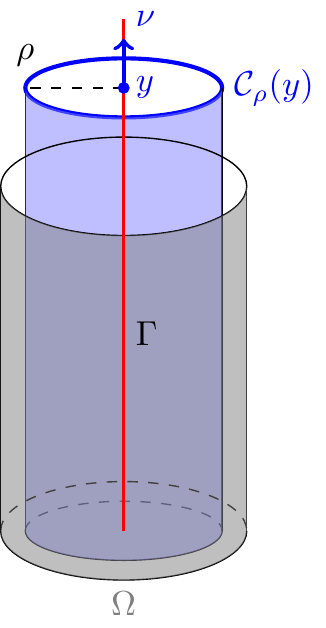}  
  \caption{
    (Left) Geometry and computational mesh from \cite{hornkjol2021} of the Darcy-Stokes model of brain clearance. Mesh and indicator functions for tracking subdomains making up the Darcy-(light blue) and the Stokes domains (dark blue and orange subregions) and their interfaces are generated with SVMTK ~\cite{mri2fem}. (Right) Model reduction from 3$d$-3$d$ to 3$d$-1$d$ problem. Dendrites (in blue) are reduced to their centerline while the coupling with the surrounding $\Omega$ is accounted for by averaging over-idealized cylindrical surfaces with radius $\rho$.
  }
  \vspace{-10pt}
  \label{fig:brain}
\end{figure}

We consider a multiphysics problem arising in modeling processes of waste clearance in the brain during sleep~\cite{xie2013sleep, eide2021sleep} with potential links to the development of Alzheimer's disease. The novel model, called the glymphatic model~\cite{iliff2012paravascular}, states that the viscous flow of cerebrospinal fluid (CSF) is tightly coupled to the porous flow in the brain tissue and that during sleep, in particular, it clears metabolic waste from the brain, for computational models see e.g. ~\cite{constanzo2020, holter2020, hornkjol2021}. To this end we will consider patient-specific geometries generated from MRI images by SVTMK library ~\cite{mri2fem} used in \cite{hornkjol2021}, see \Cref{fig:brain}. Using SVMTK, the segmented brain geometry is enclosed in a thin shell, which, together with the ventricles (the orange subregion in \Cref{fig:brain}), makes up the Stokes domain. We remark that the diameter of the Stokes domain is roughly 15\,cm while the shell thickness is on average 0.8\,mm.

In order to model the waste clearance, let $\Omega_D\subset\mathbb{R}^d$, $d=2, 3$ be the domain of the porous medium flow that represents the brain tissue\footnote{In the context of brain mechanics, the case $d=2$ is relevant, e.g., for the slices of the brain geometry.}, and let $\Omega_S\subset\mathbb{R}^d$ be the domain of viscous flow representing the subarachnoid space around it saturated by CSF. Let $\Gamma$ denote the interface between the domains, which in this case corresponds to the surface of the brain. We then consider the Darcy-Stokes model which seeks to find Stokes velocity $ \bm{u}_S : \Omega_S \to \RR^d $ and pressure $ p_S : \Omega_S \to \RR $, and Darcy velocity $ \bm{u}_D : \Omega_D \to \RR^d $ and pressure $ p_D : \Omega_D \to \RR $ that satisfy
\begin{subequations}
	\label{eq:darcy_stokes}
	\begin{align}
		- \nabla \cdot \bm{\sigma}_S (\bm{u}_S, p_S) & = \bm{f}_S & \text{ in } \Omega_S, \label{eq:darcy_stokes_a} \\
		\nabla \cdot \bm{u}_S & = 0 & \text{ in } \Omega_S, \label{eq:darcy_stokes_b} \\
		\bm{u}_D & = -K \nabla p_D & \text{ in } \Omega_D, \label{eq:darcy_stokes_c} \\
		\nabla \cdot \bm{u}_D & = f_D & \text{ in } \Omega_D, \label{eq:darcy_stokes_d} \\
		\shortintertext{ with interface conditions }
		\bm{u}_S \cdot \bm{n} - \bm{u}_D \cdot \bm{n} & = 0 & \text{ on } \Gamma, \label{eq:darcy_stokes_e} \\
		\bm{n} \cdot \bm{\sigma}_s (\bm{u}_S, p_S) \cdot \bm{n} + p_D & = 0 & \text{ on } \Gamma, \label{eq:darcy_stokes_f} \\
		\bm{n} \cdot \bm{\sigma}_s (\bm{u}_S, p_S) \cdot \bm{\tau} + D \bm{u}_S \cdot \bm{\tau} & = 0 & \text{ on } \Gamma. \label{eq:darcy_stokes_g}
	\end{align}
\end{subequations}
Here, $ \bm{\sigma}_S(\bm{u}_S, p_S) = \mu \nabla\bm{u}_S - p_S I $. We remark that for simplicity, $\bm{\sigma}_S$ is defined in terms of the full velocity gradient and not only its symmetric part, cf. \cite{layton2002coupling}. The parameters $\mu$, $K$, and $ D $ are positive constants related to the problem's physical parameters, i.e., the fluid viscosity, permeability, and the Beavers-Joseph-Saffman (BJS) coefficient. Functions $ \bm{f}_S$ and $ f_D $ represent the external forces. Additionally, $ \bm{n} $ denotes the unit outer normal of $\Omega_S$ and $ \bm{\tau} $ is any unit vector tangent to the interface. In particular, for $d=3$ the condition \eqref{eq:darcy_stokes_g} represents a pair of constraints. Finally, we assume the following boundary conditions
\begin{subequations}
	\label{eq:darcy_stokes_bc}
	\begin{align}
		\bm{u}_S & = \bm{0} & \text{ on } \partial \Omega_{S, D}, \\
		\bm{\sigma}_S (\bm{u}_S, p_S) \cdot \bm{n} & = \bm{g} & \text{ on } \partial \Omega_{S, N},
	\end{align}
\end{subequations}
for $ \partial \Omega_{S, D} \cup \partial \Omega_{S, N} = \partial \Omega_{S} \backslash \Gamma $ and
$ \partial \Omega_{S, D} \cap \partial \Omega_{S, N} = \emptyset $.

To arrive at the finite element formulation of \eqref{eq:darcy_stokes} let us introduce a Lagrange multiplier $ \lambda : \Gamma\to\mathbb{R}$, $\lambda\in \Lambda = \Lambda (\Gamma)$ for enforcing the mass conservation across the interface \eqref{eq:darcy_stokes_e}. In addition we consider conforming discrete subspaces $ \bm{V}_S \times Q_S \subset \bm{H}^1(\Omega_S) \times L^2(\Omega_S) $ and $ \bm{V}_d \times Q_D \subset \bm{H}(\div, \Omega_D) \times L^2(\Omega_D)$ for the Stokes and Darcy subproblems respectively. In the following numerical examples, such spaces are constructed by Taylor-Hood ($\mathbb{P}_{2}$-$\mathbb{P}_{1}$) elements and lowest order Raviart-Thomas elements $\mathbb{R}\mathbb{T}_{0}$ paired with discontinuous Lagrange elements $\mathbb{P}^{\text{disc}}_{0}$ for $Q_D$. The multiplier space is discretized by $\mathbb{P}^{\text{disc}}_{0}$ elements.
The variation formulation of \eqref{eq:darcy_stokes} then states to find $ (\bm{u}_S, p_S, \bm{u}_D, p_D, \lambda) \in \bm{V}_s \times Q_s \times \bm{V}_d \times Q_d  \times \Lambda $ that satisfy
\begin{equation}
	\label{eq:ds_system}
	\underbrace{
	\begin{pmatrix}
		-\mu \nabla\cdot \nabla  + D T_{\bm{\tau}}' T_{\bm{\tau}} & -\nabla &  & & T_{\bm{n}}' \\
		\nabla \cdot & & & & \\
		& {K}^{-1} I & & - \nabla & - T_{\bm{n}}' \\
		& \nabla \cdot & & & \\
		T_{\bm{n}} & -T_{\bm{n}} & & &
	\end{pmatrix}
	}_{A}
	\underbrace{
	\begin{pmatrix}
		\bm{u}_S \\ p_S \\ \bm{u}_D \\ p_D \\ \lambda
	\end{pmatrix}
	}_{x}
	=
	\underbrace{
		\begin{pmatrix}
			\bm{f}_s \\ 0 \\ 0 \\ f_d \\ 0
		\end{pmatrix}
	}_{b}.
\end{equation}
The operators $ T_{\bm{n}} $ and $ T_{\bm{\tau}} $ denote the normal and the tangential trace operators on $ \Gamma $. 
\begin{lstlisting}[style=mystylepython, caption={
      Implementation of the bilinear form of the Darcy-Stokes system \eqref{eq:ds_system}. Complete code can be
      found in scripts \texttt{HAZniCS-examples/demo\_darcy\_stokes*.py}
    }, captionpos=b, label={lst:ds_system}]
  # Mesh definitions, FEM space W declaration [...]
  uS, pS, uD, pD, lmbda = map(TrialFunction, W)
  vS, qS, vD, qD, dlmbda = map(TestFunction, W)
  
  TuS, TvS = (Trace(f, Gamma) for f in (uS, vS))
  TuD, TvD = (Trace(f, Gamma) for f in (uD, vD))
  dx_ = Measure('dx', domain=Gamma)
  
  a = block_form(W, 2)
  # Stokes
  a[0][0] = inner(mu * grad(uS), grad(vS)) * dx + 
  D * inner(dot(TvS, tau_), dot(TuS, tau_)) * dx_
  # Stabize Crouzeix-Raviart
  if VS.ufl_element().family() == 'Crouzeix-Raviart':
      tdim = meshS.topology().dim()
      hS = avg(FacetArea(meshS)) if tdim == 2 else sqrt(avg(FacetArea(meshS)))
      a[0][0] += (mu*Constant(10)/hS)*inner(jump(uS), jump(vS))*dS
        
  a[0][1] =  -inner(pS, div(vS)) * dx
  a[0][4] = inner(lmbda, dot(TvS, n_)) * dx_
  # Darcy
  a[2][2] = K ** -1 * inner(uD, vD) * dx
  a[2][3] = -inner(pD, div(vD)) * dx
  a[2][4] = -inner(lmbda, dot(TvD, n_)) * dx_
  # Symmetrize [...] 
\end{lstlisting}
In \Cref{lst:ds_system} we see that the block structure of the problem operator $A$ in \eqref{eq:ds_system} is mirrored in its implementation in FEniCS/cbc.block and that trace operators are implemented using \cite{fenicsii}. In addition, \Cref{lst:ds_system} includes interior facet stabilization employed when the $H^1$-nonconforming Crouzeix-Raviart ($\mathbb{C}\mathbb{R}_1$) elements are used to discretize the Stokes velocity. 

A parameter robust preconditioner for Darcy-Stokes problem \eqref{eq:ds_system} is derived in \cite{holter2020} within the framework of operator preconditioning \cite{mardal_winther_2011}. Specifically, \cite{holter2020} propose the following block-diagonal operator
%
\begin{equation}
	\label{eq:ds_precond}
	B = \begin{pmatrix}
		-\mu \nabla\cdot\nabla + D T_t' T_t &  &  &  &  \\
		& \!\!\!\mu^{-1} I &  &  &  \\
		& & \!\!\!K^{-1} (I - \nabla \nabla \cdot) & & \\
		& & & \!\!\!KI & \\
		& & & & \!\!\!\mu^{-1} (- \Delta+I)^{-\frac{1}{2}} + K (- \Delta+I)^{\frac{1}{2}}
	\end{pmatrix}^{-1},
\end{equation}
which is a Riesz map with respect to the inner products of parameter-weighted Sobolev spaces. In particular, the preconditioner for the $\Lambda$-block reflects posing of the Lagrange multiplier in the intersection space $ \mu^{-1/2} H^{-\half}(\Gamma) \cap K^{1/2} H^{\half}(\Gamma) $. We also note that the $\bm{V}_S$-block of the preconditioner $B$ is identical to the $(0, 0)$-component of the problem operator $A$ in \eqref{eq:ds_system}.

Implementation of the preconditioner within HAZniCS is given in \Cref{lst:ds_prec}. First, we recognize that the preconditioner extracts the relevant block from the operator $A$ to construct the Stokes velocity preconditioner while the remaining inner product operators are assembled (as they are not part of $A$). The option to extract or assemble the (auxiliary) operators to define the preconditioner is a powerful feature of HAZniCS/cbc.block. Similar functionality \cite{kirby2018solver} enables flexible specification of preconditioners in the Firedrake \cite{firedrake} finite element library. 

In \Cref{lst:ds_prec}, we use scalable algorithms from HAZniCS, and PETSc \cite{petsc-user-ref} to approximate the inverses of all the blocks. Algebraic multilevel schemes are used for the Riesz maps on $\bm{V}_S$ and $\bm{V}_D$ where in particular, in the latter, the HAZniCS preconditioner class \texttt{HXDiv} implements the auxiliary space method for $H(\div)$ problems \cite{Kolev2012}. Riesz maps due to the $L^2$ inner products on the pressure spaces are realized via simple iterative schemes such as the symmetric successive over-relaxation SSOR. Finally, the preconditioner for the Lagrange multiplier, which involves the inverse of a sum of fractional operators, is solved with a rational approximation that employs AMG internally. These components will be discussed in detail in \Cref{sec:implementation}. 

\begin{lstlisting}[style=mystylepython, caption={Scalable implementation of preconditioner for Darcy-Stokes problem \eqref{eq:ds_system}.
Complete code can be found in scripts \texttt{HAZniCS-examples/demo\_darcy\_stokes*.py}
    }, captionpos=b, label={lst:ds_prec}]
  from block.algebraic.hazmath import RA, AMG, HXDiv  
  from block.algebraic.petsc import SOR
  
  VS, QS, VD, QD, Q = W
  # Define SPD operators defining inner products on the spaces  
  # Stokes velocity block is taken from the system matrix
  B0 = AA[0][0]
  # L^2 inner product on QS
  B1 = assemble((1/mu)*inner(TrialFunction(QS), TestFunction(QS))*dx)
  # H(div) inner product on VD
  uD, vD = TrialFunction(VD), TestFunction(VD)
  B2 = assemble((1 / K) * (inner(u2, v2) * dx + inner(div(u2), div(v2)) * dx))
  # L^2 inner product on QD
  B3 = assemble(K*inner(TrialFunction(QD), TestFunction(QD))*dx)
  # Lagrange Multiplier requires -\Delta + I and I on Q 
  p, q = TrialFunction(Q), TestFunction(Q)
  h = CellDiameter(bmesh)
  A = assemble(avg(h) ** (-1) * dot(jump(p), jump(q)) * dS + inner(p, q) * dx)  
  M = assemble(inner(p, q) * dx)  # in A we use DG discretization 

  # For inversion we require parameters for RA 
  params = {'coefs': [1. / mu(0), K(0)], 'pwrs': [-0.5, 0.5], '#[...]'}
  B4 = RA(A, M, parameters=params)
  # define the approximate Riesz map
  B = block_diag_mat([AMG(B0), SOR(B1), HXDiv(B2), SOR(B3), B4])
\end{lstlisting}

\subsection{Mixed-dimensional modeling of signal propagation in neurons} \label{subsec:example_3d1d}

The interaction of slender bodies with its surrounding is of frequent interest in models of blood flow and oxygen transfer \cite{berg_davit_quintard_lorthois_2020, hartung2021mathematical}. It has recently received significant attention as it is a coupling of high dimensional gap (codimension two) which introduces mathematical difficulties \cite{gjerde2020singularity, DAngelo2008, Koppl2018, Koch2020}. Here, we consider an alternative application in neuroscience. In particular, we apply the coupled 3$d$-1$d$ model \cite{laurino2019derivation} to study electric signaling in neurons and its interaction with the extra-cellular matrix. We note that the 
complete model involves a system of partial differential equations (PDE) that represents the electrodiffusion, and a set of ordinary differential equations (ODE) representing the membrane dynamics. Our focus is on the PDE part that arises as part of the operator splitting approach to obtain the solution of the full PDE-ODE problem \cite{Jaeger2021}.

We use the reduced EMI model \cite{buccino2021improving} that represents the extracellular space as a 3$d$ domain and the neuronal body, consisting of soma, axons, and dendrites, as one-dimensional curves. This 3$d$-1$d$ coupled system states to find extracellular and intracellular potentials $(p_3, p_1) $ that satisfy 
\begin{subequations}
	\label{eq:3d1d}
	\begin{align}
			- \nabla \cdot(\sigma_{3} \nabla p_3) + \delta_{\Gamma} \frac{\rho C_m}{\Delta t} ({\Pi^{\rho}_{\Gamma}} p_3 - p_1) & = f_3 & \text{ in } \Omega, \label{eq:3d_1d_1} \\
			- \nabla \cdot( \rho^2 \sigma_{1} \nabla p_1) + \frac{\rho C_m}{\Delta t} (p_1 - {\Pi^{\rho}_{\Gamma}} p_3) & = f_1 & \text{ in } \Gamma.  \label{eq:3d_1d_b}
		\end{align}
\end{subequations}
Here, $\Omega$ is a domain in 3$d$ while $\Gamma$ is the 1$d$ networks of curves, $\Gamma$, represents the neuron by centerlines of soma, axons, and dendrites, see \Cref{fig:brain}. Coupling between the domains is realized by the averaging operator $ {\Pi^{\rho}_{\Gamma}} $ which computes the mean of functions in $\Omega$ on the idealized cylindrical surface that represents the interface between the dendrites and their surroundings. More precisely, given a point $y\in\Gamma$ and $p:\Omega\to\mathbb{R}$, ${\Pi^{\rho}_{\Gamma}}p:\Gamma\to\mathbb{R}$ is such that ${\Pi^{\rho}_{\Gamma}}p(y) = \lvert C_{\rho}(y) \rvert^{-1}\int_{C_{\rho}(y)} u\,\mathrm{d}l$ where $C_{\rho}(y)$ is a circle centered at $y$ with radius $\rho$ in plane whose normal $\nu$ is given by tangent to $\Gamma$ at $y$, cf. \Cref{fig:brain}. That is, $\rho$ represents the radius of a neuron segment and, as such, typically varies in space. However, for simplicity of the presentation, we assume $ \rho $ to be constant. Moreover, by $ \delta_{\Gamma} $ we denote the Dirac measure of $ \Gamma $. The term $ \frac{\rho C_m}{\Delta t} (p_1 - {\Pi^{\rho}_{\Gamma}} p_3) $ represents the electric current flow exchange between the domains across dimensions due to the potential differences with $\Delta t$ being the time step size.
The parameters $\sigma_3$, $\sigma_1$ and $C_m$ represent the extracellular and intracellular conductivity and the membrane capacitance respectively. We also impose boundary conditions to the system \eqref{eq:3d1d} as follows
\begin{subequations}
	\label{eq:3d1d_bc}
	\begin{align}
			p_3 & = g_3 & \text{ on } \partial \Omega_{D}, \\
			- \sigma_{3} \nabla p_3 \cdot \bm{n} & = 0 & \text{ on } \partial \Omega_{N}, \\
			- \rho^2 \sigma_1 \nabla p_1 \cdot \bm{n} & = 0 & \text{ on }  \partial \Gamma,
		\end{align}
\end{subequations}
where $  \partial \Omega_{D} \cup \partial \Omega_{N} = \partial \Omega \backslash \Gamma $ and $ \partial \Omega_{D} \cap \partial \Omega_{N} = \emptyset $.

As in the previous example, we relate a linear system of equations to \eqref{eq:3d1d} that will be used in our software to obtain reliable numerical solutions. Let $ Q_3 \subset H^1(\Omega) $ and $ Q_1 \subset H^1(\Gamma) $ be conforming finite element spaces (e.g. $\mathbb{P}_1$) on the shape-regular triangulation of $ \Omega $ and $ \Gamma $, respectively. Then, a discrete variational formulation of the problem \eqref{eq:3d1d} states to find $ (p_3, p_1) \in Q_3 \times Q_1 $ such that
\begin{equation}
	\label{eq:3d1d_system}
	\underbrace{
	\begin{pmatrix}
		-\sigma_{3} \Delta_\Omega + \tilde{\rho}_t {\Pi^{\rho}_{\Gamma}}' {\Pi^{\rho}_{\Gamma}}  &  - \tilde{\rho}_t {\Pi^{\rho}_{\Gamma}}' \\
		- \tilde{\rho}_t {\Pi^{\rho}_{\Gamma}}	 & -\rho^2 \sigma_1 \Delta_\Gamma + \tilde{\rho}_t I
	\end{pmatrix}
	}_A
	\underbrace{
	\begin{pmatrix}
			p_3 \\ p_1
		\end{pmatrix}
	}_{x}
	=
	\underbrace{
		\begin{pmatrix}
			f_3 \\ f_1 
		\end{pmatrix}
	}_{b},
\end{equation}
with $ \tilde{\rho}_t = \frac{\rho C_m}{\Delta t} $. In \Cref{lst:3d1d_system} we show the implementation of the linear system \eqref{eq:3d1d_system} in FEniCS and cbc.block.
\begin{lstlisting}[style=mystylepython, caption={Implementation of the 3$d$-1$d$ coupled system \eqref{eq:3d1d_system}. Complete code can be
      found in script \texttt{HAZniCS-examples/demo\_3d1d.py}}, captionpos=b, label={lst:3d1d_system}]
  cylinder = Circle(radius=rho, degree=10)
  Rp3, Rq3 = Average(p3, Gamma, cylinder), Average(q3, Gamma, cylinder)
  
  a = block_form(W, 2)
  # Second-order operators
  a[0][0] = sigma3 * (inner(grad(p3), grad(q3)) * dx + inner(p3, q3)) * dx
  a[1][1] = sigma1 * (inner(grad(p1), grad(q1)) * dx + inner(p1, q1)) * dx
  # Metric term
  m = block_form(W, 2)
  m[0][0] = inner(Rp3, Rq3) * dx_
  m[0][1] = -inner(p1, Rv3) * dx_
  m[1][0] = -inner(q1, Ru3) * dx_
  m[1][1] = inner(p1, q1) * dx_
  # Sources
  L = block_form(W, 1)
  L[0] = inner(f3, q3) * dx
  L[1] = inner(f1, q1) * dx
  # Assemble
  AD, M, b = map(ii_assemble, (a, m, L))
\end{lstlisting}
The operator $ A $ is symmetric positive definite and we can use the CG method to solve the system \eqref{eq:3d1d_system}. If we decompose the system as
\begin{equation}
	\label{eq:3d1d_operator}
	A =
	\underbrace{
	\begin{pmatrix}
		-\sigma_3 \Delta_\Omega  &  \\
		& - \rho^2 \sigma_1 \Delta_\Gamma
	\end{pmatrix}
	}_{A_D}
	+
	\tilde{\rho}_t
	\underbrace{
	\begin{pmatrix}
		{\Pi^{\rho}_{\Gamma}}' \\
		-I
	\end{pmatrix}
	\begin{pmatrix}
		{\Pi^{\rho}_{\Gamma}} & -I
	\end{pmatrix}
	}_{M}
\end{equation}
we can identify that the operator $ M $ induces an $L^2$-based \emph{metric space}
\begin{equation} \label{eq:3d1d_metric_space}
    \mathcal{M}(\Gamma) = \{ (q_3, q_1) \in Q_3 \times Q_1 : \int_{\Gamma} (\Pi_\Gamma q_3 - q_1 )^2 < \infty \}.
\end{equation}
We observe that the bilinear form represented by $ M $ is degenerate. More specifically, we can see that for very large values of the coupling parameter $\tilde{\rho}_t$, the semi-definite coupling part $ M $ dominates, and the system becomes nearly singular. The singular part is related to the kernel of the coupling operator, that is $ \ker(M) = \{ (q_3, q_1) \in Q_3 \times Q_1 : \Pi_\Gamma q_3 - q_1 = 0  \}$ which can be a large subspace of the solution space. Consequently, the condition number of the system grows rapidly with increasing $ \tilde{\rho}_t $, which results in slow convergence of the CG solver, even when using the standard AMG method as the preconditioner as in \Cref{subsec:example_poisson}. We remark that \cite{cerroni2019mathematical} demonstrate that (standard, smoothed aggregation) AMG leads to robust solvers when the coupling is weak ($\tilde{\rho}_t\ll1$).

To ensure uniform convergence of the AMG in the parameter $ \tilde{\rho}_t $, we follow the theory of subspace correction method in \cite{Lee2007} to construct block Schwarz smoothers for the AMG method. The blocks are chosen specifically to obtain the $ \tilde{\rho}_t $-uniformly convergent method. We call the linear systems induced by operators such as \eqref{eq:3d1d_operator} \emph{metric-perturbed} problems. In \Cref{subsec:mamg} we demonstrate how to solve the system \eqref{eq:3d1d_system} with HAZniCS methods based on the AMG with specialized block Schwarz smoothers. In \Cref{subsec:3d1d_results} we showcase some key performance points of the solver.


\section{Implementation} \label{sec:implementation}

\begin{figure}[htbp]
	\includegraphics[width=0.7\textwidth]{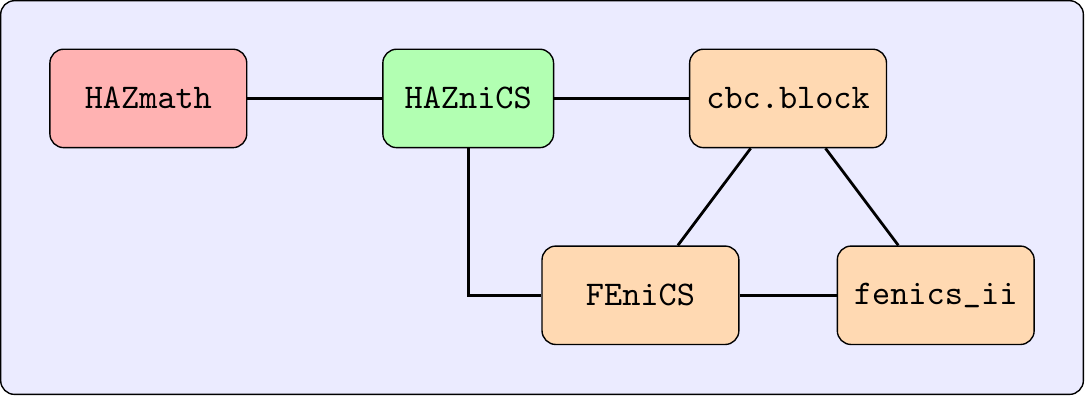}
	\caption{Structure of the HAZniCS framework and relevant components.}
	\label{fig:flowchart}
\end{figure}

The software module HAZniCS combines several libraries, each providing a key functionality for multiphysics simulations. The main components include:
\begin{enumerate}[label=(\roman*)]
	\item HAZmath \cite{hazmath} - a finite element, graph, and solver library built in C;
	\item FEniCS \cite{fenics_book} - a computing platform in Python for solving PDE;
	\item cbc.block \cite{mardal2012block} - an extension to FEniCS that enables assembling and solving block-partitioned problems;          
	\item FEniCS$_{ii}$ \cite{fenicsii} - an extension to FEniCS that enables assembling systems of equations posed on domains
          with different dimensionality (that are not necesarrily embedded manifolds).
\end{enumerate}

We note that while Python and FEniCS use memory management systems, HAZmath requires that the users keep track of the memory themselves. As such, any object transferred between the two systems is copied to make the interactions between FEniCS and HAZmath as simple as possible. Hence, pointers to the underlying data are not passed across the interface, even though this would decrease memory usage. In particular, we then reduce the risk of segmentation fault caused by a pointer in HAZmath that points to some data that Python has deleted. Furthermore, while SWIG provides the means to create Pythonic interfaces to C libraries, e.g., by specifying the input and output of functions, we have decided on making the interface as close as possible to the underlying C code.

In HAZniCS, each of the approximation methods for preconditioners mentioned in the \Cref{sec:examples} is implemented in HAZmath as a C function with the same signature - it takes in a vector (an array of \texttt{double} values), applies a set of operations and returns a solution vector. To be able to use it in Python, the HAZniCS Python library is generated using SWIG \cite{swig} and, in turn, can be imported simply as
\begin{lstlisting}[style=mystylepython, label={lst:import_haznics}]
	import haznics
\end{lstlisting}
In the following code snippets, we demonstrate how this interface is built.

For each preconditioner, HAZniCS stores two functions - a setup and an application function. The setup functions take in different variables depending on the type of the preconditioner but always return a pointer to HAZmath data type \texttt{precond} of a general preconditioner. This data type has two components: data \texttt{data} and matrix-vector operation function \texttt{fct()}, see \Cref{lst:ds_precondtype},
\begin{lstlisting}[style=mystyleC, caption={HAZmath structure for type \texttt{precond}.}, captionpos=b, label={lst:ds_precondtype}]
	typedef struct {
		void *data;
		void (*fct)(REAL *, REAL *, void *);
	} precond;
\end{lstlisting}
where type \texttt{REAL} is a macro of the standard C type \texttt{double}. During setup, HAZmath saves all data necessary for applying the preconditioner and points to the right function that executes the application algorithm. Hence, the matrix-vector function \texttt{fct()} serves as the application function for the preconditioner. It always has the same signature - it takes in two arrays of \texttt{REAL} values (one store's input and the other output vector) and any data related to the matrix-vector operation as \texttt{void*}. 

Now, using the generated HAZniCS Python library, we wrap the HAZmath preconditioner functions as class methods in cbc.block. In this way, efficient HAZmath preconditioners (in C) can be used with FEniCS (or PETSc) operators and cbc.block iterative methods (in Python) in a code that is easily readable and simple to utilize. 

We show an example of the implementation of the AMG preconditioner class in \Cref{lst:amg_block}. Before calling the preconditioner setup function, some input FEniCS data types need to be converted to HAZmath data types. For example, an auxiliary function \texttt{PETSC\_to\_dCSRmat()} converts types \texttt{dolfin.GenericMatrix} or \texttt{dolfin.PETScMatrix} to HAZmath matrix type \texttt{dCSRmat}. This conversion is simple, as PETSc and HAZmath utilize compressed sparse row (CSR) format for matrices where each non-zero element is of double-precision floating-point (double) format. Note that all PETSc-HAZmath conversion functions copy the matrix data rather than copying references to data.
\begin{lstlisting}[style=mystylepython, caption={Preconditioner class \texttt{AMG} implemented in HAZmath backend of cbc.block.}, captionpos=b, label={lst:amg_block}]
	class AMG(Precond):
		def __init__(self, A, parameters=None):
			# change data type for the matrix (to dCSRmat pointer)
			A_ptr = PETSc_to_dCSRmat(A)
			# initialize amg parameters (AMG_param pointer)
			amgparam = haznics.AMG_param()
			# set extra amg parameters
			if parameters:
			    haznics.param_amg_set_dict(parameters, amgparam)
			# set AMG preconditioner
			precond = haznics.create_precond_amg(A_ptr, amgparam)
			#[...]
			Precond.__init__(self, A, "AMG", parameters, amgparam, precond)
\end{lstlisting}
We comment that in HAZmath, all preconditioner application functions have the same signature. On the other hand, in any cbc.block iterative method, all preconditioners are applied through a matrix-vector product method \texttt{matvec()}. Therefore, we define a base  class \texttt{Precond} equipped with a \texttt{matvec()} method designed specifically to call of HAZmath preconditioners. The class is derived from cbc.block data \texttt{block\_base}, making the HAZmath preconditioners fully integrated with other classes and methods of the cbc.block library.
\begin{lstlisting}[style=mystylepython, caption={Baseclass \texttt{Precond} with \texttt{matvec()} implemented in HAZmath backend of cbc.block.}, captionpos=b, label={lst:ds_precondmatvec}]
class Precond(block_base):
  #[...]	
  def matvec(self, b):
    #[...]
    # create solution vector 
    x = self.A.create_vec(dim=1)
    x = df.Vector(df.MPI.comm_self, x.size())		
    #[...]
    # convert rhs and lhs to numpy arrays
    b_np = b[:]
    x_np = x[:]
    # apply the preconditioner (solution saved in x_np)
    haznics.apply_precond(b_np, x_np, self.precond)
    # convert x_np to GenericVector
    x.set_local(x_np)
    return x
\end{lstlisting}

Bridging matrix-vector operation functions of cbc.block (in Python) and HAZmath (in C) is also done using SWIG \cite{swig}. In the file \texttt{haznics.i}, we make a typemap for function \texttt{apply\_precond()} that, before applying the preconditioner matrix-vector function, casts \texttt{numpy} arrays to C arrays of doubles with an additional integer variable indicating array length. This is demonstrated in \Cref{lst:ds_swig}.
\begin{lstlisting}[style=mystyleC, caption={SWIG interface of HAZmath function \texttt{apply\_precond()} that takes in \texttt{numpy} arrays.}, captionpos=b, label={lst:ds_swig}]
	%include "numpy.ii" 
	%include "hazmath.h"
	%numpy_typemaps(double, NPY_DOUBLE, REAL)
	//[...]
	%apply (int DIM1, double* IN_ARRAY1) {(int len1, double* vec1),
		(int len2, double* vec2)}
	//[...]
	%inline %{
		void my_apply_precond(int len1, double* vec1, int len2, double* vec2, 
		precond* pc) {
			//[...]
			apply_precond(vec1, vec2, pc);
		}
		%}
\end{lstlisting}
Finally, within HAZmath, the \verb|apply_precond()| function performs the application of the preconditioner using the data and matrix-vector function that are passed in the input variable \texttt{precond *pc}, see \Cref{lst:ds_haz_applyprecond}.
\begin{lstlisting}[style=mystyleC, caption={HAZmath function \texttt{apply\_precond()}.}, captionpos=b, label={lst:ds_haz_applyprecond}]
	void apply_precond(REAL *r, REAL *z, precond *pc) {
		pc->fct(r, z, pc->data);
	}
\end{lstlisting}

In the following, we detail the implementation of the preconditioners used in \Cref{sec:examples}. In \Cref{subsec:amg} and \Cref{subsec:ra} we show the AMG and the rational approximation preconditioners that can be used through cbc.block extension. On the other hand, what we have described previously is only one of the ways we can use HAZmath solvers and preconditioners with FEniCS. We can also directly call HAZmath functions within FEniCS without relying on cbc.block since we already have a compiled Python library HAZniCS. This use case is shown in \Cref{subsec:mamg} where we describe solvers for metric-perturbed problems.

\subsection{Algebraic multigrid method} \label{subsec:amg}


As the main preconditioning routine, we use the Algebraic MultiGrid Method (AMG)~\cite{Brandt.A;McCormick.S;Ruge.J.1982a}, which constructs a multilevel hierarchy of vector spaces, each of which is responsible for correcting different components of the error. More specifically, our approach is based on the \emph{Unsmoothed Aggregation} (UA-AMG) and the \emph{Smoothed Aggregation} (SA-AMG) method. The UA-AMG method was proposed in~\cite{Vakhutinsky.I;Dudkin.L;Ryvkin.A.1979a} and further developed in~\cite{Blaheta.R.1986a,Marek.I.1991a}. Some popular UA-AMG algorithms are based on graph matching (or pairwise aggregation). Such algorithms with different level of sophistication are found in several works~\cite{Vanek.P;Mandel.J;Brezina.M.1996a,Notay.Y.2010b,Livne.O;Brandt.A.2012a,DAmbra.P;Vassilevski.P.2014a, Kim.H;Xu.J;Zikatanov.L.2003b,HuLinZikatanov2019,HuWuZikatanov2020,Urschel2015}. The SA-AMG method was first proposed in~\cite{Mika.S;Vanek.P.1992a,Mika.S;Vanek.P.1992b} and later extended and analyzed in~\cite{1996Vanek_Mandel_Brezina,Vanek.P;Mandel.J;Brezina.M.1998a,HuVassilevskiXu2016a}.

Compared to classical AMG, one advantage of the aggregation-based AMG methods is that several approximations of near kernel components of the matrix describing the linear system can be preserved as elements of every subspace in the hierarchy. We now briefly explain the basic constructions involved in obtaining multilevel hierarchies of spaces via aggregation, which is known as the \emph{setup phase} of an AMG algorithm.

For a linear system with symmetric and positive definite matrix $\mathsf{A} \in \RR^{n\times n}$, we introduce the undirected graph $\mathcal{G}(\mathsf{A})$ associated with the sparsity pattern of $\mathsf{A}$. The vertices of $\mathcal{G}(\mathsf{A})$  are labeled as $\{1, \ldots, n\}$ and for the set of edges $ \mathcal{E} $ we have $ (i,j) \in \mathcal{E} \iff a_{ij} \neq 0 $. A typical aggregation method consists of four steps stated in the \Cref{alg:setup_amg}. The near kernel components needed in algorithm~\Cref{alg:setup_amg} are often known from the differential operator in hand. When solving a discretized elliptic equation, usually only one near kernel component (the constant function/vector) is used, while for linear elasticity the rigid body modes are utilized in the setup phase. 
\begin{algorithm} \caption{Setup phase of the two-level aggregation-based AMG method.} \label{alg:setup_amg}
	\begin{algorithmic}[1]
		\STATE \textbf{Filter values}: Set $\mathcal{G}(\mathsf{A}) \coloneq \mathcal{G}(\widetilde{\mathsf{A}})$, where $\widetilde{\mathsf{A}}$ is the matrix obtained from $\mathsf{A}$ obtained after filtering out all entries of $\mathsf{A}$ for which $\frac{|a_{ij}|}{\sqrt{a_{ii}a_{jj}}}$ is smaller than a given threshold.
		\STATE \textbf{Create aggregates}: Split the set of vertices $ \{1, 2, \dots, n \} $ as a union of $n_c$ non-overlapping subsets $ \{ \mathfrak{a}_i \}_{i = 1}^{n_c}$.
		\STATE \textbf{Construct coarse space}: Let $ \bm{1}_{\mathfrak{a}_i} $ be the indicator vectors of the aggregates $ \mathfrak{a}_i$, $ i = 1, \dots, n_c $. For given $ k $ near kernel components $ [\psi_1, \ldots, \psi_k] \in \RR^{ n \times k}$, define vectors $\phi_i = [\diag (\psi_1) \bm{1}_{\mathfrak{a}_i}, \ldots, \diag(\psi_k) \bm{1}_{\mathfrak{a}_i}]$ for each aggregate $i = 1, \dots, n_c$. With that, define the coarse space $ V_c \subset V = \RR^n $ of $\dim V_c = k n_c$ as the span of the columns of the matrix $ \mathsf{P} = [\phi_1,\ldots,\phi_{n_c}]\in \RR^{n \times (k n_c)} $.
		\STATE \textbf{Construct course level matrix}: Compute $ \mathsf{A}_c = \mathsf{P}^T \mathsf{A} \mathsf{P} $. 
	\end{algorithmic}
\end{algorithm}
For the multilevel methods, we can repeat the steps in the \Cref{alg:setup_amg} recursively by applying it to $\mathsf{A}_c$ in place of $\mathsf{A}$. The recursive process is halted if the maximum number of levels is reached or the dimension $ n_c $ is smaller than a minimal coarse space dimension.

Furthermore, the SA-AMG algorithm adds a \emph{smoother} to the definition of $ \mathsf{P} $, i.e. in Step 3 of \Cref{alg:setup_amg} we have $ \mathsf{P} = p( \mathsf{S} \mathsf{A} )[\phi_1, \ldots, \phi_{n_c}] \in \RR^{n \times (kn_c)} $. Here, $ p(\cdot) $ is a fixed degree polynomial. Most often the polynomial is chosen to be $ p(t) = 1 - t $, while in the UA-AMG it is set to $ p(t) = 1 $. We note that the (polynomial) smoothing of the basis vectors improves the stability of the coarse spaces. However, unlike in the UA-AMG, the smoothing necessarily results in a larger number of nonzeroes per row in the coarse grid matrices, while smoothing with higher degree polynomials may lead to an inefficient setup algorithm. Thus, the appropriate choice of the smoother $\mathsf{S}$ and the polynomial $ p(\cdot) $ is essential to the stable and fast convergence of the SA-AMG method.

The applications of the aggregation-based AMG preconditioners are summarized in the \Cref{alg:two-level}. We state only the two-level preconditioning iteration, which utilizes a multiplicative preconditioner $ \mathsf{B} = \operatorname{AMG}(\mathsf{A}) \approx \mathsf{A}^{-1}$. 
\begin{algorithm} \caption{Two-level AMG algorithm.} \label{alg:two-level}
	\begin{algorithmic}[1]
		\REQUIRE {Given $ \mathsf{g} \in \RR^n$, do}
		\STATE {Pre-smoothing: $ \mathsf{v} =\mathsf{S} \mathsf{g} $.}
		\STATE {Coarse grid correction: $ \mathsf{w} = \mathsf{v} + \mathsf{P} \mathsf{B}_{c} \mathsf{P}^T(\mathsf{g} - \mathsf{A} \mathsf{v})$.}
		\STATE {Post-smoothing: $\mathsf{B} \mathsf{g} = \mathsf{w} + \mathsf{S}^{T}(\mathsf{g} - \mathsf{A} \mathsf{w})$.}
	\end{algorithmic}
\end{algorithm}

The action of $ \mathsf{B}_{c} $ is determined by the coarse space solver, which can be a direct or another iterative method. In multilevel setting, $\mathsf{B}_c$ represents the recursive application of the \Cref{alg:two-level} where we replace the fine-level matrix $\mathsf{A}$ with the coarse-level matrix $\mathsf{A}_c$. The recursion stops when reaching the maximal coarsest level. This multilevel algorithm is called the V-cycle, but more sophisticated cycling procedures can often be employed.
In HAZmath, other cycles can be used, such as the linear Algebraic Multilevel Iteration (AMLI) methods~\cite{AV_89,AV_90,Johannesbook} and the nonlinear AMLI methods~\cite{AV_91,AV_94,K_02,Panayotsbook,Notay.Y.2010b,HuVassilevskiXu2013} which correspond to optimized polynomial accelerations. 

The implementation of \Cref{alg:setup_amg} and \Cref{alg:two-level} can be found in HAZmath in files\\ \texttt{src/solver/amg\_setup\_ua.c} and \texttt{src/solver/mgcycle.c}, respectively. Due to the extensive length, we skip the implementation code in this paper, but rather show the interface of the HAZmath's AMG method in HAZniCS. 

\begin{lstlisting}[style=mystylepython, caption={Call of the AMG preconditioner for the linear elliptic problem. Complete code can be
      found in script \texttt{HAZniCS-examples/demo\_elliptic\_test.py}}, captionpos=b, label={lst:amg_call_haznics}]
  from haznics import AMG
  # AMG setup parameters
  params = {
    	"AMG_type": haznics.UA_AMG,             
    	"cycle_type": haznics.NL_AMLI_CYCLE,     
    	"smoother": haznics.SMOOTHER_GS,   
    	"coarse_solver": haznics.DIRECT,
    	"aggregation_type": haznics.VMB,
    	"strong_coupled": 0.0, 
    	"max_aggregation": 100,
  }
  # Solver setup
  B = AMG(A, params)
  Ainv = ConjGrad(A, precond=B, tolerance=1e-10)
  # Solve
  x = Ainv * b
\end{lstlisting}

In \Cref{lst:amg_call_haznics} we showcase how to use the AMG method from HAZmath as the preconditioner in FEniCS-related examples. We import the preconditioner class \texttt{AMG} from cbc.block implementation of which has been shown in \Cref{lst:amg_block}. It takes the coefficient matrix $\mathsf{A}$ and an optional dictionary of setup parameters. Such parameters are set through HAZmath macros and are integrated within the HAZniCS Python library through a dictionary. For example, we can specify the type of the AMG method we will apply using the keyword \verb|"AMG_type"| and the value \verb|haznics.UA_AMG|. Other listed keywords determine
\verb|"cycle_type"| (cycling algorithm),
\verb|"smoother"| (type of smoother),   
\verb|"coarse_solver"| (coarse grid solver),
\verb|"aggregation_type"| (type of aggregation), 
\verb|"strong_coupled"| (the filtering threshold in Step 1 of \Cref{alg:setup_amg}) and
\verb|"max_aggregation"| (maximum number of vertices in an aggregate).
Full list of parameters is found in the structure \verb|AMG_param| in the HAZmath's \verb|include/params.h| and values of different macros are given in HAZmath's \verb|include/macro.h|. 

\begin{lstlisting}[style=mystyleC, caption={Implementation of the AMG preconditioner in HAZmath.}, captionpos=b, label={lst:amg_precond_hazmath}]
	void precond_amg(REAL *r, REAL *z, void *data) {
		precond_data *pcdata=(precond_data *)data; // data for the preconditioner
		const INT m = pcdata->mgl_data[0].A.row; // general size of the system
		const INT maxit = pcdata->maxit; // how many times to apply AMG
		INT i;
		
		AMG_param amgparam; param_amg_init(&amgparam);
		param_prec_to_amg(&amgparam, pcdata); // set up AMG parameters
		
		AMG_data *mgl = pcdata->mgl_data; // data for the AMG
		mgl->b.row = m; array_cp(m, r, mgl->b.val); // residual is the rhs
		mgl->x.row = m; dvec_set(m, &mgl->x, 0.0);
		
		for (i = 0; i < maxit; ++i) mgcycle(mgl, &amgparam); // apply AMG
		
		array_cp(m, mgl->x.val, z); // copy the result to z
	}
\end{lstlisting}

Furthermore, the AMG preconditioner is passed to the CG iterative solver \texttt{ConjGrad} from cbc.block to act on the residual in each iteration. As shown in the previous section, the application of the preconditioner is made as a matrix-vector operation, which in the case of the AMG method corresponds to the function \texttt{precond\_amg()} stated in \Cref{lst:amg_precond_hazmath}. It is a simple function that reads the AMG setup data through the variable \texttt{pcdata->mgl\_data}, sets up the right-hand side vector (the residual variable \texttt{r} of the outer iterative method) and initializes the solution vector, applies the AMG algorithm from \Cref{alg:two-level} through the function \texttt{mgcycle()} and returns the computed solution through the increment variable \texttt{z}.

Using HAZmath's implementation of the AMG method through the function \texttt{mgcycle()} gives the flexibility to apply and modify the algorithm to other relevant methods and applications. The following two sections present how we use it in algorithms that approximate inverses of fractional and metric-perturbed operators.

\subsection{Rational approximation} \label{subsec:ra}

In \Cref{subsec:example_ds}, we have introduced a preconditioner based on the (sum of) fractional powers of SPD operators. In particular, in solving the Darcy-Stokes system \eqref{eq:ds_system} iteratively the operator $ B = \left( \mu^{-1} (- \Delta)^{-\frac{1}{2}} + K (- \Delta)^{\frac{1}{2}} \right)^{-1} $ is used in the preconditioner. That means that in each iteration, we need to compute $ z = B r $. We discuss in this section how to use and implement rational approximation \cite{Hofreither2020} that acts as an application of the \emph{inverse} of a fractional operator $ (\alpha A^s + \beta A^t) $, for $ A $ a symmetric positive definite operator, $ s,t \in [-1, 1]$ and $\alpha, \beta \geq 0 $.  

The basic idea is to find a rational function approximating $ f(x) = (\alpha x^{s} + \beta x^{t})^{-1} $ for $ x > 0 $, $\alpha, \beta \geq 0 $ and $ s,t \in [-1, 1] $, that is,
\begin{equation}
	(\alpha x^{s} + \beta x^{t})^{-1} \approx R(x) = \frac{P_{k'}(x)}{Q_{k}(x)}, 
\end{equation}
where $P_{k'}$ and $Q_{k}$ are polynomials of degree $k'$ and $k$, respectively. Assuming $k' \leq k$, the rational function can be given in partial fraction form
\begin{equation}
	R(x) = c_0 + \sum_{i=1}^{n_p} \frac{c_i}{x - p_i},
\end{equation}
for $ c_0 \in \RR $, $ c_i, p_i \in \CC $, $ i = 1, 2, \dots, n_p$. Let $ A $ be a symmetric positive definite operator. Then, the rational function $ R(\cdot) $ can be used to approximate $ f(A) $ as follows,
\begin{equation}
\label{eq:rational}
	z = f(A) r \approx c_0 r + \sum_{i=1}^{n_p} c_i \left( A - p_i I \right)^{-1}r.
\end{equation}
The overall algorithm is shown in \Cref{alg:rational-approx}.

\begin{algorithm}
	\caption{Compute $ z = f(A) r $ using rational approximation.} \label{alg:rational-approx}
	\begin{algorithmic}[1]
		\STATE Solve for $ w_i $:
		$ \left( A - p_i I \right) w_i = r, \quad i = 1, 2, \dots, n_p. $
		\STATE Compute: $ z = c_0 r + \sum\limits_{i=1}^{n_p} c_i w_i $
	\end{algorithmic}
\end{algorithm}
In our case, the operator $ A $ is a discretization of the Laplacian operator $ -\Delta $, and $ I $ is the discrete operator of the $L^2$ inner product. Therefore, the equations in Step 1 of \Cref{alg:rational-approx} can be viewed as discretizations of the shifted Laplacian problems $ -\Delta \, w_i - p_i \, w_i  = r $. For real non-positive poles, the problem is SPD, so we may define fractions or functions of the operator $ -\Delta - p_i I$.

Let $ \mathsf{A} $ be the stiffness matrix associated with $ -\Delta$ and $\mathsf{M}$ a corresponding mass matrix. Consider the following generalized eigenvalue problem
\begin{equation}\label{eq:ra_gep}
	\mathsf{A} \mathsf{U} = \mathsf{M} \mathsf{U} \Lambda, \quad \mathsf{U}^T \mathsf{M} \mathsf{U} = \mathsf{I} \quad \Longrightarrow\quad \mathsf{U}^T \mathsf{A} \mathsf{U} = \Lambda. 
\end{equation}
For any continuous function $ F(x) $, $ x \in [0, \rho] $ we define
\begin{equation}\label{eq:def-f}
	F(\mathsf{A}) \coloneq \mathsf{M} \mathsf{U} f(\Lambda) \mathsf{U}^T \mathsf{M},
\end{equation}
where $ \rho \coloneq \rho \left(\mathsf{M}^{-1} \mathsf{A} \right) $ is the spectral radius of the matrix $ \mathsf{M}^{-1} \mathsf{A} $. We would like to approximate $ f(\mathsf{A}) = (F(\mathsf{A}))^{-1} $ using the rational approximation $ R(x)$ of $ f(x) = \frac{1}{F(x)}$. We note that, if we have a function $ g(t) = f(\rho t) $ defined on the unit interval $ [0, 1] $ and $ r(t) $ is the best rational approximation to $ g(t) $, then
\begin{equation}
	f(x) \approx R(x) = r \left( \frac{x}{\rho} \right) \approx g \left(\frac{x}{\rho} \right), \quad r \left( \frac{x}{\rho} \right) = c_0 + \sum\limits_{i = 1}^{n_p} \frac{c_i}{\frac{x}{\rho} - p_i}.
\end{equation}
Therefore, if we know $ c_i $ and $ p_i $ for $ g(t) $ on the interval $ [0, 1] $ we immediately get 
\begin{equation} \label{eq:ra_scaled}
	f(x) \approx c_0 + \sum\limits_{i = 1}^{n_p} \frac{\rho c_i}{x - \rho p_i}.
\end{equation}
Then, using \eqref{eq:ra_gep}, \eqref{eq:def-f} and \eqref{eq:ra_scaled}, the rational approximation of $f(\mathsf{A})$ is 
\begin{equation}
	f(\mathsf{A}) \approx c_0 \mathsf{M}^{-1} + \sum\limits_{i = 1}^{n_p} \rho c_i \left(\mathsf{A} - \rho p_i \mathsf{M} \right)^{-1}.
\end{equation}
We remark that $f(\mathsf{A})$ is a dual to nodal mapping.

In summary, to apply the rational approximation, we need to find solvers to apply $\mathsf{M}^{-1}$ and each $\left(\mathsf{A} - \rho p_i \mathsf{M} \right)^{-1}$. If $ p_i \in \RR, \, p_i \leq 0$, we end up solving a series of elliptic problems where multigrid methods are very efficient. As mentioned in \Cref{subsec:amg}, the HAZmath library contains several fast-performing implementations of the AMG method, such as SA-AMG and UA-AMG methods.

Furthermore, many methods compute the coefficients $ c_i $ and $ p_i $, see e.g., an overview in \cite{Hofreither2020}. In HAZmath, we have implemented the Adaptive Antoulas-Anderson (AAA) algorithm proposed in \cite{Nakatsukasa2018}. The AAA method is based on a representation of the rational approximation in barycentric form and greedy selection of the interpolation points. In most cases, this approach leads to $ p_i \leq 0 $. Thus we can use the AMG method to solve each problem in Step 1 of \Cref{alg:rational-approx}. We show in the following how we use the rational approximation and other methods from the HAZmath library to solve the Darcy-Stokes problem in \Cref{subsec:example_ds}.

In the demo examples \texttt{HAZniCS-examples/demo\_darcy\_stokes*.py} we have specified the block problem and the preconditioner using FEniCS extensions FEniCS$_{ii}$ and cbc.block, see also \Cref{lst:ds_system} and \Cref{lst:ds_prec}. For the fractional block in \eqref{eq:ds_precond}, we use the rational approximation from HAZmath.
\begin{lstlisting}[style=mystylepython, caption={Call of the HAZmath rational approximation preconditioner in the demo examples \texttt{HAZniCS-examples/demo\_darcy\_stokes*.py}.}, captionpos=b, label={lst:ds_b4}]
	from block.algebraic.hazmath import RA
	#[...]
	parameters = {'coefs': [1./mu(0), K(0)], 'pwrs': [-0.5, 0.5], [...]}
	B4 = RA(A, M, parameters)
\end{lstlisting}
First, we import the preconditioner class \texttt{RA} representing the rational approximation method from the cbc.block backend designated for HAZmath methods. It takes in two matrices, $\mathsf{A}$ and $\mathsf{M}$, that are the discretizations of $H^1$ and $L^2$ inner products on the solution function space. It also takes in an optional dictionary of parameters that, among others, specify weights $\alpha, \beta$ and fractional powers $s, t$. The call of \texttt{RA} sets up the data from the preconditioner, see \Cref{lst:ds_ra}. That is, it computes:
\begin{itemize}
	\item the coefficients $ c_i, p_i $ with the AAA algorithm based on matrices $\mathsf{A}$ and $\mathsf{M}$ and parameters $\alpha, \beta $ in keyword \texttt{'coefs'} and $s, t$ in keyword \texttt{'pwrs'};
	\item AMG levels for each $ \mathsf{A} - p_i \mathsf{M} $ based on optional additional parameters in the \texttt{parameters} dictionary.
\end{itemize}
These two steps are performed in the function \texttt{create\_precond\_ra()} in HAZmath.
\begin{lstlisting}[style=mystylepython, caption={Class \texttt{RA} implemented in HAZmath backend of cbc.block.}, captionpos=b, label={lst:ds_ra}]
	class RA(Precond):
		def __init__(self, A, M, parameters=None):
			# change data type for the matrices (to dCSRmat pointer)
			A_ptr, M_ptr = map(PETSc_to_dCSRmat, (A, M))
			# initialize amg parameters (AMG_param pointer)
			amgparam = haznics.AMG_param()
			#[...]
			haznics.param_amg_set_dict(parameters, amgparam)
			# get scalings
			scaling_a = 1. / A.norm("linf")
			scaling_m = 1. / df.as_backend_type(M).mat().getDiagonal().min()[1]
			# get coefs and powers
			alpha, beta = parameters['coefs']
			s_power, t_power = parameters['pwrs']
			# set RA preconditioner #
			precond = haznics.create_precond_ra(A_ptr, M_ptr, s_power, t_power,
																					alpha, beta, scaling_a, scaling_m,
																					amgparam)
			# [...]
			Precond.__init__(self, A, "RA", parameters, precond)
\end{lstlisting}
Additionally, we need to compute the upper bound on $ \rho = \rho \left(\mathsf{M}^{-1} \mathsf{A} \right) $. In case of $ \mathbb{P}_1 $ finite elements, we have
\begin{equation}\label{diag}
	\rho\left(\mathsf{M}^{-1} \mathsf{A}\right)
	\le \frac{1}{\lambda_{\min{}}(\mathsf{M})} \| \mathsf{A} \|_{\infty}
	\le \frac{d (d+1) }{\min\{\diag(\mathsf{M})\}} \| \mathsf{A} \|_{\infty},
\end{equation}
where $ d $ is the topological dimension of the problem. Thus, the function \texttt{create\_precond\_ra()} also takes scaling parameters to approximate the spectral radius. We note that, in practice, $\mathsf{M}^{-1}$ scales as (at least) inverse of the discretization parameter $h^{-1}$, so the dimension $d$ is not an important factor in the scalings.

The rational approximation preconditioner is then applied in each iteration through a matrix-vector function, as explained at the beginning of \Cref{sec:implementation}. In the case of the rational approximation preconditioner, the matrix-vector function is the HAZmath function \texttt{precond\_ra\_fenics()} that applies the two steps from \Cref{alg:rational-approx}. In \Cref{lst:ds_haz_precond_ra} we show the implementation snippet of the key parts of the preconditioner algorithm from \Cref{alg:rational-approx}.
\begin{lstlisting}[style=mystyleC, caption={HAZmath function \texttt{precond\_ra\_fenics()}.}, captionpos=b, label={lst:ds_haz_precond_ra}]
	void precond_ra_fenics(REAL *r, REAL *z, void *data) {
		//[...]
		// z = z + residues[0] * M^{-1} r
		if(fabs(residues->val[0]) > 0.) { 
			status = dcsr_pcg(scaled_M, &r_vec, &z_vec, &pc_scaled_M, 1e-6, 100, 1, 0);
		}
		array_ax(n, residues->val[0], z_vec.val);
		// [...]
		for(i = 0; i < npoles; ++i) {
			//[...]
			dvec_set(update.row, &update, 0.0);
			// solve (A - poles[i] * M) update = r
			status = dcsr_pcg(&(mgl[i][0].A), &r_vec, &update, &pc_frac_A,1e-6,100,1,0);
			//[...]
			// z = z + residues[i+1]*update
			array_axpy(n, residues->val[i+1], update.val, z_vec.val);
		}
	}
\end{lstlisting}

\subsection{Solvers for interface metric-perturbed problems} \label{subsec:mamg}
In this section, we continue with presenting the implementation of the solver for the 3$d$-1$d$ coupled problem \eqref{eq:3d1d_system} in \Cref{subsec:example_3d1d}. Additionally, we introduce an alternative way to use HAZmath solvers in FEniCS. In the previous section, we have bridged the two libraries via a class of preconditioners implemented in cbc.block, while here we directly call functions from HAZmath through the generated Python interface.

In the demo example \texttt{demo\_3d1d.py} we have specified the block problem \eqref{eq:3d1d_system} using FEniCS extensions FEniCS$_{ii}$ and cbc.block, see also \Cref{lst:3d1d_system}. Next, we display in \Cref{lst:3d1d_solver} the steps necessary to use HAZmath solver for this block problem directly through the library \texttt{haznics}. 
\begin{lstlisting}[style=mystylepython, caption={Call of HAZniCS solver for the 3$d$-1$d$ coupled system \eqref{eq:3d1d_system}. Complete code can be
      found in script \texttt{HAZniCS-examples/demo\_3d1d.py}}, captionpos=b, label={lst:3d1d_solver}]
  # convert vectors
  bb = ii_convert(b)
  b_np = bb[:]
  bhaz = haznics.create_dvector(b_np)
  xhaz = haznics.dvec_create_p(n)
  # convert matrices; A = AD + rho * M
  Ahaz = block_to_haz(A)
  Mhaz = block_to_haz(M)
  ADhaz = block_to_haz(AD)

  # call solver
  niters = haznics.fenics_metric_amg_solver(Ahaz, bhaz, xhaz, ADhaz, Mhaz)
\end{lstlisting}
The listing consists of three parts: data conversion, a wrapper for the solver function and specifying solver parameters. First, after assembly, the system matrix $ \mathsf{A} $ and the right hand side $ \mathsf{b} $ are of type \texttt{block\_mat} and \texttt{block\_vec}, respectively. We convert them to HAZmath data types \texttt{dvector} and \texttt{block\_dCSRmat} so we are able to use them in the solver that is called through the HAZniCS function \texttt{fenics\_metric\_amg\_solver()}. This auxiliary function acts as an intermediary to set solver data and parameters and to run the solver. An excerpt from the function is given in \Cref{lst:wrapper}. We remark that the signature of the wrapper function needs to be added to the interface file \texttt{haznics.i} to be able to use it through the HAZniCS Python library since it is not a part of the standard HAZmath library.
\begin{lstlisting}[style=mystyleC, caption={Wrapper function for the solver of the system \eqref{eq:3d1d_system}.}, captionpos=b, label={lst:wrapper}]
	INT fenics_metric_amg_solver(block_dCSRmat *A, dvector *b, dvector *x, 
															 block_dCSRmat *AD, block_dCSRmat *M)
	{
		/* set Parameters from Reading in Input File */
		input_param inparam;
		param_input_init(&inparam);
		param_input("./input_metric.dat", &inparam);
		
		//[...]
		/* Use Krylov Iterative Solver */
		if ( (linear_itparam.linear_precond_type >= 10) && \
				 (linear_itparam.linear_precond_type < 15) ){
			solver_flag = linear_solver_bdcsr_krylov_metric_amg(A, b, x,&linear_itparam, 
																													&amgparam, AD, M);
		}
		/* No preconditioner */
		else{
			solver_flag = linear_solver_bdcsr_krylov(A, b, x, &linear_itparam);
		}
		return solver_flag;
	}
\end{lstlisting}

Unless we want to use default values, it is required to set relevant parameters for the HAZmath solver, such as the tolerance of the iterative method or the type of the preconditioner. This can be done by creating an input file that passes the specific parameters to HAZmath to a variable of type \texttt{input\_param}. A snippet of the input file \texttt{input\_metric.dat} for the 3$d$-1$d$ coupled problem can be found in \Cref{lst:input_file}. 


\begin{lstlisting}[style=mystyleC, caption={Input file example. Complete code can be found in script \texttt{HAZniCS-examples/input\_metric.dat}.}, captionpos=b, label={lst:input_file}]
  //[...]
  linear_itsolver_type    = 1     % 1: CG
  linear_itsolver_maxit   = 1000
  linear_itsolver_tol     = 1e-8
  linear_stop_type        = 1     % 1: ||r||/||b||
  linear_precond_type     = 14    % 14: Schwarz on interface part (symm multipl)
                                  %    + AMG on the whole matrix 
  //[...]
  AMG_type        = SA
  AMG_cycle_type  = V 
  AMG_levels      = 10
  AMG_maxit       = 1
\end{lstlisting}
Finally, this parameter setup allows to apply the solver \texttt{linear\_solver\_bdcsr\_krylov\_metric\_amg()}. We recall that we have chosen to solve the 3$d$-1$d$ problem \eqref{eq:3d1d_system} by the CG method preconditioned with AMG that uses block Schwarz smoothers to obtain robustness in the coupling parameter $ \tilde{\rho}_t \gg 1 $. The HAZmath implementation of that solver has a slight modification that uses a combination of the block Schwarz and Gauss-Seidel smoothers. We give a few details on the algorithm and its implementation in the following section.

\subsubsection{Metric-perturbed algebraic multigrid method}

Let us go back to the operator \eqref{eq:3d1d_operator} and set $ V = Q_3 \times Q_1 $. The general subspace correction method looks for a stable space decomposition 
\begin{equation} \label{eq:space-decomposition}
	V = V_0 + V_1 + \ldots + V_J 
\end{equation}
to divide solving the system on the whole space $ V $ to solving smaller problems on each subspace and summing up the contributions in additive or multiplicative fashion \cite{Xu1992SIAMReview,XuZikatanov2002}. Furthermore, the following condition from \cite{Lee2007} is sufficient to obtain a robust subspace correction method to solve nearly singular system such as \eqref{eq:3d1d_system}:
\begin{equation} \label{eq:kernel-decomposition}
	\operatorname{Ker}(M)\cap V =
	(\operatorname{Ker}(M)\cap V_0)+                           
	(\operatorname{Ker}(M)\cap V_1)+\ldots+                           
	(\operatorname{Ker}(M)\cap V_J).
\end{equation}
More specifically, we can employ this space decomposition to create a robust AMG method where $ V_0 $ represents the coarse space and $ V_i, \; i = 1, \dots, J$, define a Schwarz-type smoother. By robustness, we imply that the convergence of the method is independent of the values of the coupling parameter $\tilde{\rho}_t$ and mesh size parameter $h$. To construct subspace splitting satisfying~\eqref{eq:kernel-decomposition} it is necessary to choose the subspaces so that the following holds: For each element of a frame spanning the null-space of $M$, there exists a subspace containing this frame element. Notice that this is a requirement that does not assume that the frame element is known, but rather, the assumption is that a subspace where this element is contained is known. 



\begin{algorithm}
	\caption{Compute $ z = B r $ using metric-perturbed AMG} \label{alg:metric_amg}
	\begin{algorithmic}[1]
		\REQUIRE Given $r$ and $z \gets 0$
		\STATE Solve on the interface using forward Schwarz smoother: $ z \gets  z +  {\Pi^{\rho}_{\Gamma}} B_{\text{Schwarz}} \Pi'_{\Gamma} r$ 
		
		\STATE Solve on the whole space using AMG method:  $z \gets z + B_{\text{AMG}}(r - A z)$
		
		\STATE Solve on the interface using backward Schwarz smoother: $z \gets z +  {\Pi^{\rho}_{\Gamma}} B'_{\text{Schwarz}} \Pi'_{\Gamma} (r- A z) $ 

	\end{algorithmic}
\end{algorithm}

From Algorithm~\ref{alg:metric_amg}, we can see that $B$ is defined as
\begin{equation*}
	I - BA : = (I - {\Pi^{\rho}_{\Gamma}} B'_{\text{Schwarz}} \Pi'_{\Gamma}A) (I - B_{\text{AMG}} A) (I - {\Pi^{\rho}_{\Gamma}} B_{\text{Schwarz}} \Pi'_{\Gamma}A).
\end{equation*}
It is easy to see that $B$ is symmetric and, following the theory developed in~\cite{hu2013combined}, $B$ is also positive definite if $B_{\text{AMG}}$ is symmetric positive definite and $B_{\text{Schwarz}}$ is nonexpansive. Therefore, it can be used as a preconditioner for the CG method.  This preconditioner is implemented in HAZmath and its excerpt from the function is given in \Cref{lst:metric-AMG}.

\begin{lstlisting}[style=mystyleC, caption={metric AMG preconditioner.}, captionpos=b, label={lst:metric-AMG}]
  void precond_bdcsr_metric_amg_symmetric(REAL *r, REAL *z, void *data)	
  {
    //[...] 
    // Schwarz method on the interface part
    Schwarz_param *schwarz_param = predata->schwarz_param;
    Schwarz_data *schwarz_data = predata->schwarz_data;
    smoother_dcsr_Schwarz_forward(schwarz_data, schwarz_param, &zz, &rr);		

    //[...]
    // AMG solve on the whole matrix
    AMG_data_bdcsr *mgl = predata->mgl_data;
    mgl->b.row  total_row; array_cp(total_row, r, mgl->b.val); 
    mgl->x.row = total_col; array_cp(total_row, z, mgl->x.val);
    for ( i=maxit; i--; ) mgcycle_bdcsr(mgl,&amgparam);

    //[...]
    // Schwarz method on the interface part
    smoother_dcsr_Schwarz_backward(schwarz_data, schwarz_param, &zz, &rr);

    //[...]
  }
\end{lstlisting}

\section{Results} \label{sec:results}

In this section, we show the performance of the solvers and preconditioners developed for the examples in \Cref{sec:examples}. We recall that their complete code can be found in \cite{github_haznics}.

\subsection{Linear elliptic problem} \label{subsec:poisson_results}

We use the 3$d$ elliptic problem \eqref{eq:poisson} to compare the HAZniCS solvers to already established solver libraries. This way, we demonstrate that HAZniCS, specifically the AMG solver within, shows a fast and reliable performance when solving common PDE problems. For the comparison, we use the AMG method BoomerAMG from the HYPRE library \cite{hypre} of scalable linear solvers and multigrid method that is already integrated within FEniCS software through PETSc. We note that all the computations are performed in serial on a workstation with an 11th Gen Intel(R) Core(TM) i7-1165G7 @ 2.80GHz (8 cores) and 40GB of RAM.

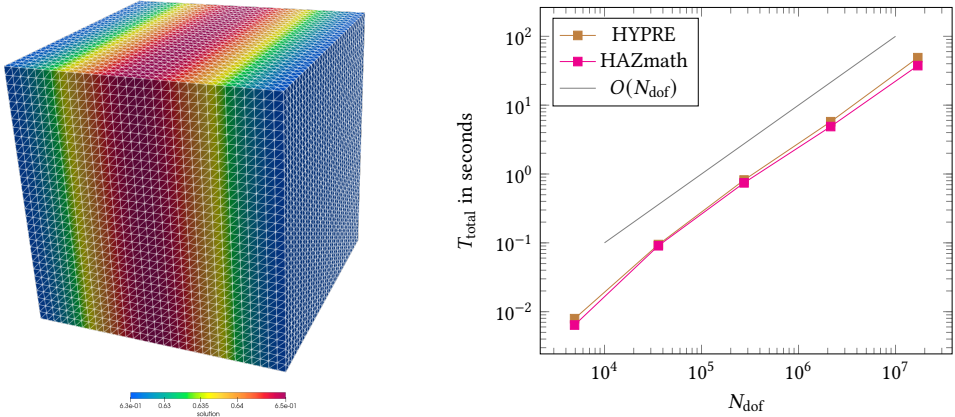
\begin{figure}[htbp]
	\centering
	\includegraphics[width=0.45\textwidth]{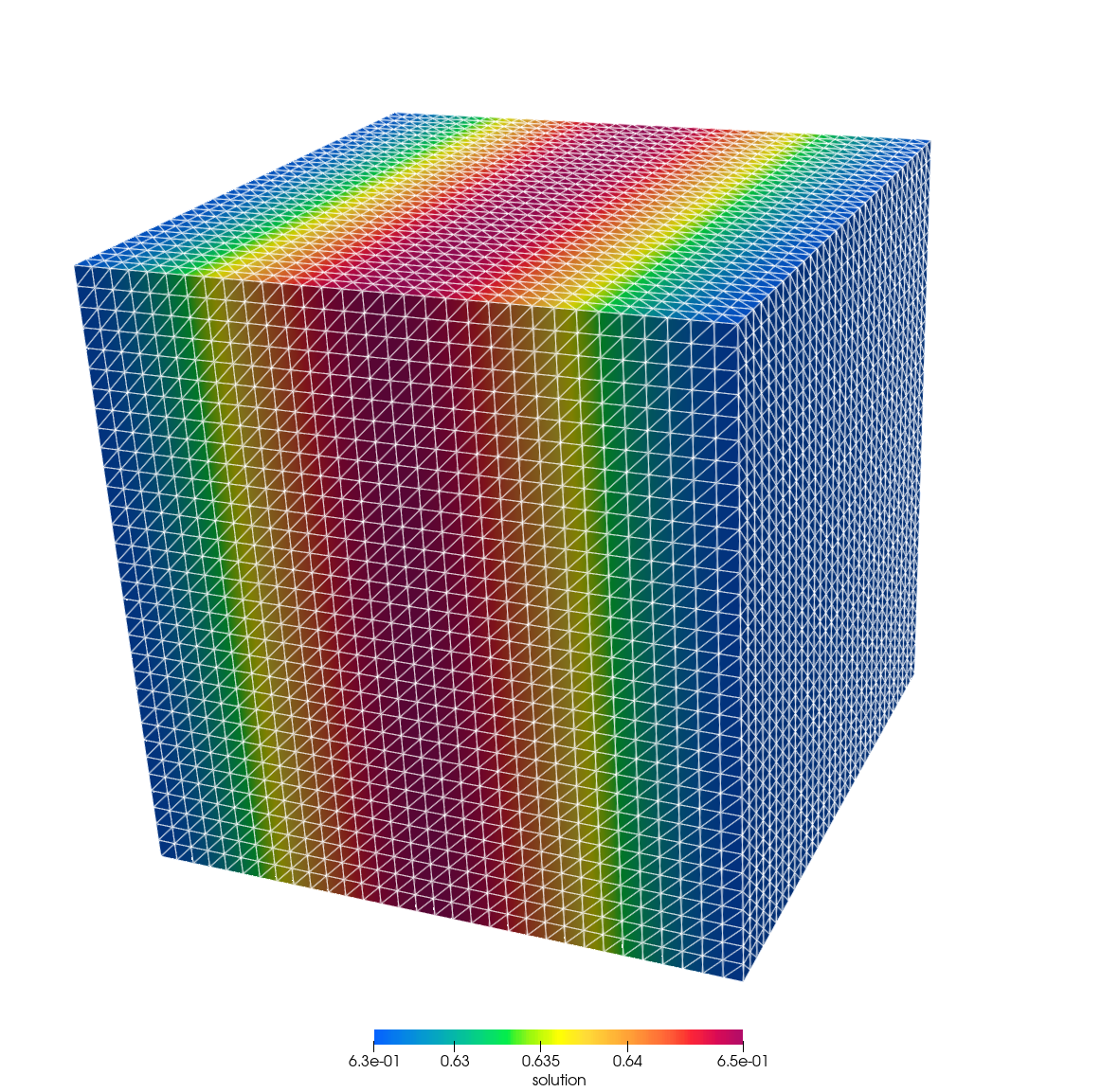}
	\begin{tikzpicture}[thick, scale=0.8]
		\begin{axis}[
			xmode=log,
			ymode=log,
			legend pos=north west,
			xlabel=$N_{\text{dof}}$,
			ylabel=$T_{\text{total}}$ in seconds
			]
			\addplot[brown, mark=square*] table {
				4913 0.007913589477539062
				35937 0.09441876411437988
				274625 0.817718505859375
				2146689 5.748968124389648
				16974593 48.952526807785034
			};
			\addlegendentry{HYPRE}
			\addplot[magenta, mark=square*] table {
				4913 0.006397651626586898
				35937 0.09070765194702157
				274625 0.7436500520935061
				2146689 4.889569152023313
				16974593 37.51694211402889
			};
			\addlegendentry{HAZmath}
			\addplot[gray] table {
				10000 0.1
				100000 1
				1000000 10
				10000000 100
			};
			\addlegendentry{\(O(N_{\text{dof}})\)}
		\end{axis}
	\end{tikzpicture}
	\caption{(Left) Illustration of the domain of the elliptic problem \eqref{eq:poisson} and its solution profile. (Right) Total CPU time required to solve \eqref{eq:poisson} with the CG method up to relative residual tolerance $10^{-6}$. The brown data points represent the total elapsed time of the solver (setup + solve) when CG is preconditioned with the HAZmath AMG method, while the magenta data points represent the total elapsed time when CG is preconditioned with HYPRE AMG. Results are obtained by running \texttt{HAZniCS-examples/demo\_elliptic\_test.py}.
        }
	\label{fig:poisson_cpu}
\end{figure}

The results are given in \Cref{tab:cpu_ndof} and the right part of \Cref{fig:poisson_cpu}. It is clear that the AMG methods of HYPRE and HAZmath show similar performance. While the HYPRE BoomerAMG method gives fewer total CG iterations and consequently less solving time, the setup of the HAZmath's UA-AMG method is multiple times faster while still taking comparable solving time. Therefore, we are confident about using the methods from HAZniCS in our multiphysics solvers, namely the HAZmath's AMG method as a component of the rational approximation and metric-perturbed preconditioners from Sections \ref{subsec:ra} and \ref{subsec:mamg}.

\begin{table}[htpb]
	\begin{tabular}{r||rrrr||rrrr}
		\hline 
		\multicolumn{1}{c||}{} & \multicolumn{4}{c||}{HAZmath} & \multicolumn{4}{c}{HYPRE} \\ \cline{2-9}
		\multicolumn{1}{c||}{${N}_{\text{dof}}$} & \multicolumn{1}{c}{$N_{\text{iter}}$} & \multicolumn{1}{c}{Setup (s)} & \multicolumn{1}{c}{Solve (s)} & \multicolumn{1}{c||}{\textbf{Total (s)}} & \multicolumn{1}{c}{$N_{\text{iter}}$} & \multicolumn{1}{c}{Setup (s)} & \multicolumn{1}{c}{Solve (s)} & \multicolumn{1}{c}{\textbf{Total (s)}} \\ \hline \hline
		\textbf{729} 		&  9 & 0.0004 & 0.0010 & \textbf{0.0014} & 7 & 0.0007 & 0.0007 & \textbf{0.0014} \\
		\textbf{4913} 		& 10 & 0.0012 & 0.0052 & \textbf{0.0064} & 8 & 0.0040 & 0.0039 & \textbf{0.0079} \\
		\textbf{35937} 		& 11 & 0.0077 & 0.0830 & \textbf{0.0907} & 8 & 0.0393 & 0.0551 & \textbf{0.0944} \\
		\textbf{274625} 	& 11 & 0.0631 & 0.6805 & \textbf{0.7436} & 9 & 0.2507 & 0.5670 & \textbf{0.8177} \\
		\textbf{2146689} 	& 12 & 0.5308 & 4.4215 & \textbf{4.9523} & 9 & 1.9855 & 3.8539 & \textbf{5.8394} \\
		\textbf{16974593}   & 12 & 4.6096 & 32.907 & \textbf{37.517} & 9 & 19.255 & 29.697 & \textbf{48.952} \\
		\hline	
	\end{tabular}
	\caption{Performance of the CG method preconditioned with either HAZmath AMG or HYPRE AMG, with regards to the number of degrees of freedom $N_{\text{dof}}$. We measure number of iterations ($N_{\text{iter}}$) and CPU time in seconds for setup and solve part of the solver required to solve the elliptic problem \eqref{eq:poisson} with relative residual tolerance $10^{-6}$. Results are obtained by running \texttt{HAZniCS-examples/demo\_elliptic\_test.py}.}
	\label{tab:cpu_ndof}
\end{table}

\subsection{Darcy-Stokes problem} \label{subsec:ds_results}

%
%
%
To demonstrate the performance of the rational approximation algorithms of HAZniCS, we next focus on the Darcy-Stokes problem \eqref{eq:darcy_stokes} and its preconditioner \eqref{eq:ds_precond}. Implementation of the preconditioner in HAZniCS can be found in \Cref{lst:ds_prec}, and we recall that we utilize multilevel methods for the Stokes velocity and Darcy flux blocks while the multiplier block uses rational approximation detailed in \Cref{subsec:ra}.

Let us first showcase the robustness and scalability of implementing the Darcy-Stokes preconditioner. Here we focus on the (more challenging) case $\Omega_S$, $\Omega_D\subset\mathbb{R}^3$ while results for a similar study in two dimensions are given in \Cref{sec:ds_2d}. Let now $\Omega_S=\left[0, \tfrac{1}{2}\right]\times \left[0, 1\right]^2$ and $\Omega_D=\left[\tfrac{1}{2}, 1\right]\times \left[0, 1\right]^2$. We consider discretization of \eqref{eq:ds_system} by (stabilized, cf. \Cref{lst:ds_system}), $\mathbb{C}\mathbb{R}_1$-$\mathbb{P}_0$ elements in the Stokes domain, $\mathbb{R}\mathbb{T}_0$-$\mathbb{P}_0$ elements in the Darcy domain and $\mathbb{P}_0$ elements on the interface. Using gradually refined meshes of $\Omega_S\cup\Omega_D$, which match on the interface $\Gamma$, the choice of elements leads to linear systems with $2\cdot 10^{3} < N_{\text{dofs}} < 11 \cdot 10^{6}$. Furthermore, we shall vary the model parameters such that $10^{-6}\leq \mu, K \leq 1$ while $D=0.1$ is fixed. 

\begin{figure}[htbp]
	\includegraphics[width=0.8\textwidth]{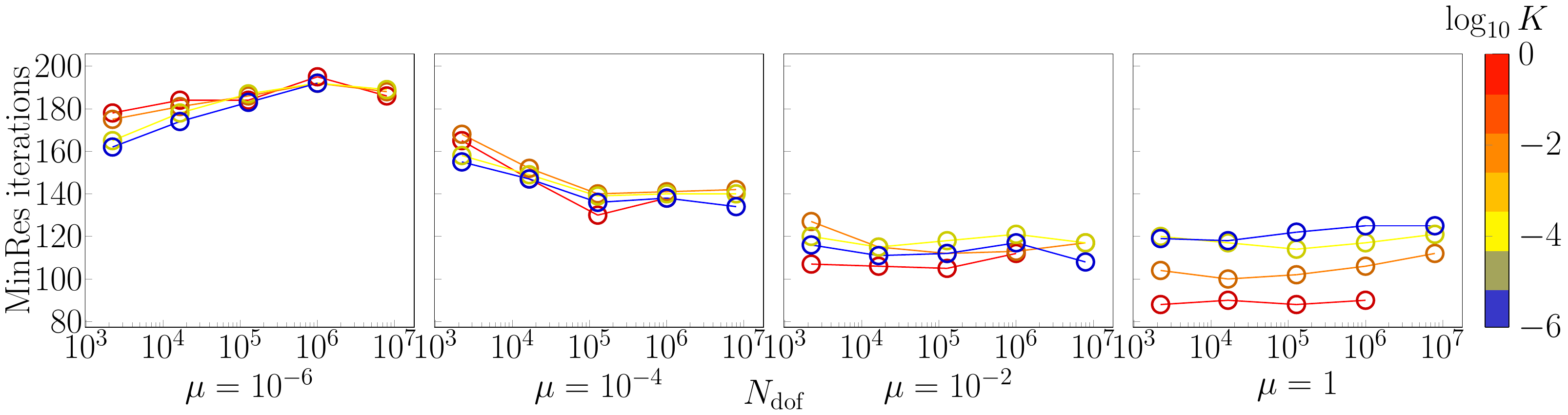}
	\includegraphics[width=0.8\textwidth]{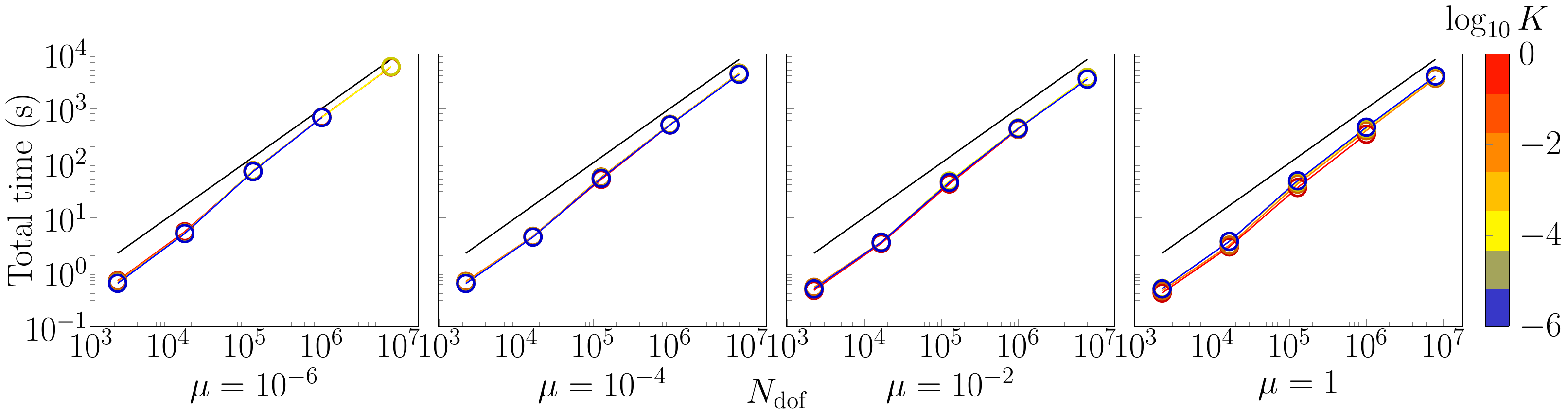}
	\caption{
          Performance of Darcy-Stokes preconditioner \eqref{eq:ds_precond} implemented
          in \Cref{lst:ds_prec}.
Discretization by $(\mathbb{C}\mathbb{R}_1-\mathbb{P}_0)-(\mathbb{R}\mathbb{T}_0-\mathbb{P}_0)-\mathbb{P}_0$
          elements with $D=0.1$.           
          (Top) Number of MinRes iterations until convergence
          for different values of $\mu, K$ and mesh sizes. (Bottom) Total solution
          time for solving \eqref{eq:ds_system} including the setup time of the
          preconditioner. Black line indicates linear scaling.
          Results are obtained by running 
          \texttt{HAZniCS-examples/demo\_darcy\_stokes\_3d\_flat.py}.
        }
	\label{fig:darcy-stokes-precond}
\end{figure}

The performance of the preconditioner is summarized in \Cref{fig:darcy-stokes-precond}, where we list the dependence of the solution time and the number of MinRes iterations on mesh size and model parameters. Here the convergence criterion is the reduction of the preconditioned residual norm by $10^{12}$. Moreover, the tolerance in the rational approximation is set to $10^{-12}$ yielding roughly $n_p\approx 20$ poles in \eqref{eq:rational}. However, numerical experiments \cite{budisa2022rational} suggest that a less accurate approximation, leading to as little as $6$ poles, could be sufficient. In  \Cref{fig:darcy-stokes-precond}, it can be seen that iteration counts are bounded in the parameters, that is, \eqref{eq:ds_precond} defines a parameter robust Darcy-Stokes preconditioner. Moreover, the implementation in \Cref{lst:ds_prec} leads to optimal, $\mathcal{O}(N_{\text{dofs}})$, scaling.

\begin{figure}[htbp]
	\includegraphics[height=0.3\textwidth]{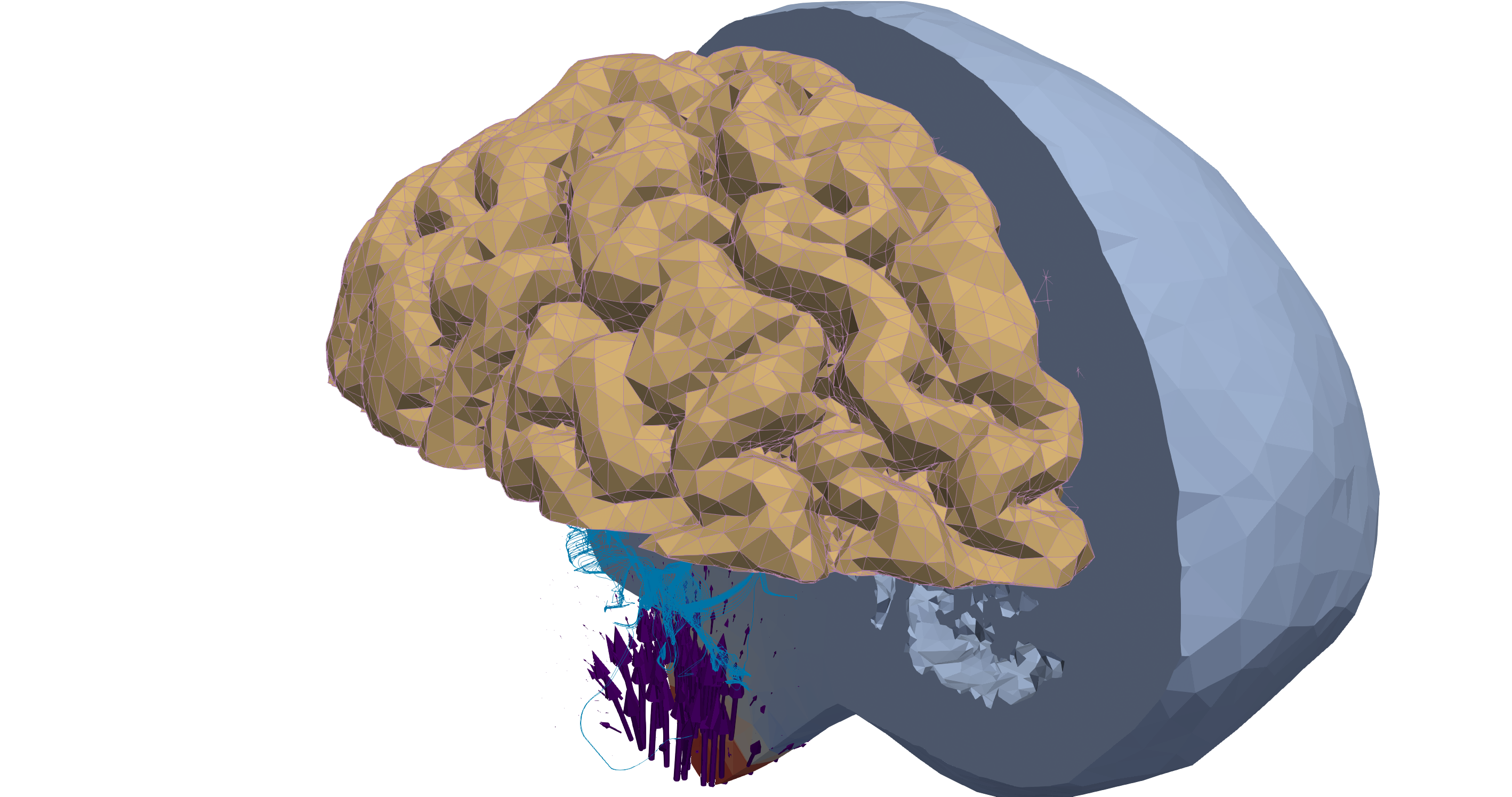}
	\includegraphics[height=0.3\textwidth]{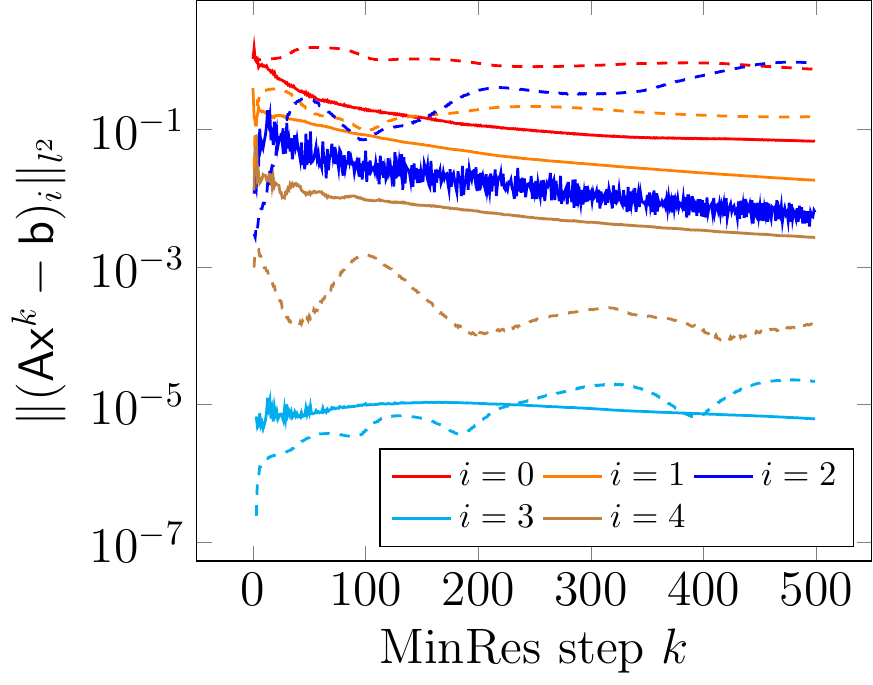}
	\caption{
          Darcy-Stokes model on realistic brain geometry. (Left) Solution
          field when no-slip boundary conditions are considered everywhere except
          on the small region on the base (cf. larger pressure in red). Here, the traction locally increases pressure
          and induces flow in the Stokes domain. Pressure in the Darcy domain is rather
          uniform. Flow in the Darcy domain is visualized by streamlines.
          (Right) Convergence of the solution components (denoted by $0\leq i \leq 4$
          for subspaces $\bm{V}_S$, $Q_S$, $\bm{V}_D$, $Q_D$ and $\Lambda$)
          of \eqref{eq:ds_system} when
          using preconditioner \eqref{eq:ds_precond} (solid lines). Dashed lines
          show (diverging) behavior when using a simpler preconditioner which
          utilized $(-\Delta+I)^{1/2}$ on the interface. 
          Results are obtained by running 
          \texttt{HAZniCS-examples/demo\_darcy\_stokes\_brain.py}.          
        }
	\label{fig:darcy-stokes-results}
\end{figure}

The experimental setup of our previous example led to rather small multiplier spaces with $\text{dim}\Lambda=2048$ for the finest mesh considered. To get larger interfaces and multiplier spaces, we finally turn to the brain geometry in \Cref{fig:brain}. Although realistic, the geometry is still largely simplified as we have excluded the cerebellum, the aqueduct, and the central canal and expanded the subarachnoid space to allow more visible CSF pathways. Nevertheless, the geometry fully represents the complexity of the interface (gyric and sulcal brain surface), which is an important part and an additional difficulty when solving the coupled viscous-porous flow problem. Using the same discretization as before the computational mesh leads to $N_{\text{dofs}}\approx 11\cdot 10^{6}$ with $\text{dim}\Lambda \approx 50\cdot 10^3$ For the purpose of illustration we set $\mu=3$, $K=10^{-4}$, $D=0.5$ and consider most of the outer surface of $\Omega_S$ with no-slip boundary condition except for a small region on the bottom where traction is prescribed. The flow field computed after 500 iterations of MinRes is plotted in \Cref{fig:darcy-stokes-results}. Therein we also compare convergence of MinRes solver using preconditioner \eqref{eq:ds_precond} with a simpler one which uses in the $\Lambda$ block the operator $K(-\Delta+I)^{1/2}$, cf. the analysis in \cite{layton2002coupling}. Importantly, we observe that the new preconditioner, which ignores the intersection structure of the multiplier space, leads to very slow convergence or even divergence of the unknowns. In contrast, with \eqref{eq:ds_precond} MinRes appears to be converging. We remark that the rather slow (in comparison to \Cref{fig:darcy-stokes-precond}) convergence with block diagonal preconditioner \eqref{eq:ds_precond} is related to the thin-shell geometry of the Stokes domain. In particular, the performance of block-diagonal Stokes preconditioner using the mass matrix approximation for the Schur complement is known to deteriorate for certain boundary conditions when the aspect ratio of the domain is large \cite{sogn2022stable}.

%
%


\subsection{3d-1d coupled problem} \label{subsec:3d1d_results}
\begin{figure}[htbp]
	\includegraphics[width=0.45\textwidth]{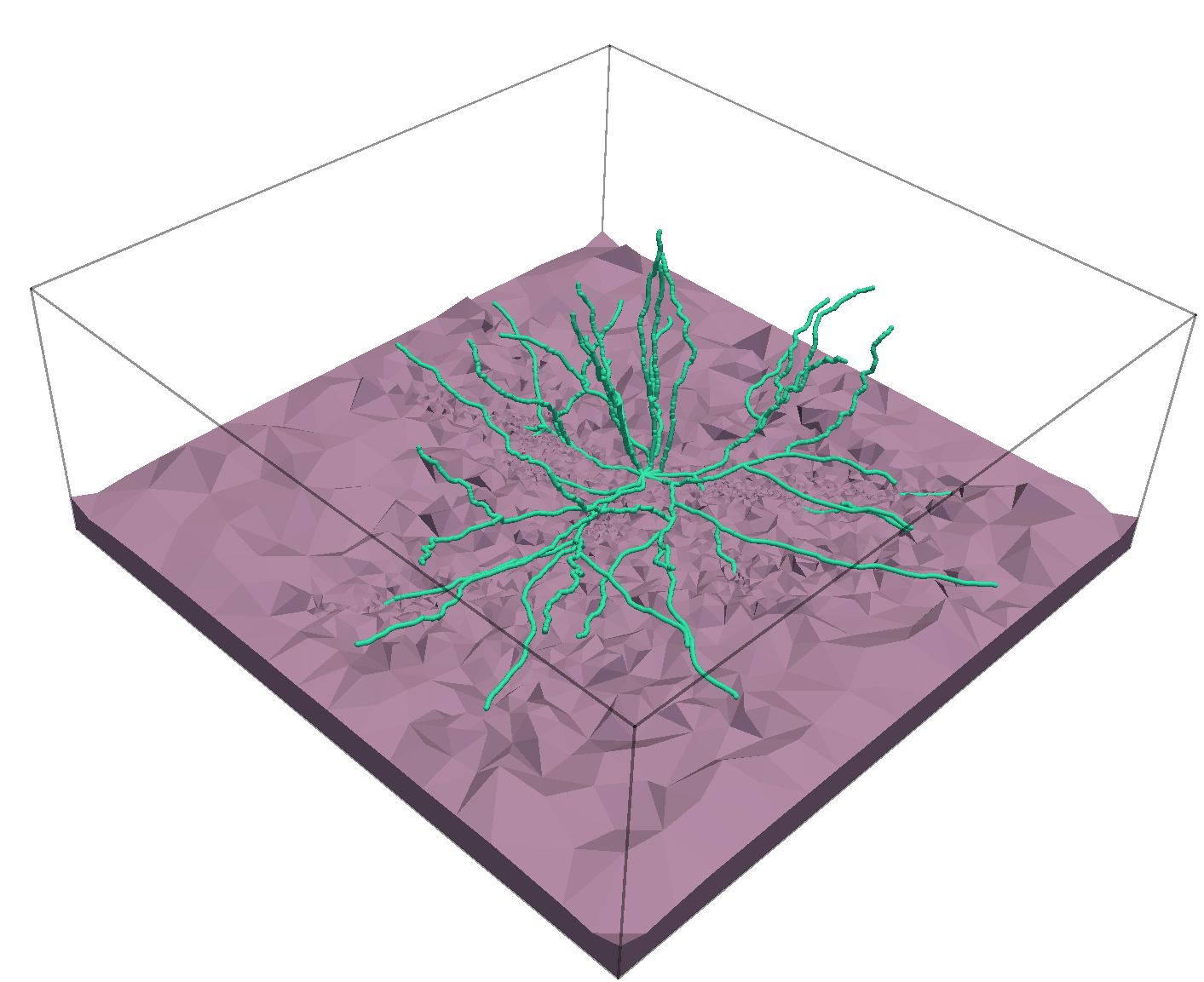}
	\includegraphics[width=0.45\textwidth]{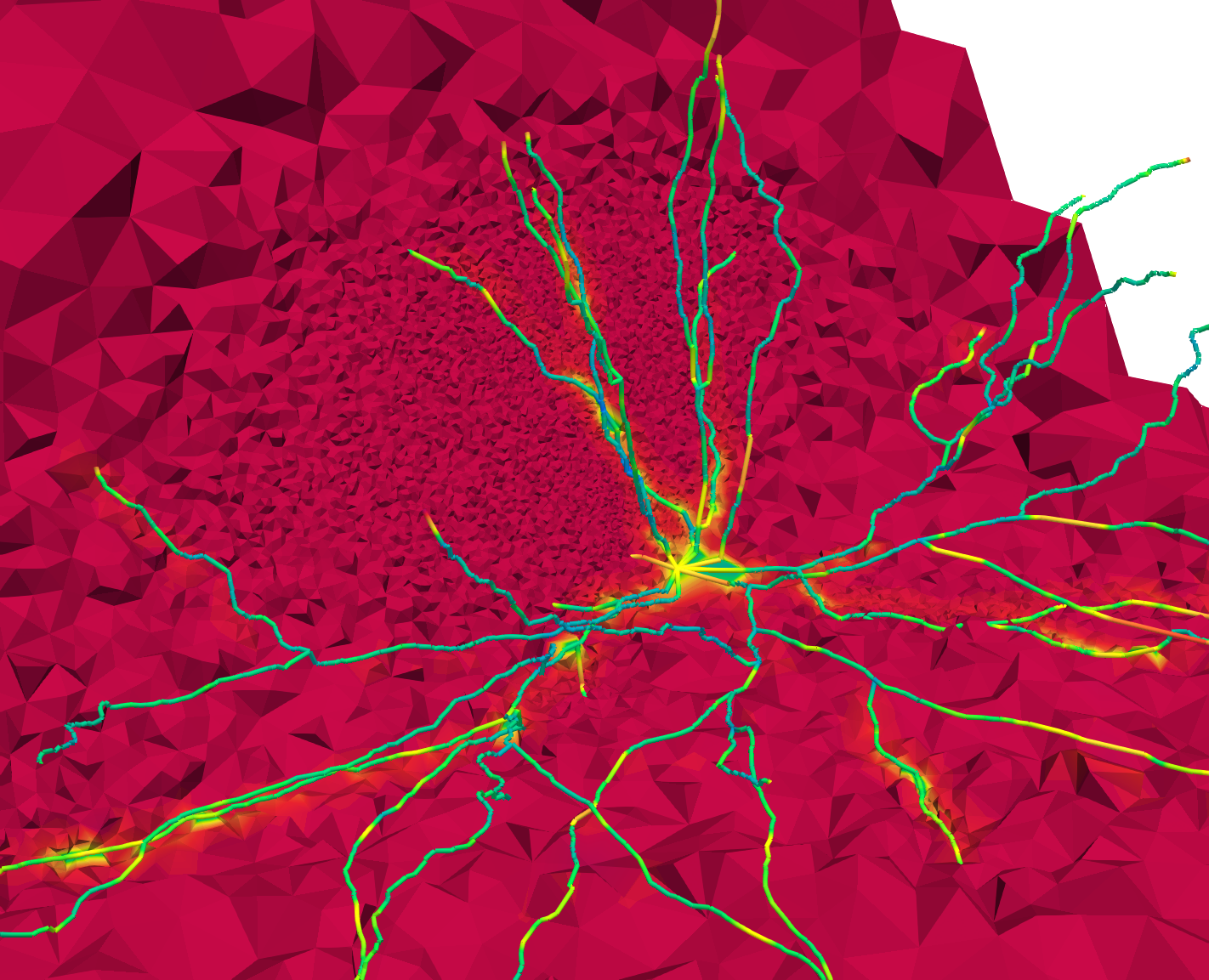}
	\caption{Domain geometry of the 3$d$-1$d$ problem \eqref{eq:3d1d}. (Left) The 1$d$ domain as the neuron and the network of neuronal dendrites is marked in blue and a shallow clip of the 3$d$ brain tissue domain is marked in grey. The outline of the 3$d$ domain is marked with black lines. (Right) A clip of the solution of potentials $(p_3, p_1)$. Results are obtained by running \texttt{HAZniCS-examples/demo\_3d1d.py}.}
	\label{fig:3d1d-geo}
\end{figure}

Lastly, we demonstrate how the mixed-dimensional flow problem from \Cref{subsec:example_3d1d} is solved using the HAZniCS solver for metric-perturbed problems. The problem is defined by the geometry illustrated in \Cref{fig:3d1d-geo}. The neuron geometry is obtained from the NeuroMorpho.Org inventory of digitally reconstructed neurons, and glia \cite{neuron_geo}. The neuron from a rat's brain includes a soma and 72 dendrite branches. It is embedded in a rectangular box of approximate dimensions $281 \mu m \times 281 \mu m \times 106 \mu m $. Then, the mixed-dimensional geometry is discretized with an unstructured tetrahedron in a way conforming to $\Gamma$, i.e., the $1d$ neuron mesh consists of the $3d$ edges lying on $\Gamma$. As discretization, we have $\mathbb{P}_1$ finite elements for both the 3$d$ and 1$d$ function spaces. Overall we end up with 641 788 degrees of freedom in 3$d$ and $3156$ degrees of freedom for the 1$d$ problem. Additionally, we enforce homogeneous Neumann conditions on the outer boundary of both subdomains. 

To obtain the numerical solution, we use the CG method to solve the system \eqref{eq:3d1d_system} preconditioned with the metric-perturbed AMG method described in \Cref{subsec:mamg}. The solver is executed through the call of the HAZniCS wrapper function \texttt{fenics\_metric\_amg\_solver()}, as presented in \Cref{lst:wrapper}. The solver parameters are set through the input file \texttt{input.dat}. The convergence is considered reached if the $l_2$ relative residual norm is less than $10^{-6}$. We choose the SA-AMG that uses the block Schwarz smoother (defined by the kernel decomposition \eqref{eq:kernel-decomposition}) for the interface degrees of freedom and standard Gauss-Seidel smoother on the interior degrees of freedom. We note by the interface degrees of freedom the sub-components of the 3$d$ variable that contribute to the interface current flow exchange, i.e., the nonzero components of ${\Pi^{\rho}_{\Gamma}} q_3 $ for $q_3 \in Q_3$. The application of the block Schwarz smoothers is done in a symmetric multiplicative way.

\begin{table}[htpb]
	\begin{tabular}{l||rrrrr||rrrrr}
		\hline 
		\multicolumn{1}{c||}{} & \multicolumn{10}{c}{$\log_{10}(\Delta t)$  $[s]$} \\ \cline{2-11}
		\multicolumn{1}{c||}{} & \multicolumn{1}{r}{-10} & \multicolumn{1}{r}{-8} & \multicolumn{1}{r}{-6} & \multicolumn{1}{r}{-4} & \multicolumn{1}{r||}{-2} & \multicolumn{1}{r}{-10} & \multicolumn{1}{r}{-8} & \multicolumn{1}{r}{-6} & \multicolumn{1}{r}{-4} & \multicolumn{1}{r}{-2} \\ \cline{2-11}
		\multicolumn{1}{c||}{$\rho [\mu m]$} & \multicolumn{5}{c||}{$N_{\text{iter}}$} & \multicolumn{5}{c}{$\kappa(BA)$} \\ \hline \hline
		5.0 	& 8 & 8 & 8 & 7 & 7 	& 2.793 & 2.778 & 2.667 & 2.061 & 1.671 \\
		1.0 	& 7 & 7 & 7 & 6 & 6 	& 1.941 & 1.952 & 1.947 & 1.413 & 1.408 \\
		0.5 	& 8 & 8 & 7 & 6 & 6 	& 2.541 & 2.515 & 2.276 & 1.527 & 1.410 \\
		0.1 	& 8 & 8 & 7 & 6 & 6 	& 2.583 & 2.562 & 2.257 & 1.492 & 1.413 \\
		\hline	
	\end{tabular}
	\caption{Performance of the CG method preconditioned with metric-perturbed AMG method from HAZniCS, with regards to parameters $\rho$ and $\Delta t$. We measure number of iterations ($N_{\text{iter}}$) of the solver with relative residual tolerance $10^{-6}$ and the approximate condition number $\kappa(BA)$ of the preconditioned system \eqref{eq:3d1d_system}. Results are obtained by running \texttt{HAZniCS-examples/demo\_3d1d.py}.}
	\label{tab:3d1d_eps_rho}
\end{table}

We study the performance of our solver with regard to the time step size $\Delta t$ and the coupling/dendrite radius $\rho$. Specifically, we are interested in the solver performance for very small time steps since this results in the metric term in the system \eqref{eq:3d1d_system} to dominate. The conductivity and membrane capacitance parameters remain constant and fixed throughout their respective domains to $\sigma_3 = 3 $ mS cm$^{-1} $, $\sigma_1 = 7 $ mS cm$^{-1} $ and $C_m = 1\, \mu$F cm$^{-2}$ \cite{buccino2019does}. The results given in \Cref{tab:3d1d_eps_rho} show stable number of iterations $N_{\text{iter}}$ and condition numbers $\kappa(\mathsf{B}\mathsf{A})$, where $\mathsf{A}$ is the matrix of the mixed-dimensional system \eqref{eq:3d1d_system} and $\mathsf{B} = \operatorname{metricAMG(\mathsf{A})}$ is the metric-perturbed AMG preconditioner. This is especially important in realistic cases when $\rho = 5 \mu$m and the time step magnitude is in nanoseconds, which is represented in the top left part of the table. In summary, the results show the method is robust with regard to problem parameters. Therefore we can confidently incorporate the method as part of the solver for the full EMI model \cite{Jaeger2021, buccino2021improving}.

\section{Conclusion} \label{sec:conclusion}
This paper introduces a collection of software solutions, HAZniCS, for solving interface-coupled multiphysics problems. The software combines two frameworks, HAZmath and FEniCS, into a flexible and powerful tool to obtain reliable and efficient simulators for various coupled problems. The focus of this work has been the (3$d$-2$d$ coupled) Darcy-Stokes model and the 3$d$-1$d$ coupled diffusion model for which we have presented the implementation and illustrated the performance of our solvers. In addition, we believe that the results shown in the paper demonstrate a great potential to utilize our framework in other relevant applications. The solver library allows interfacing with other finite element libraries which support the discretization of multiphysics problems, such as the new generation FEniCS platform FEniCSx or the Julia library Gridap.jl \cite{gridap}.

\begin{acks}
AB, MK, and KAM acknowledge the financial support from the SciML project funded by the Norwegian Research Council grant 102155. The work of XH is partially supported by the National Science Foundation under grant DMS-2208267. MK acknowledges support from Norwegian Research Council grant 303362. The research of LZ is supported in part by the U.~S.-Norway Fulbright Foundation and the U.~S. National Science Foundation grant DMS-2208249. The collaborative efforts of XH and LZ were supported in part by the NSF DMS-2132710 through the Workshop on Numerical Modeling with Neural Networks, Learning, and Multilevel FE.
\end{acks}

\bibliographystyle{ACM-Reference-Format}
\bibliography{multiphy_ref,haznics_AMG}


\begin{thebibliography}{84}


\ifx \showCODEN    \undefined \def \showCODEN     #1{\unskip}     \fi
\ifx \showDOI      \undefined \def \showDOI       #1{#1}\fi
\ifx \showISBNx    \undefined \def \showISBNx     #1{\unskip}     \fi
\ifx \showISBNxiii \undefined \def \showISBNxiii  #1{\unskip}     \fi
\ifx \showISSN     \undefined \def \showISSN      #1{\unskip}     \fi
\ifx \showLCCN     \undefined \def \showLCCN      #1{\unskip}     \fi
\ifx \shownote     \undefined \def \shownote      #1{#1}          \fi
\ifx \showarticletitle \undefined \def \showarticletitle #1{#1}   \fi
\ifx \showURL      \undefined \def \showURL       {\relax}        \fi
\providecommand\bibfield[2]{#2}
\providecommand\bibinfo[2]{#2}
\providecommand\natexlab[1]{#1}
\providecommand\showeprint[2][]{arXiv:#2}

\bibitem[\protect\citeauthoryear{Adler, Hu, and Zikatanov}{Adler
  et~al\mbox{.}}{2009}]%
        {hazmath}
\bibfield{author}{\bibinfo{person}{J.~Adler}, \bibinfo{person}{X.~Hu}, {and}
  \bibinfo{person}{L.~Zikatanov}.} \bibinfo{year}{2009}\natexlab{}.
\newblock \bibinfo{booktitle}{\emph{HAZMATH: A Simple Finite Element, Graph,
  and Solver Library}}.
\newblock
\urldef\tempurl%
\url{https://hazmathteam.github.io/hazmath/}
\showURL{%
\tempurl}


\bibitem[\protect\citeauthoryear{Agoshkov}{Agoshkov}{1988}]%
        {agoshkov1988poincare}
\bibfield{author}{\bibinfo{person}{V.~I. Agoshkov}.}
  \bibinfo{year}{1988}\natexlab{}.
\newblock \showarticletitle{Poincar{\'e}-Steklov operators and domain
  decomposition methods in finite dimensional spaces}. In
  \bibinfo{booktitle}{\emph{First International Symposium on Domain
  Decomposition Methods for Partial Differential Equations}}.
  \bibinfo{pages}{73--112}.
\newblock


\bibitem[\protect\citeauthoryear{Axelsson and Vassilevski}{Axelsson and
  Vassilevski}{1989}]%
        {AV_89}
\bibfield{author}{\bibinfo{person}{O.~Axelsson} {and}
  \bibinfo{person}{P.~Vassilevski}.} \bibinfo{year}{1989}\natexlab{}.
\newblock \showarticletitle{Algebraic multilevel preconditioning methods. {I}}.
\newblock \bibinfo{journal}{\emph{Numer. Math.}}  \bibinfo{volume}{56}
  (\bibinfo{year}{1989}), \bibinfo{pages}{157--177}.
\newblock


\bibitem[\protect\citeauthoryear{Axelsson and Vassilevski}{Axelsson and
  Vassilevski}{1990}]%
        {AV_90}
\bibfield{author}{\bibinfo{person}{O.~Axelsson} {and} \bibinfo{person}{P.~S.
  Vassilevski}.} \bibinfo{year}{1990}\natexlab{}.
\newblock \showarticletitle{Algebraic multilevel preconditioning methods.
  {II}}.
\newblock \bibinfo{journal}{\emph{SIAM J. Numer. Anal.}} \bibinfo{volume}{27},
  \bibinfo{number}{6} (\bibinfo{year}{1990}), \bibinfo{pages}{1569--1590}.
\newblock


\bibitem[\protect\citeauthoryear{Axelsson and Vassilevski}{Axelsson and
  Vassilevski}{1991}]%
        {AV_91}
\bibfield{author}{\bibinfo{person}{O.~Axelsson} {and} \bibinfo{person}{P.~S.
  Vassilevski}.} \bibinfo{year}{1991}\natexlab{}.
\newblock \showarticletitle{A black box generalized conjugate gradient solver
  with inner iterations and variable-step preconditioning}.
\newblock \bibinfo{journal}{\emph{SIAM J. Matrix Anal. Appl.}}
  \bibinfo{volume}{12}, \bibinfo{number}{4} (\bibinfo{year}{1991}),
  \bibinfo{pages}{625--644}.
\newblock
\showCODEN{SJMAEL}
\showISSN{0895-4798}


\bibitem[\protect\citeauthoryear{Axelsson and Vassilevski}{Axelsson and
  Vassilevski}{1994}]%
        {AV_94}
\bibfield{author}{\bibinfo{person}{O.~Axelsson} {and} \bibinfo{person}{P.~S.
  Vassilevski}.} \bibinfo{year}{1994}\natexlab{}.
\newblock \showarticletitle{Variable-step multilevel preconditioning methods.
  {I}. {S}elfadjoint and positive definite elliptic problems}.
\newblock \bibinfo{journal}{\emph{Numer. Linear Algebra Appl.}}
  \bibinfo{volume}{1}, \bibinfo{number}{1} (\bibinfo{year}{1994}),
  \bibinfo{pages}{75--101}.
\newblock
\showCODEN{NLAAEM}
\showISSN{1070-5325}


\bibitem[\protect\citeauthoryear{Babu\v{s}ka}{Babu\v{s}ka}{1971}]%
        {babuska1971}
\bibfield{author}{\bibinfo{person}{I.~Babu\v{s}ka}.}
  \bibinfo{year}{1971}\natexlab{}.
\newblock \showarticletitle{Error-bounds for finite element method}.
\newblock \bibinfo{journal}{\emph{Numer. Math.}}  \bibinfo{volume}{16}
  (\bibinfo{year}{1971}), \bibinfo{pages}{322--333}.
\newblock
Issue 4.
\urldef\tempurl%
\url{https://doi.org/10.1007/BF02165003}
\showDOI{\tempurl}


\bibitem[\protect\citeauthoryear{Babu\v{s}ka and Aziz}{Babu\v{s}ka and
  Aziz}{1972}]%
        {babuska_aziz}
\bibfield{author}{\bibinfo{person}{I.~Babu\v{s}ka} {and} \bibinfo{person}{A.~K.
  Aziz}.} \bibinfo{year}{1972}\natexlab{}.
\newblock \bibinfo{booktitle}{\emph{Survey lectures on the mathematical
  foundation of the finite element method}}.
\newblock \bibinfo{publisher}{Academic Press}, \bibinfo{address}{New York,
  London}, \bibinfo{pages}{3--345}.
\newblock
\showISBNx{978-0-12-068650-6}
\urldef\tempurl%
\url{https://doi.org/10.1016/C2013-0-10319-4}
\showDOI{\tempurl}


\bibitem[\protect\citeauthoryear{Badia, Nobile, and Vergara}{Badia
  et~al\mbox{.}}{2009}]%
        {badia2009robin}
\bibfield{author}{\bibinfo{person}{S.~Badia}, \bibinfo{person}{F.~Nobile},
  {and} \bibinfo{person}{C.~Vergara}.} \bibinfo{year}{2009}\natexlab{}.
\newblock \showarticletitle{Robin--Robin preconditioned Krylov methods for
  fluid--structure interaction problems}.
\newblock \bibinfo{journal}{\emph{Computer Methods in Applied Mechanics and
  Engineering}} \bibinfo{volume}{198}, \bibinfo{number}{33-36}
  (\bibinfo{year}{2009}), \bibinfo{pages}{2768--2784}.
\newblock


\bibitem[\protect\citeauthoryear{B{\ae}rland}{B{\ae}rland}{2019}]%
        {baerland2019auxiliary}
\bibfield{author}{\bibinfo{person}{T.~B{\ae}rland}.}
  \bibinfo{year}{2019}\natexlab{}.
\newblock \showarticletitle{An auxiliary space preconditioner for fractional
  {L}aplacian of negative order}.
\newblock \bibinfo{journal}{\emph{arXiv preprint arXiv:1908.04498}}
  (\bibinfo{year}{2019}).
\newblock


\bibitem[\protect\citeauthoryear{B{\ae}rland, Kuchta, and Mardal}{B{\ae}rland
  et~al\mbox{.}}{2019}]%
        {Baerland2019}
\bibfield{author}{\bibinfo{person}{T.~B{\ae}rland},
  \bibinfo{person}{M.~Kuchta}, {and} \bibinfo{person}{K.-A. Mardal}.}
  \bibinfo{year}{2019}\natexlab{}.
\newblock \showarticletitle{{Multigrid methods for discrete fractional Sobolev
  spaces}}.
\newblock \bibinfo{journal}{\emph{SIAM Journal on Scientific Computing}}
  \bibinfo{volume}{41}, \bibinfo{number}{2} (\bibinfo{year}{2019}),
  \bibinfo{pages}{A948--A972}.
\newblock
\showISSN{10957197}
\urldef\tempurl%
\url{https://doi.org/10.1137/18M1191488}
\showDOI{\tempurl}


\bibitem[\protect\citeauthoryear{Balay, Abhyankar, Adams, Benson, Brown, Brune,
  Buschelman, and et. al.}{Balay et~al\mbox{.}}{2022}]%
        {petsc-user-ref}
\bibfield{author}{\bibinfo{person}{S.~Balay}, \bibinfo{person}{S.~Abhyankar},
  \bibinfo{person}{M.~F. Adams}, \bibinfo{person}{S.~Benson},
  \bibinfo{person}{J.~Brown}, \bibinfo{person}{P.~Brune},
  \bibinfo{person}{K.~Buschelman}, {and} \bibinfo{person}{et. al.}}
  \bibinfo{year}{2022}\natexlab{}.
\newblock \bibinfo{booktitle}{\emph{{PETSc/TAO} Users Manual}}.
\newblock \bibinfo{type}{{T}echnical {R}eport} ANL-21/39 - Revision 3.18.
  \bibinfo{institution}{Argonne National Laboratory}.
\newblock


\bibitem[\protect\citeauthoryear{Beazley}{Beazley}{1996}]%
        {swig}
\bibfield{author}{\bibinfo{person}{D.~M. Beazley}.}
  \bibinfo{year}{1996}\natexlab{}.
\newblock \showarticletitle{SWIG: An Easy to Use Tool for Integrating Scripting
  Languages with C and C++}. In \bibinfo{booktitle}{\emph{Proceedings of the
  4th Conference on USENIX Tcl/Tk Workshop, 1996 - Volume 4}} (Monterey,
  California) \emph{(\bibinfo{series}{TCLTK'96})}. \bibinfo{publisher}{USENIX
  Association}, \bibinfo{address}{USA}, \bibinfo{pages}{15}.
\newblock
\urldef\tempurl%
\url{https://www.swig.org/}
\showURL{%
\tempurl}


\bibitem[\protect\citeauthoryear{Berg, Davit, Quintard, and Lorthois}{Berg
  et~al\mbox{.}}{2020}]%
        {berg_davit_quintard_lorthois_2020}
\bibfield{author}{\bibinfo{person}{M.~Berg}, \bibinfo{person}{Y.~Davit},
  \bibinfo{person}{M.~Quintard}, {and} \bibinfo{person}{S.~Lorthois}.}
  \bibinfo{year}{2020}\natexlab{}.
\newblock \showarticletitle{Modelling solute transport in the brain
  microcirculation: is it really well mixed inside the blood vessels?}
\newblock \bibinfo{journal}{\emph{Journal of Fluid Mechanics}}
  \bibinfo{volume}{884} (\bibinfo{year}{2020}), \bibinfo{pages}{A39}.
\newblock
\urldef\tempurl%
\url{https://doi.org/10.1017/jfm.2019.866}
\showDOI{\tempurl}


\bibitem[\protect\citeauthoryear{Blaheta}{Blaheta}{1986}]%
        {Blaheta.R.1986a}
\bibfield{author}{\bibinfo{person}{R.~Blaheta}.}
  \bibinfo{year}{1986}\natexlab{}.
\newblock \showarticletitle{{A multilevel method with correction by aggregation
  for solving discrete elliptic problems}}.
\newblock \bibinfo{journal}{\emph{Aplikace matematiky}} \bibinfo{volume}{31},
  \bibinfo{number}{5} (\bibinfo{year}{1986}), \bibinfo{pages}{365--378}.
\newblock


\bibitem[\protect\citeauthoryear{Boon, Hornkj{\o}l, Kuchta, Mardal, and
  Ruiz-Baier}{Boon et~al\mbox{.}}{2022a}]%
        {hornkjol2021}
\bibfield{author}{\bibinfo{person}{W.~M. Boon},
  \bibinfo{person}{M.~Hornkj{\o}l}, \bibinfo{person}{M.~Kuchta},
  \bibinfo{person}{K.-A. Mardal}, {and} \bibinfo{person}{R.~Ruiz-Baier}.}
  \bibinfo{year}{2022}\natexlab{a}.
\newblock \showarticletitle{Parameter-robust methods for the {B}iot--{S}tokes
  interfacial coupling without {L}agrange multipliers}.
\newblock \bibinfo{journal}{\emph{J. Comput. Phys.}}  \bibinfo{volume}{467}
  (\bibinfo{year}{2022}), \bibinfo{pages}{111464}.
\newblock


\bibitem[\protect\citeauthoryear{Boon, Koch, Kuchta, and Mardal}{Boon
  et~al\mbox{.}}{2022b}]%
        {boon2021robust}
\bibfield{author}{\bibinfo{person}{W.~M. Boon}, \bibinfo{person}{T.~Koch},
  \bibinfo{person}{M.~Kuchta}, {and} \bibinfo{person}{K.-A. Mardal}.}
  \bibinfo{year}{2022}\natexlab{b}.
\newblock \showarticletitle{Robust Monolithic Solvers for the {S}tokes--{D}arcy
  Problem with the {D}arcy Equation in Primal Form}.
\newblock \bibinfo{journal}{\emph{SIAM Journal on Scientific Computing}}
  \bibinfo{volume}{44}, \bibinfo{number}{4} (\bibinfo{year}{2022}),
  \bibinfo{pages}{B1148--B1174}.
\newblock


\bibitem[\protect\citeauthoryear{Bramble, Pasciak, and Vassilevski}{Bramble
  et~al\mbox{.}}{2000}]%
        {bramble2000computational}
\bibfield{author}{\bibinfo{person}{J.~Bramble}, \bibinfo{person}{J.~Pasciak},
  {and} \bibinfo{person}{P.~Vassilevski}.} \bibinfo{year}{2000}\natexlab{}.
\newblock \showarticletitle{Computational scales of {S}obolev norms with
  application to preconditioning}.
\newblock \bibinfo{journal}{\emph{Math. Comp.}} \bibinfo{volume}{69},
  \bibinfo{number}{230} (\bibinfo{year}{2000}), \bibinfo{pages}{463--480}.
\newblock


\bibitem[\protect\citeauthoryear{Bramble}{Bramble}{2019}]%
        {bramble2019multigrid}
\bibfield{author}{\bibinfo{person}{J.~H. Bramble}.}
  \bibinfo{year}{2019}\natexlab{}.
\newblock \bibinfo{booktitle}{\emph{Multigrid methods}}.
\newblock \bibinfo{publisher}{Chapman and Hall/CRC}.
\newblock


\bibitem[\protect\citeauthoryear{Brandt, McCormick, and Ruge}{Brandt
  et~al\mbox{.}}{1982}]%
        {Brandt.A;McCormick.S;Ruge.J.1982a}
\bibfield{author}{\bibinfo{person}{A.~Brandt}, \bibinfo{person}{S.~F.
  McCormick}, {and} \bibinfo{person}{J.~W. Ruge}.}
  \bibinfo{year}{1982}\natexlab{}.
\newblock \bibinfo{booktitle}{\emph{Algebraic multigrid {(AMG)} for automatic
  multigrid solutions with application to geodetic computations}}.
\newblock \bibinfo{type}{{T}echnical {R}eport}. \bibinfo{institution}{Inst. for
  Computational Studies}, \bibinfo{address}{Fort Collins, CO}.
\newblock


\bibitem[\protect\citeauthoryear{Buccino, Kuchta, Jæger, Ness, Berthet,
  Mardal, Cauwenberghs, and Tveito}{Buccino et~al\mbox{.}}{2019}]%
        {buccino2019does}
\bibfield{author}{\bibinfo{person}{A.~P. Buccino}, \bibinfo{person}{M.~Kuchta},
  \bibinfo{person}{K.~H. Jæger}, \bibinfo{person}{T.~V. Ness},
  \bibinfo{person}{P.~Berthet}, \bibinfo{person}{K.-A. Mardal},
  \bibinfo{person}{G.~Cauwenberghs}, {and} \bibinfo{person}{A.~Tveito}.}
  \bibinfo{year}{2019}\natexlab{}.
\newblock \showarticletitle{How does the presence of neural probes affect
  extracellular potentials?}
\newblock \bibinfo{journal}{\emph{Journal of Neural Engineering}}
  \bibinfo{volume}{16} (\bibinfo{date}{4} \bibinfo{year}{2019}),
  \bibinfo{pages}{026030}.
\newblock
Issue 2.
\showISSN{1741-2560}
\urldef\tempurl%
\url{https://doi.org/10.1088/1741-2552/ab03a1}
\showDOI{\tempurl}


\bibitem[\protect\citeauthoryear{Buccino, Kuchta, Schreiner, and
  Mardal}{Buccino et~al\mbox{.}}{2021}]%
        {buccino2021improving}
\bibfield{author}{\bibinfo{person}{A.~P. Buccino}, \bibinfo{person}{M.~Kuchta},
  \bibinfo{person}{J.~Schreiner}, {and} \bibinfo{person}{K.-A. Mardal}.}
  \bibinfo{year}{2021}\natexlab{}.
\newblock \bibinfo{booktitle}{\emph{Improving Neural Simulations with the EMI
  Model}}.
\newblock \bibinfo{publisher}{Springer International Publishing},
  \bibinfo{address}{Cham}, \bibinfo{pages}{87--98}.
\newblock
\showISBNx{978-3-030-61157-6}
\urldef\tempurl%
\url{https://doi.org/10.1007/978-3-030-61157-6_7}
\showDOI{\tempurl}


\bibitem[\protect\citeauthoryear{Budi\v{s}a, Hu, Kuchta, Mardal, and
  Zikatanov}{Budi\v{s}a et~al\mbox{.}}{2022a}]%
        {github_haznics}
\bibfield{author}{\bibinfo{person}{A.~Budi\v{s}a}, \bibinfo{person}{X.~Hu},
  \bibinfo{person}{M.~Kuchta}, \bibinfo{person}{K.-A. Mardal}, {and}
  \bibinfo{person}{L.~Zikatanov}.} \bibinfo{year}{2022}\natexlab{a}.
\newblock \bibinfo{booktitle}{\emph{HAZniCS examples.}}
\newblock
\urldef\tempurl%
\url{https://doi.org/10.5281/zenodo.7220688}
\showDOI{\tempurl}


\bibitem[\protect\citeauthoryear{Budi\v{s}a, Hu, Kuchta, Mardal, and
  Zikatanov}{Budi\v{s}a et~al\mbox{.}}{2022b}]%
        {budisa2022rational}
\bibfield{author}{\bibinfo{person}{A.~Budi\v{s}a}, \bibinfo{person}{X.~Hu},
  \bibinfo{person}{M.~Kuchta}, \bibinfo{person}{K.-A. Mardal}, {and}
  \bibinfo{person}{L.~Zikatanov}.} \bibinfo{year}{2022}\natexlab{b}.
\newblock \showarticletitle{Rational approximation preconditioners for
  multiphysics problems}.
\newblock \bibinfo{journal}{\emph{arXiv preprint arXiv:2209.11659}}
  (\bibinfo{year}{2022}).
\newblock


\bibitem[\protect\citeauthoryear{Cerroni, Laurino, and Zunino}{Cerroni
  et~al\mbox{.}}{2019}]%
        {cerroni2019mathematical}
\bibfield{author}{\bibinfo{person}{D.~Cerroni}, \bibinfo{person}{F.~Laurino},
  {and} \bibinfo{person}{P.~Zunino}.} \bibinfo{year}{2019}\natexlab{}.
\newblock \showarticletitle{Mathematical analysis, finite element approximation
  and numerical solvers for the interaction of 3d reservoirs with 1d wells}.
\newblock \bibinfo{journal}{\emph{GEM-International Journal on Geomathematics}}
  \bibinfo{volume}{10}, \bibinfo{number}{1} (\bibinfo{year}{2019}),
  \bibinfo{pages}{1--27}.
\newblock


\bibitem[\protect\citeauthoryear{D'Ambra and Vassilevski}{D'Ambra and
  Vassilevski}{2014}]%
        {DAmbra.P;Vassilevski.P.2014a}
\bibfield{author}{\bibinfo{person}{P.~D'Ambra} {and}
  \bibinfo{person}{P.~Vassilevski}.} \bibinfo{year}{2014}\natexlab{}.
\newblock \bibinfo{title}{Adaptive {AMG} based on weighted matching for systems
  of elliptic {PDE}s arising from displacements and mixed methods}.
\newblock
\newblock
\newblock
\shownote{Submitted.}


\bibitem[\protect\citeauthoryear{Deparis, Discacciati, Fourestey, and
  Quarteroni}{Deparis et~al\mbox{.}}{2006}]%
        {deparis2006fluid}
\bibfield{author}{\bibinfo{person}{S.~Deparis},
  \bibinfo{person}{M.~Discacciati}, \bibinfo{person}{G.~Fourestey}, {and}
  \bibinfo{person}{A.~Quarteroni}.} \bibinfo{year}{2006}\natexlab{}.
\newblock \showarticletitle{Fluid--structure algorithms based on
  Steklov--Poincar{\'e} operators}.
\newblock \bibinfo{journal}{\emph{Computer Methods in Applied Mechanics and
  Engineering}} \bibinfo{volume}{195}, \bibinfo{number}{41-43}
  (\bibinfo{year}{2006}), \bibinfo{pages}{5797--5812}.
\newblock


\bibitem[\protect\citeauthoryear{D’Angelo and Quarteroni}{D’Angelo and
  Quarteroni}{2008}]%
        {DAngelo2008}
\bibfield{author}{\bibinfo{person}{C.~D’Angelo} {and}
  \bibinfo{person}{A.~Quarteroni}.} \bibinfo{year}{2008}\natexlab{}.
\newblock \showarticletitle{{On the coupling of 1d and 3d diffusion-reaction
  equations: Application to tissue perfusion problems}}.
\newblock \bibinfo{journal}{\emph{Mathematical Models and Methods in Applied
  Sciences}} \bibinfo{volume}{18}, \bibinfo{number}{8} (\bibinfo{year}{2008}),
  \bibinfo{pages}{1481--1504}.
\newblock


\bibitem[\protect\citeauthoryear{Eide, Vinje, Pripp, Mardal, and Ringstad}{Eide
  et~al\mbox{.}}{2021}]%
        {eide2021sleep}
\bibfield{author}{\bibinfo{person}{P.~K. Eide}, \bibinfo{person}{V.~Vinje},
  \bibinfo{person}{A.~H. Pripp}, \bibinfo{person}{K.-A. Mardal}, {and}
  \bibinfo{person}{G.~Ringstad}.} \bibinfo{year}{2021}\natexlab{}.
\newblock \showarticletitle{Sleep deprivation impairs molecular clearance from
  the human brain}.
\newblock \bibinfo{journal}{\emph{Brain}} \bibinfo{volume}{144},
  \bibinfo{number}{3} (\bibinfo{year}{2021}), \bibinfo{pages}{863--874}.
\newblock


\bibitem[\protect\citeauthoryear{Falgout and Yang}{Falgout and Yang}{2002}]%
        {hypre}
\bibfield{author}{\bibinfo{person}{R.~D. Falgout} {and} \bibinfo{person}{U.~M.
  Yang}.} \bibinfo{year}{2002}\natexlab{}.
\newblock \showarticletitle{hypre: A Library of High Performance
  Preconditioners}. In \bibinfo{booktitle}{\emph{Computational Science --- ICCS
  2002}}, \bibfield{editor}{\bibinfo{person}{Peter M.~A. Sloot},
  \bibinfo{person}{Alfons~G. Hoekstra}, \bibinfo{person}{C.~J.~Kenneth Tan},
  {and} \bibinfo{person}{Jack~J. Dongarra}} (Eds.).
  \bibinfo{publisher}{Springer Berlin Heidelberg}, \bibinfo{address}{Berlin,
  Heidelberg}, \bibinfo{pages}{632--641}.
\newblock


\bibitem[\protect\citeauthoryear{F{\"u}hrer}{F{\"u}hrer}{2022}]%
        {fuhrer2022multilevel}
\bibfield{author}{\bibinfo{person}{T.~F{\"u}hrer}.}
  \bibinfo{year}{2022}\natexlab{}.
\newblock \showarticletitle{Multilevel decompositions and norms for negative
  order Sobolev spaces}.
\newblock \bibinfo{journal}{\emph{Math. Comp.}} \bibinfo{volume}{91},
  \bibinfo{number}{333} (\bibinfo{year}{2022}), \bibinfo{pages}{183--218}.
\newblock


\bibitem[\protect\citeauthoryear{Gjerde, Kumar, and Nordbotten}{Gjerde
  et~al\mbox{.}}{2020}]%
        {gjerde2020singularity}
\bibfield{author}{\bibinfo{person}{I.~G. Gjerde}, \bibinfo{person}{K.~Kumar},
  {and} \bibinfo{person}{J.~M. Nordbotten}.} \bibinfo{year}{2020}\natexlab{}.
\newblock \showarticletitle{A singularity removal method for coupled 1{D}--3{D}
  flow models}.
\newblock \bibinfo{journal}{\emph{Computational Geosciences}}
  \bibinfo{volume}{24}, \bibinfo{number}{2} (\bibinfo{year}{2020}),
  \bibinfo{pages}{443--457}.
\newblock


\bibitem[\protect\citeauthoryear{Harizanov, Lazarov, Margenov, and
  Marinov}{Harizanov et~al\mbox{.}}{2020}]%
        {harizanov2020numerical}
\bibfield{author}{\bibinfo{person}{S.~Harizanov}, \bibinfo{person}{R.~Lazarov},
  \bibinfo{person}{S.~Margenov}, {and} \bibinfo{person}{P.~Marinov}.}
  \bibinfo{year}{2020}\natexlab{}.
\newblock \showarticletitle{Numerical solution of fractional
  diffusion--reaction problems based on BURA}.
\newblock \bibinfo{journal}{\emph{Computers \& Mathematics with Applications}}
  \bibinfo{volume}{80}, \bibinfo{number}{2} (\bibinfo{year}{2020}),
  \bibinfo{pages}{316--331}.
\newblock


\bibitem[\protect\citeauthoryear{Harizanov, Lazarov, Margenov, Marinov, and
  Vutov}{Harizanov et~al\mbox{.}}{2018}]%
        {harizanov2018optimal}
\bibfield{author}{\bibinfo{person}{S.~Harizanov}, \bibinfo{person}{R.~Lazarov},
  \bibinfo{person}{S.~Margenov}, \bibinfo{person}{P.~Marinov}, {and}
  \bibinfo{person}{Y.~Vutov}.} \bibinfo{year}{2018}\natexlab{}.
\newblock \showarticletitle{Optimal solvers for linear systems with fractional
  powers of sparse SPD matrices}.
\newblock \bibinfo{journal}{\emph{Numerical Linear Algebra with Applications}}
  \bibinfo{volume}{25}, \bibinfo{number}{5} (\bibinfo{year}{2018}),
  \bibinfo{pages}{e2167}.
\newblock


\bibitem[\protect\citeauthoryear{Harizanov, Lirkov, and Margenov}{Harizanov
  et~al\mbox{.}}{2022}]%
        {harizanov2022rational}
\bibfield{author}{\bibinfo{person}{S.~Harizanov}, \bibinfo{person}{I.~Lirkov},
  {and} \bibinfo{person}{S.~Margenov}.} \bibinfo{year}{2022}\natexlab{}.
\newblock \showarticletitle{Rational Approximations in Robust Preconditioning
  of Multiphysics Problems}.
\newblock \bibinfo{journal}{\emph{Mathematics}} \bibinfo{volume}{10},
  \bibinfo{number}{5} (\bibinfo{year}{2022}), \bibinfo{pages}{780}.
\newblock


\bibitem[\protect\citeauthoryear{Hartung, Badr, Mihelic, Dunn, Cheng, Kura,
  Boas, Kleinfeld, Alaraj, and Linninger}{Hartung et~al\mbox{.}}{2021}]%
        {hartung2021mathematical}
\bibfield{author}{\bibinfo{person}{G.~Hartung}, \bibinfo{person}{S.~Badr},
  \bibinfo{person}{S.~Mihelic}, \bibinfo{person}{A.~Dunn},
  \bibinfo{person}{X.~Cheng}, \bibinfo{person}{S.~Kura}, \bibinfo{person}{D.~A.
  Boas}, \bibinfo{person}{D.~Kleinfeld}, \bibinfo{person}{A.~Alaraj}, {and}
  \bibinfo{person}{A.~A. Linninger}.} \bibinfo{year}{2021}\natexlab{}.
\newblock \showarticletitle{Mathematical synthesis of the cortical circulation
  for the whole mouse brain—part II: Microcirculatory closure}.
\newblock \bibinfo{journal}{\emph{Microcirculation}} \bibinfo{volume}{28},
  \bibinfo{number}{5} (\bibinfo{year}{2021}), \bibinfo{pages}{e12687}.
\newblock


\bibitem[\protect\citeauthoryear{Hofreither}{Hofreither}{2020}]%
        {Hofreither2020}
\bibfield{author}{\bibinfo{person}{C.~Hofreither}.}
  \bibinfo{year}{2020}\natexlab{}.
\newblock \showarticletitle{{A unified view of some numerical methods for
  fractional diffusion}}.
\newblock \bibinfo{journal}{\emph{Computers and Mathematics with Applications}}
  \bibinfo{volume}{80}, \bibinfo{number}{2} (\bibinfo{year}{2020}),
  \bibinfo{pages}{332--350}.
\newblock
\showISSN{08981221}
\urldef\tempurl%
\url{https://doi.org/10.1016/j.camwa.2019.07.025}
\showDOI{\tempurl}


\bibitem[\protect\citeauthoryear{Holter, Kuchta, and Mardal}{Holter
  et~al\mbox{.}}{2020}]%
        {holter2020}
\bibfield{author}{\bibinfo{person}{K.-E. Holter}, \bibinfo{person}{M.~Kuchta},
  {and} \bibinfo{person}{K.-A. Mardal}.} \bibinfo{year}{2020}\natexlab{}.
\newblock \showarticletitle{Robust preconditioning of monolithically coupled
  multiphysics problems}.
\newblock \bibinfo{journal}{\emph{arXiv preprint arXiv:2001.05527}}
  (\bibinfo{year}{2020}).
\newblock


\bibitem[\protect\citeauthoryear{Holter, Kuchta, and Mardal}{Holter
  et~al\mbox{.}}{2021}]%
        {holter2021robust}
\bibfield{author}{\bibinfo{person}{K.~E. Holter}, \bibinfo{person}{M.~Kuchta},
  {and} \bibinfo{person}{K.-A. Mardal}.} \bibinfo{year}{2021}\natexlab{}.
\newblock \showarticletitle{Robust preconditioning for coupled Stokes--Darcy
  problems with the Darcy problem in primal form}.
\newblock \bibinfo{journal}{\emph{Computers \& Mathematics with Applications}}
  \bibinfo{volume}{91} (\bibinfo{year}{2021}), \bibinfo{pages}{53--66}.
\newblock


\bibitem[\protect\citeauthoryear{Hu, Lin, and Zikatanov}{Hu
  et~al\mbox{.}}{2019}]%
        {HuLinZikatanov2019}
\bibfield{author}{\bibinfo{person}{X.~Hu}, \bibinfo{person}{J.~Lin}, {and}
  \bibinfo{person}{L.~T. Zikatanov}.} \bibinfo{year}{2019}\natexlab{}.
\newblock \showarticletitle{An {Adaptive} {Multigrid} {Method} {Based} on
  {Path} {Cover}}.
\newblock \bibinfo{journal}{\emph{SIAM Journal on Scientific Computing}}
  \bibinfo{volume}{41}, \bibinfo{number}{5} (\bibinfo{date}{Jan.}
  \bibinfo{year}{2019}), \bibinfo{pages}{S220--S241}.
\newblock
\showISSN{1064-8275, 1095-7197}
\urldef\tempurl%
\url{https://doi.org/10.1137/18M1194493}
\showDOI{\tempurl}
\newblock
\shownote{Citation Key Alias: HuLinZikatanov2018.}


\bibitem[\protect\citeauthoryear{Hu, Vassilevski, and Xu}{Hu
  et~al\mbox{.}}{2013a}]%
        {HuVassilevskiXu2013}
\bibfield{author}{\bibinfo{person}{X.~Hu}, \bibinfo{person}{P.~S. Vassilevski},
  {and} \bibinfo{person}{J.~Xu}.} \bibinfo{year}{2013}\natexlab{a}.
\newblock \showarticletitle{Comparative {Convergence} {Analysis} of {Nonlinear}
  {AMLI}-{Cycle} {Multigrid}}.
\newblock \bibinfo{journal}{\emph{SIAM J. Numer. Anal.}} \bibinfo{volume}{51},
  \bibinfo{number}{2} (\bibinfo{date}{Jan.} \bibinfo{year}{2013}),
  \bibinfo{pages}{1349--1369}.
\newblock
\showISSN{0036-1429, 1095-7170}
\urldef\tempurl%
\url{https://doi.org/10.1137/110850049}
\showDOI{\tempurl}


\bibitem[\protect\citeauthoryear{Hu, Vassilevski, and Xu}{Hu
  et~al\mbox{.}}{2016}]%
        {HuVassilevskiXu2016a}
\bibfield{author}{\bibinfo{person}{X.~Hu}, \bibinfo{person}{P.~S. Vassilevski},
  {and} \bibinfo{person}{J.~Xu}.} \bibinfo{year}{2016}\natexlab{}.
\newblock \showarticletitle{A two-grid {SA}-{AMG} convergence bound that
  improves when increasing the polynomial degree: {I}mproving {TG} convergence
  with increasing smoothing steps}.
\newblock \bibinfo{journal}{\emph{Numerical Linear Algebra with Applications}}
  \bibinfo{volume}{23}, \bibinfo{number}{4} (\bibinfo{date}{Aug.}
  \bibinfo{year}{2016}), \bibinfo{pages}{746--771}.
\newblock
\showISSN{10705325}
\urldef\tempurl%
\url{https://doi.org/10.1002/nla.2053}
\showDOI{\tempurl}


\bibitem[\protect\citeauthoryear{Hu, Wu, and Zikatanov}{Hu
  et~al\mbox{.}}{2020}]%
        {HuWuZikatanov2020}
\bibfield{author}{\bibinfo{person}{X.~Hu}, \bibinfo{person}{K.~Wu}, {and}
  \bibinfo{person}{L.~T. Zikatanov}.} \bibinfo{year}{2020}\natexlab{}.
\newblock \showarticletitle{A {Posteriori} {Error} {Estimates} for {Multilevel}
  {Methods} for {Graph} {Laplacians}}.
\newblock \bibinfo{journal}{\emph{arXiv:2007.00189 [cs, math]}}
  (\bibinfo{date}{June} \bibinfo{year}{2020}).
\newblock
\urldef\tempurl%
\url{http://arxiv.org/abs/2007.00189}
\showURL{%
\tempurl}
\newblock
\shownote{arXiv: 2007.00189.}


\bibitem[\protect\citeauthoryear{Hu, Wu, Wu, Xu, Zhang, Zhang, and
  Zikatanov}{Hu et~al\mbox{.}}{2013b}]%
        {hu2013combined}
\bibfield{author}{\bibinfo{person}{X.~Hu}, \bibinfo{person}{S.~Wu},
  \bibinfo{person}{X.-H. Wu}, \bibinfo{person}{J.~Xu}, \bibinfo{person}{C.-S.
  Zhang}, \bibinfo{person}{S.~Zhang}, {and} \bibinfo{person}{L.~Zikatanov}.}
  \bibinfo{year}{2013}\natexlab{b}.
\newblock \showarticletitle{Combined preconditioning with applications in
  reservoir simulation}.
\newblock \bibinfo{journal}{\emph{Multiscale Modeling \& Simulation}}
  \bibinfo{volume}{11}, \bibinfo{number}{2} (\bibinfo{year}{2013}),
  \bibinfo{pages}{507--521}.
\newblock


\bibitem[\protect\citeauthoryear{Iliff, Wang, Liao, Plogg, Peng, Gundersen,
  Benveniste, Vates, Deane, Goldman, Nagelhus, and Nedergaard}{Iliff
  et~al\mbox{.}}{2012}]%
        {iliff2012paravascular}
\bibfield{author}{\bibinfo{person}{J.~J. Iliff}, \bibinfo{person}{M.~Wang},
  \bibinfo{person}{Y.~Liao}, \bibinfo{person}{B.~A. Plogg},
  \bibinfo{person}{W.~Peng}, \bibinfo{person}{G.~A. Gundersen},
  \bibinfo{person}{H.~Benveniste}, \bibinfo{person}{G.~E. Vates},
  \bibinfo{person}{R.~Deane}, \bibinfo{person}{S.~A. Goldman},
  \bibinfo{person}{E.~A. Nagelhus}, {and} \bibinfo{person}{M.~Nedergaard}.}
  \bibinfo{year}{2012}\natexlab{}.
\newblock \showarticletitle{A Paravascular Pathway Facilitates CSF Flow Through
  the Brain Parenchyma and the Clearance of Interstitial Solutes, Including
  Amyloid $\beta$}.
\newblock \bibinfo{journal}{\emph{Science Translational Medicine}}
  \bibinfo{volume}{4}, \bibinfo{number}{147} (\bibinfo{year}{2012}),
  \bibinfo{pages}{147ra111}.
\newblock
\urldef\tempurl%
\url{https://doi.org/10.1126/scitranslmed.3003748}
\showDOI{\tempurl}


\bibitem[\protect\citeauthoryear{J{\ae}ger, Hustad, Cai, and Tveito}{J{\ae}ger
  et~al\mbox{.}}{2021}]%
        {Jaeger2021}
\bibfield{author}{\bibinfo{person}{K.~H. J{\ae}ger}, \bibinfo{person}{K.~G.
  Hustad}, \bibinfo{person}{X.~Cai}, {and} \bibinfo{person}{A.~Tveito}.}
  \bibinfo{year}{2021}\natexlab{}.
\newblock \bibinfo{booktitle}{\emph{Operator Splitting and Finite Difference
  Schemes for Solving the EMI Model}}.
\newblock \bibinfo{publisher}{Springer International Publishing},
  \bibinfo{address}{Cham}, \bibinfo{pages}{44--55}.
\newblock
\showISBNx{978-3-030-61157-6}
\urldef\tempurl%
\url{https://doi.org/10.1007/978-3-030-61157-6_4}
\showDOI{\tempurl}


\bibitem[\protect\citeauthoryear{Kedarasetti, Turner, Echagarruga, Gluckman,
  Drew, and Constanzo}{Kedarasetti et~al\mbox{.}}{2020}]%
        {constanzo2020}
\bibfield{author}{\bibinfo{person}{R.~T. Kedarasetti}, \bibinfo{person}{K.~L.
  Turner}, \bibinfo{person}{C.~Echagarruga}, \bibinfo{person}{B.~J. Gluckman},
  \bibinfo{person}{P.~J. Drew}, {and} \bibinfo{person}{F.~Constanzo}.}
  \bibinfo{year}{2020}\natexlab{}.
\newblock \showarticletitle{Functional hyperemia drives fluid exchange in the
  paravascular space}.
\newblock \bibinfo{journal}{\emph{Fluids Barriers CNS}} \bibinfo{volume}{17},
  \bibinfo{number}{52} (\bibinfo{year}{2020}).
\newblock
\urldef\tempurl%
\url{https://doi.org/10.1186/s12987-020-00214-3}
\showDOI{\tempurl}


\bibitem[\protect\citeauthoryear{Kim, Xu, and Zikatanov}{Kim
  et~al\mbox{.}}{2003}]%
        {Kim.H;Xu.J;Zikatanov.L.2003b}
\bibfield{author}{\bibinfo{person}{H.~Kim}, \bibinfo{person}{J.~Xu}, {and}
  \bibinfo{person}{L.~Zikatanov}.} \bibinfo{year}{2003}\natexlab{}.
\newblock \showarticletitle{{A multigrid method based on graph matching for
  convection diffusion equations}}.
\newblock \bibinfo{journal}{\emph{Numerical linear algebra with applications}}
  (\bibinfo{year}{2003}).
\newblock
\urldef\tempurl%
\url{http://onlinelibrary.wiley.com/doi/10.1002/nla.317/abstract}
\showURL{%
\tempurl}


\bibitem[\protect\citeauthoryear{Kirby and Mitchell}{Kirby and
  Mitchell}{2018}]%
        {kirby2018solver}
\bibfield{author}{\bibinfo{person}{R.~C. Kirby} {and}
  \bibinfo{person}{L.~Mitchell}.} \bibinfo{year}{2018}\natexlab{}.
\newblock \showarticletitle{Solver composition across the {PDE}/linear algebra
  barrier}.
\newblock \bibinfo{journal}{\emph{SIAM Journal on Scientific Computing}}
  \bibinfo{volume}{40}, \bibinfo{number}{1} (\bibinfo{year}{2018}),
  \bibinfo{pages}{C76--C98}.
\newblock


\bibitem[\protect\citeauthoryear{Koch, Schneider, Helmig, and Jenny}{Koch
  et~al\mbox{.}}{2020}]%
        {Koch2020}
\bibfield{author}{\bibinfo{person}{T.~Koch}, \bibinfo{person}{M.~Schneider},
  \bibinfo{person}{R.~Helmig}, {and} \bibinfo{person}{P.~Jenny}.}
  \bibinfo{year}{2020}\natexlab{}.
\newblock \showarticletitle{{Modeling tissue perfusion in terms of 1d-3d
  embedded mixed-dimension coupled problems with distributed sources}}.
\newblock \bibinfo{journal}{\emph{J. Comput. Phys.}}  \bibinfo{volume}{410}
  (\bibinfo{year}{2020}).
\newblock
\urldef\tempurl%
\url{https://doi.org/10.1016/j.jcp.2020.109370}
\showDOI{\tempurl}


\bibitem[\protect\citeauthoryear{Kolev and Vassilevski}{Kolev and
  Vassilevski}{2012}]%
        {Kolev2012}
\bibfield{author}{\bibinfo{person}{T.~V. Kolev} {and} \bibinfo{person}{P.~S.
  Vassilevski}.} \bibinfo{year}{2012}\natexlab{}.
\newblock \showarticletitle{Parallel auxiliary space amg solver for H(div)
  problems}.
\newblock \bibinfo{journal}{\emph{SIAM Journal on Scientific Computing}}
  \bibinfo{volume}{34} (\bibinfo{year}{2012}), \bibinfo{pages}{1--21}.
\newblock
Issue 6.
\showISSN{10648275}
\urldef\tempurl%
\url{https://doi.org/10.1137/110859361}
\showDOI{\tempurl}


\bibitem[\protect\citeauthoryear{K{\"o}ppl, Vidotto, Wohlmuth, and
  Zunino}{K{\"o}ppl et~al\mbox{.}}{2018}]%
        {Koppl2018}
\bibfield{author}{\bibinfo{person}{T.~K{\"o}ppl}, \bibinfo{person}{E.~Vidotto},
  \bibinfo{person}{B.~Wohlmuth}, {and} \bibinfo{person}{P.~Zunino}.}
  \bibinfo{year}{2018}\natexlab{}.
\newblock \showarticletitle{{Mathematical modeling, analysis and numerical
  approximation of second-order elliptic problems with inclusions}}.
\newblock \bibinfo{journal}{\emph{Mathematical Models and Methods in Applied
  Sciences}} \bibinfo{volume}{28}, \bibinfo{number}{5} (\bibinfo{year}{2018}),
  \bibinfo{pages}{953--978}.
\newblock


\bibitem[\protect\citeauthoryear{Kraus and Margenov}{Kraus and
  Margenov}{2009}]%
        {Johannesbook}
\bibfield{author}{\bibinfo{person}{J.~Kraus} {and}
  \bibinfo{person}{S.~Margenov}.} \bibinfo{year}{2009}\natexlab{}.
\newblock \bibinfo{booktitle}{\emph{Robust algebraic multilevel methods and
  algorithms}}. \bibinfo{series}{Radon series on computational and applied
  mathematics}, Vol.~\bibinfo{volume}{5}.
\newblock \bibinfo{publisher}{de Gruyter}, \bibinfo{address}{Berlin, New York}.
\newblock
\showISBNx{978-3-11-019365-7}


\bibitem[\protect\citeauthoryear{Kraus}{Kraus}{2002}]%
        {K_02}
\bibfield{author}{\bibinfo{person}{J.~K. Kraus}.}
  \bibinfo{year}{2002}\natexlab{}.
\newblock \showarticletitle{An algebraic preconditioning method for
  {$M$}-matrices: linear versus non-linear multilevel iteration}.
\newblock \bibinfo{journal}{\emph{Numer. Linear Algebra Appl.}}
  \bibinfo{volume}{9}, \bibinfo{number}{8} (\bibinfo{year}{2002}),
  \bibinfo{pages}{599--618}.
\newblock
\showCODEN{NLAAEM}
\showISSN{1070-5325}
\urldef\tempurl%
\url{https://doi.org/10.1002/nla.281}
\showDOI{\tempurl}


\bibitem[\protect\citeauthoryear{Kuchta}{Kuchta}{2021}]%
        {fenicsii}
\bibfield{author}{\bibinfo{person}{M.~Kuchta}.}
  \bibinfo{year}{2021}\natexlab{}.
\newblock \showarticletitle{Assembly of Multiscale Linear PDE Operators}. In
  \bibinfo{booktitle}{\emph{Numerical Mathematics and Advanced Applications
  ENUMATH 2019}}, \bibfield{editor}{\bibinfo{person}{Fred~J. Vermolen} {and}
  \bibinfo{person}{Cornelis Vuik}} (Eds.). \bibinfo{publisher}{Springer
  International Publishing}, \bibinfo{address}{Cham},
  \bibinfo{pages}{641--650}.
\newblock
\showISBNx{978-3-030-55874-1}
\urldef\tempurl%
\url{https://doi.org/10.1007/978-3-030-55874-1_63}
\showDOI{\tempurl}


\bibitem[\protect\citeauthoryear{Kuchta, Laurino, Mardal, and Zunino}{Kuchta
  et~al\mbox{.}}{2021}]%
        {Kuchta2021}
\bibfield{author}{\bibinfo{person}{M.~Kuchta}, \bibinfo{person}{F.~Laurino},
  \bibinfo{person}{K.~A. Mardal}, {and} \bibinfo{person}{P.~Zunino}.}
  \bibinfo{year}{2021}\natexlab{}.
\newblock \showarticletitle{{Analysis and approximation of mixed-dimensional
  PDEs on 3D-1D domains coupled with Lagrange multipliers}}.
\newblock \bibinfo{journal}{\emph{SIAM J. Numer. Anal.}} \bibinfo{volume}{59},
  \bibinfo{number}{1} (\bibinfo{year}{2021}), \bibinfo{pages}{558--582}.
\newblock
\urldef\tempurl%
\url{https://doi.org/10.1137/20M1329664}
\showDOI{\tempurl}


\bibitem[\protect\citeauthoryear{Laurino and Zunino}{Laurino and
  Zunino}{2019}]%
        {laurino2019derivation}
\bibfield{author}{\bibinfo{person}{F.~Laurino} {and}
  \bibinfo{person}{P.~Zunino}.} \bibinfo{year}{2019}\natexlab{}.
\newblock \showarticletitle{Derivation and analysis of coupled {PDE}s on
  manifolds with high dimensionality gap arising from topological model
  reduction}.
\newblock \bibinfo{journal}{\emph{ESAIM: Mathematical Modelling and Numerical
  Analysis}} \bibinfo{volume}{53}, \bibinfo{number}{6} (\bibinfo{year}{2019}),
  \bibinfo{pages}{2047--2080}.
\newblock


\bibitem[\protect\citeauthoryear{Layton, Schieweck, and Yotov}{Layton
  et~al\mbox{.}}{2002}]%
        {layton2002coupling}
\bibfield{author}{\bibinfo{person}{W.~J. Layton},
  \bibinfo{person}{F.~Schieweck}, {and} \bibinfo{person}{I.~Yotov}.}
  \bibinfo{year}{2002}\natexlab{}.
\newblock \showarticletitle{Coupling fluid flow with porous media flow}.
\newblock \bibinfo{journal}{\emph{SIAM J. Numer. Anal.}} \bibinfo{volume}{40},
  \bibinfo{number}{6} (\bibinfo{year}{2002}), \bibinfo{pages}{2195--2218}.
\newblock


\bibitem[\protect\citeauthoryear{Lee, Wu, Xu, and Zikatanov}{Lee
  et~al\mbox{.}}{2007}]%
        {Lee2007}
\bibfield{author}{\bibinfo{person}{Y.~Lee}, \bibinfo{person}{J.~Wu},
  \bibinfo{person}{J.~Xu}, {and} \bibinfo{person}{L.~T. Zikatanov}.}
  \bibinfo{year}{2007}\natexlab{}.
\newblock \showarticletitle{{Robust subspace correction methods for nearly
  singular systems}}.
\newblock \bibinfo{journal}{\emph{Mathematical Models and Methods in Applied
  Sciences}} \bibinfo{volume}{17}, \bibinfo{number}{11} (\bibinfo{year}{2007}),
  \bibinfo{pages}{1937--1963}.
\newblock


\bibitem[\protect\citeauthoryear{Livne and Brandt}{Livne and Brandt}{2012}]%
        {Livne.O;Brandt.A.2012a}
\bibfield{author}{\bibinfo{person}{O.~Livne} {and}
  \bibinfo{person}{A.~Brandt}.} \bibinfo{year}{2012}\natexlab{}.
\newblock \showarticletitle{Lean algebraic multigrid {(LAMG)}: fast graph
  {L}aplacian linear solver}.
\newblock \bibinfo{journal}{\emph{SIAM Journal on Scientific Computing}}
  \bibinfo{volume}{34}, \bibinfo{number}{4} (\bibinfo{year}{2012}),
  \bibinfo{pages}{499--523}.
\newblock


\bibitem[\protect\citeauthoryear{Logg, Mardal, Wells, et~al\mbox{.}}{Logg
  et~al\mbox{.}}{2012}]%
        {fenics_book}
\bibfield{author}{\bibinfo{person}{A.~Logg}, \bibinfo{person}{K.-A. Mardal},
  \bibinfo{person}{G.~N. Wells}, {et~al\mbox{.}}}
  \bibinfo{year}{2012}\natexlab{}.
\newblock \bibinfo{booktitle}{\emph{Automated Solution of Differential
  Equations by the Finite Element Method}}.
\newblock \bibinfo{publisher}{Springer}, \bibinfo{address}{Berlin Heidelberg}.
\newblock
\showISBNx{978-3-642-23098-1}
\urldef\tempurl%
\url{https://doi.org/10.1007/978-3-642-23099-8}
\showDOI{\tempurl}


\bibitem[\protect\citeauthoryear{Mardal and Haga}{Mardal and Haga}{2012}]%
        {mardal2012block}
\bibfield{author}{\bibinfo{person}{K.-A. Mardal} {and} \bibinfo{person}{J.~B.
  Haga}.} \bibinfo{year}{2012}\natexlab{}.
\newblock \showarticletitle{{Block preconditioning of systems of PDEs}}.
\newblock In \bibinfo{booktitle}{\emph{Automated solution of differential
  equations by the finite element method}}. \bibinfo{publisher}{Springer},
  \bibinfo{address}{Berlin Heidelberg}, \bibinfo{pages}{643--655}.
\newblock
\urldef\tempurl%
\url{https://doi.org/10.1007/978-3-642-23099-8_35}
\showDOI{\tempurl}


\bibitem[\protect\citeauthoryear{Mardal, Rognes, Thompson, and Valnes}{Mardal
  et~al\mbox{.}}{2022}]%
        {mri2fem}
\bibfield{author}{\bibinfo{person}{K.-A. Mardal}, \bibinfo{person}{M.~E.
  Rognes}, \bibinfo{person}{T.~B. Thompson}, {and} \bibinfo{person}{L.~M.
  Valnes}.} \bibinfo{year}{2022}\natexlab{}.
\newblock \bibinfo{booktitle}{\emph{Mathematical modeling of the human brain:
  from magnetic resonance images to finite element simulation}}.
\newblock \bibinfo{publisher}{Springer}, \bibinfo{address}{Berlin Heidelberg}.
\newblock


\bibitem[\protect\citeauthoryear{Mardal and Winther}{Mardal and
  Winther}{2011}]%
        {mardal_winther_2011}
\bibfield{author}{\bibinfo{person}{K.-A. Mardal} {and}
  \bibinfo{person}{R.~Winther}.} \bibinfo{year}{2011}\natexlab{}.
\newblock \showarticletitle{{Preconditioning discretizations of systems of
  partial differential equations}}.
\newblock \bibinfo{journal}{\emph{Numerical Linear Algebra with Applications}}
  \bibinfo{volume}{18}, \bibinfo{number}{1} (\bibinfo{year}{2011}),
  \bibinfo{pages}{1--40}.
\newblock
\showISSN{10705325}
\urldef\tempurl%
\url{https://doi.org/10.1002/nla.716}
\showDOI{\tempurl}


\bibitem[\protect\citeauthoryear{Marek}{Marek}{1991}]%
        {Marek.I.1991a}
\bibfield{author}{\bibinfo{person}{I.~Marek}.} \bibinfo{year}{1991}\natexlab{}.
\newblock \showarticletitle{{Aggregation methods of computing stationary
  distributions of Markov processes}}.
\newblock In \bibinfo{booktitle}{\emph{Numerical Treatment of Eigenvalue
  Problems Vol. 5/Numerische Behandlung von Eigenwertaufgaben Band 5}}.
  \bibinfo{publisher}{Springer}, \bibinfo{pages}{155--169}.
\newblock


\bibitem[\protect\citeauthoryear{M\'{i}ka and Van{\v{e}}k}{M\'{i}ka and
  Van{\v{e}}k}{1992}]%
        {Mika.S;Vanek.P.1992b}
\bibfield{author}{\bibinfo{person}{S.~M\'{i}ka} {and}
  \bibinfo{person}{P.~Van{\v{e}}k}.} \bibinfo{year}{1992}\natexlab{}.
\newblock \showarticletitle{Acceleration of convergence of a two level
  algebraic algorithm by aggregation in smoothing process}.
\newblock \bibinfo{journal}{\emph{Appl. Math.}}  \bibinfo{volume}{37}
  (\bibinfo{year}{1992}), \bibinfo{pages}{343--356}.
\newblock


\bibitem[\protect\citeauthoryear{M{\'\i}ka and Van{\v{e}}k}{M{\'\i}ka and
  Van{\v{e}}k}{1992}]%
        {Mika.S;Vanek.P.1992a}
\bibfield{author}{\bibinfo{person}{S.~M{\'\i}ka} {and}
  \bibinfo{person}{P.~Van{\v{e}}k}.} \bibinfo{year}{1992}\natexlab{}.
\newblock \showarticletitle{A Modification of the two-level algorithm with
  overcorrection}.
\newblock \bibinfo{journal}{\emph{Appl. Math.}}  \bibinfo{volume}{37}
  (\bibinfo{year}{1992}), \bibinfo{pages}{13--28}.
\newblock


\bibitem[\protect\citeauthoryear{Nakatsukasa, S{\`e}te, and
  Trefethen}{Nakatsukasa et~al\mbox{.}}{2018}]%
        {Nakatsukasa2018}
\bibfield{author}{\bibinfo{person}{Y.~Nakatsukasa},
  \bibinfo{person}{O.~S{\`e}te}, {and} \bibinfo{person}{L.~N. Trefethen}.}
  \bibinfo{year}{2018}\natexlab{}.
\newblock \showarticletitle{The {{AAA Algorithm}} for {{Rational
  Approximation}}}.
\newblock \bibinfo{journal}{\emph{SIAM Journal on Scientific Computing}}
  \bibinfo{volume}{40}, \bibinfo{number}{3} (\bibinfo{year}{2018}),
  \bibinfo{pages}{A1494--A1522}.
\newblock
\showISSN{1064-8275}
\urldef\tempurl%
\url{https://doi.org/10.1137/16M1106122}
\showDOI{\tempurl}


\bibitem[\protect\citeauthoryear{NeuroMorpho}{NeuroMorpho}{2017}]%
        {neuron_geo}
NeuroMorpho \bibinfo{year}{2017}\natexlab{}.
\newblock \bibinfo{booktitle}{\emph{Digital reconstruction of a neuron, ID
  NMO\_72183.}}
\newblock
\urldef\tempurl%
\url{https://neuromorpho.org/neuron_info.jsp?neuron_name=P14_rat1_layerIII_cell1}
\showURL{%
\tempurl}


\bibitem[\protect\citeauthoryear{Notay}{Notay}{2010}]%
        {Notay.Y.2010b}
\bibfield{author}{\bibinfo{person}{Y.~Notay}.} \bibinfo{year}{2010}\natexlab{}.
\newblock \showarticletitle{{An aggregation-based algebraic multigrid method}}.
\newblock \bibinfo{journal}{\emph{Electronic transactions on numerical
  analysis}}  \bibinfo{volume}{37} (\bibinfo{year}{2010}),
  \bibinfo{pages}{123--146}.
\newblock
\urldef\tempurl%
\url{http://www.emis.ams.org/journals/ETNA/vol.37.2010/pp123-146.dir/pp123-146.pdf}
\showURL{%
\tempurl}


\bibitem[\protect\citeauthoryear{Quarteroni and Valli}{Quarteroni and
  Valli}{1991}]%
        {quarteroni1991theory}
\bibfield{author}{\bibinfo{person}{A.~Quarteroni} {and}
  \bibinfo{person}{A.~Valli}.} \bibinfo{year}{1991}\natexlab{}.
\newblock \showarticletitle{Theory and application of Steklov-Poincar{\'e}
  operators for boundary-value problems}.
\newblock In \bibinfo{booktitle}{\emph{Applied and Industrial Mathematics}}.
  \bibinfo{publisher}{Springer}, \bibinfo{pages}{179--203}.
\newblock


\bibitem[\protect\citeauthoryear{Rathgeber, Ham, Mitchell, Lange, Luporini,
  Mcrae, Bercea, Markall, and Kelly}{Rathgeber et~al\mbox{.}}{2016}]%
        {firedrake}
\bibfield{author}{\bibinfo{person}{F.~Rathgeber}, \bibinfo{person}{D.~A. Ham},
  \bibinfo{person}{L.~Mitchell}, \bibinfo{person}{M.~Lange},
  \bibinfo{person}{F.~Luporini}, \bibinfo{person}{A.~T.~T. Mcrae},
  \bibinfo{person}{G.-T. Bercea}, \bibinfo{person}{G.~R. Markall}, {and}
  \bibinfo{person}{P.~H.~J. Kelly}.} \bibinfo{year}{2016}\natexlab{}.
\newblock \showarticletitle{Firedrake: Automating the Finite Element Method by
  Composing Abstractions}.
\newblock \bibinfo{journal}{\emph{ACM Trans. Math. Softw.}}
  \bibinfo{volume}{43}, \bibinfo{number}{3}, Article \bibinfo{articleno}{24}
  (\bibinfo{date}{dec} \bibinfo{year}{2016}), \bibinfo{numpages}{27}~pages.
\newblock
\showISSN{0098-3500}


\bibitem[\protect\citeauthoryear{Sogn and Takacs}{Sogn and Takacs}{2022}]%
        {sogn2022stable}
\bibfield{author}{\bibinfo{person}{J.~Sogn} {and} \bibinfo{person}{S.~Takacs}.}
  \bibinfo{year}{2022}\natexlab{}.
\newblock \showarticletitle{Stable discretizations and {IETI-DP} solvers for
  the {S}tokes system in multi-patch {I}sogeometric {A}nalysis}.
\newblock \bibinfo{journal}{\emph{arXiv preprint arXiv:2202.13707}}
  (\bibinfo{year}{2022}).
\newblock


\bibitem[\protect\citeauthoryear{Urschel, Xu, Hu, and Zikatanov}{Urschel
  et~al\mbox{.}}{2015}]%
        {Urschel2015}
\bibfield{author}{\bibinfo{person}{J.~C. Urschel}, \bibinfo{person}{J.~Xu},
  \bibinfo{person}{X.~Hu}, {and} \bibinfo{person}{L.~T. Zikatanov}.}
  \bibinfo{year}{2015}\natexlab{}.
\newblock \showarticletitle{A {Cascadic} {Multigrid} {Algorithm} for
  {Computing} the {Fiedler} {Vector} of {Graph} {Laplacians}}.
\newblock \bibinfo{journal}{\emph{Journal of Computational Mathematics}}
  \bibinfo{volume}{33}, \bibinfo{number}{2} (\bibinfo{date}{June}
  \bibinfo{year}{2015}), \bibinfo{pages}{209--226}.
\newblock
\showISSN{0254-9409, 1991-7139}
\urldef\tempurl%
\url{https://doi.org/10.4208/jcm.1412-m2014-0041}
\showDOI{\tempurl}


\bibitem[\protect\citeauthoryear{Vakhutinsky, Dudkin, and Ryvkin}{Vakhutinsky
  et~al\mbox{.}}{1979}]%
        {Vakhutinsky.I;Dudkin.L;Ryvkin.A.1979a}
\bibfield{author}{\bibinfo{person}{I.~Vakhutinsky},
  \bibinfo{person}{L.~Dudkin}, {and} \bibinfo{person}{A.~Ryvkin}.}
  \bibinfo{year}{1979}\natexlab{}.
\newblock \showarticletitle{{Iterative aggregation--A new approach to the
  solution of large-scale problems}}.
\newblock \bibinfo{journal}{\emph{Econometrica: Journal of the Econometric
  Society}} (\bibinfo{year}{1979}), \bibinfo{pages}{821--841}.
\newblock


\bibitem[\protect\citeauthoryear{Van{\v{e}}k, Mandel, and Brezina}{Van{\v{e}}k
  et~al\mbox{.}}{1996}]%
        {Vanek.P;Mandel.J;Brezina.M.1996a}
\bibfield{author}{\bibinfo{person}{P.~Van{\v{e}}k},
  \bibinfo{person}{J.~Mandel}, {and} \bibinfo{person}{M.~Brezina}.}
  \bibinfo{year}{1996}\natexlab{}.
\newblock \showarticletitle{Algebraic multigrid based on smoothed aggregation
  for second and fourth order problems}.
\newblock \bibinfo{journal}{\emph{Computing}}  \bibinfo{volume}{56}
  (\bibinfo{year}{1996}), \bibinfo{pages}{179--196}.
\newblock


\bibitem[\protect\citeauthoryear{Van{\v{e}}k, Mandel, and Brezina}{Van{\v{e}}k
  et~al\mbox{.}}{1998}]%
        {Vanek.P;Mandel.J;Brezina.M.1998a}
\bibfield{author}{\bibinfo{person}{P.~Van{\v{e}}k},
  \bibinfo{person}{J.~Mandel}, {and} \bibinfo{person}{M.~Brezina}.}
  \bibinfo{year}{1998}\natexlab{}.
\newblock \showarticletitle{Convergence of algebraic multigrid based on
  smoothed aggregation}.
\newblock \bibinfo{journal}{\emph{Computing}}  \bibinfo{volume}{56}
  (\bibinfo{year}{1998}), \bibinfo{pages}{179--196}.
\newblock


\bibitem[\protect\citeauthoryear{Van\v{e}k, Mandel, and Brezina}{Van\v{e}k
  et~al\mbox{.}}{1996}]%
        {1996Vanek_Mandel_Brezina}
\bibfield{author}{\bibinfo{person}{P.~Van\v{e}k}, \bibinfo{person}{J.~Mandel},
  {and} \bibinfo{person}{M.~Brezina}.} \bibinfo{year}{1996}\natexlab{}.
\newblock \showarticletitle{Algebraic multigrid by smoothed aggregation for
  second and fourth order elliptic problems}.
\newblock Vol.~\bibinfo{volume}{56}. \bibinfo{pages}{179--196}.
\newblock
\showISSN{0010-485X}
\urldef\tempurl%
\url{https://doi.org/10.1007/BF02238511}
\showDOI{\tempurl}
\newblock
\shownote{International GAMM-Workshop on Multi-level Methods (Meisdorf, 1994).}


\bibitem[\protect\citeauthoryear{Vassilevski}{Vassilevski}{2008}]%
        {Panayotsbook}
\bibfield{author}{\bibinfo{person}{P.~S. Vassilevski}.}
  \bibinfo{year}{2008}\natexlab{}.
\newblock \bibinfo{booktitle}{\emph{Multilevel block factorization
  preconditioners}}.
\newblock \bibinfo{publisher}{Springer}, \bibinfo{address}{New York}. xiv+529
  pages.
\newblock
\showISBNx{978-0-387-71563-6}


\bibitem[\protect\citeauthoryear{Verdugo and Badia}{Verdugo and Badia}{2022}]%
        {gridap}
\bibfield{author}{\bibinfo{person}{F.~Verdugo} {and}
  \bibinfo{person}{S.~Badia}.} \bibinfo{year}{2022}\natexlab{}.
\newblock \showarticletitle{The software design of Gridap: A Finite Element
  package based on the Julia JIT compiler}.
\newblock \bibinfo{journal}{\emph{Computer Physics Communications}}
  \bibinfo{volume}{276} (\bibinfo{year}{2022}), \bibinfo{pages}{108341}.
\newblock
\showISSN{0010-4655}


\bibitem[\protect\citeauthoryear{Xie, Kang, Xu, Chen, Liao, Thiyagarajan,
  O’Donnell, Christensen, Nicholson, Iliff, et~al\mbox{.}}{Xie
  et~al\mbox{.}}{2013}]%
        {xie2013sleep}
\bibfield{author}{\bibinfo{person}{L.~Xie}, \bibinfo{person}{H.~Kang},
  \bibinfo{person}{Q.~Xu}, \bibinfo{person}{M.~J. Chen},
  \bibinfo{person}{Y.~Liao}, \bibinfo{person}{M.~Thiyagarajan},
  \bibinfo{person}{J.~O’Donnell}, \bibinfo{person}{D.~J. Christensen},
  \bibinfo{person}{C.~Nicholson}, \bibinfo{person}{J.~J. Iliff},
  {et~al\mbox{.}}} \bibinfo{year}{2013}\natexlab{}.
\newblock \showarticletitle{Sleep drives metabolite clearance from the adult
  brain}.
\newblock \bibinfo{journal}{\emph{science}} \bibinfo{volume}{342},
  \bibinfo{number}{6156} (\bibinfo{year}{2013}), \bibinfo{pages}{373--377}.
\newblock


\bibitem[\protect\citeauthoryear{Xu}{Xu}{1992}]%
        {Xu1992SIAMReview}
\bibfield{author}{\bibinfo{person}{J.~Xu}.} \bibinfo{year}{1992}\natexlab{}.
\newblock \showarticletitle{Iterative methods by space decomposition and
  subspace correction}.
\newblock \bibinfo{journal}{\emph{SIAM Rev.}} \bibinfo{volume}{34},
  \bibinfo{number}{4} (\bibinfo{year}{1992}), \bibinfo{pages}{581--613}.
\newblock
\showISSN{0036-1445}
\urldef\tempurl%
\url{https://doi.org/10.1137/1034116}
\showDOI{\tempurl}


\bibitem[\protect\citeauthoryear{Xu and Zikatanov}{Xu and Zikatanov}{2002}]%
        {XuZikatanov2002}
\bibfield{author}{\bibinfo{person}{J.~Xu} {and}
  \bibinfo{person}{L.~Zikatanov}.} \bibinfo{year}{2002}\natexlab{}.
\newblock \showarticletitle{The method of alternating projections and the
  method of subspace corrections in {H}ilbert space}.
\newblock \bibinfo{journal}{\emph{J. Amer. Math. Soc.}} \bibinfo{volume}{15},
  \bibinfo{number}{3} (\bibinfo{year}{2002}), \bibinfo{pages}{573--597}.
\newblock
\showISSN{0894-0347}
\urldef\tempurl%
\url{https://doi.org/10.1090/S0894-0347-02-00398-3}
\showDOI{\tempurl}


\bibitem[\protect\citeauthoryear{Zhao, Hu, Cai, and Karniadakis}{Zhao
  et~al\mbox{.}}{2017}]%
        {ZhaoHuCaiKarniadakis2017a}
\bibfield{author}{\bibinfo{person}{X.~Zhao}, \bibinfo{person}{X.~Hu},
  \bibinfo{person}{W.~Cai}, {and} \bibinfo{person}{G.~E. Karniadakis}.}
  \bibinfo{year}{2017}\natexlab{}.
\newblock \showarticletitle{Adaptive finite element method for fractional
  differential equations using hierarchical matrices}.
\newblock \bibinfo{journal}{\emph{Computer Methods in Applied Mechanics and
  Engineering}}  \bibinfo{volume}{325} (\bibinfo{date}{Oct.}
  \bibinfo{year}{2017}), \bibinfo{pages}{56--76}.
\newblock
\showISSN{0045-7825}
\urldef\tempurl%
\url{https://doi.org/10.1016/j.cma.2017.06.017}
\showDOI{\tempurl}


\end{thebibliography}

\appendix
\section{Darcy-Stokes problem in 2d}\label{sec:ds_2d}

To further illustrate robustness of Darcy-Stokes preconditioner \eqref{eq:ds_precond}
and scalability of its HAZniCS implementation (\Cref{lst:ds_prec}) we consider the experimental
and solver setup of \Cref{subsec:ds_results} in two-dimensions. Namely, we let
$\Omega_S=\left[0, \tfrac{1}{2}\right]\times \left[0, 1\right]$ and
$\Omega_D=\left[\tfrac{1}{2}, 1\right]\times \left[0, 1\right]$. In addition,
\eqref{eq:ds_system} will be discretized in terms
of $(\mathbb{P}_2-\mathbb{P}_1)-(\mathbb{R}\mathbb{T}_0-\mathbb{P}_0)-\mathbb{P}_0$
elements as well as by the non-conforming (stabilized) 
$(\mathbb{C}\mathbb{R}_1-\mathbb{P}_0)-(\mathbb{R}\mathbb{T}_0-\mathbb{P}_0)-\mathbb{P}_0$
elements. As in \Cref{subsec:ds_results} we employ triangulations of $\Omega_S$, $\Omega_D$
whose trace meshes match on the interface $\Gamma$. We remark that on the finest
level of refinement the two discretizations lead to similar number of unknowns
with $N_{\text{dof}}\approx 1.84\cdot 10^{6}$ and $N_{\text{dof}}\approx 1.71\cdot 10^{6}$ for
Taylor-Hood and Crouzeix-Raviart based spaces respectively. Finally, the two-dimensional
setting allows for comparison between the inexact/multilevel based approximation of the
Darcy-Stokes preconditioner, cf. \Cref{lst:ds_prec}, and its realization using LU decomposition
for the leading blocks. Such preconditioner can be defined in HAZniCS as shown in \Cref{lst:ds_prec_lu}.
We note that in both cases the multiplier block uses the rational approximation.

\begin{lstlisting}[style=mystylepython, caption={Implementation of preconditioner \eqref{eq:ds_precond} for Darcy-Stokes problem \eqref{eq:ds_system} using exact inverses for the leading/bulk blocks. Complete code can be found in scripts \texttt{HAZniCS-examples/demo\_darcy\_stokes*.py}
    }, captionpos=b, label={lst:ds_prec_lu}]
  from block.algebraic.petsc import LU
  
  # [...] Setup blocks as of the preconditioner
  # B0, B1, B2, B3 = ....
  # Only Multiplier block will be inexact
  B4 = RA(A, M, parameters=params)
  
  # define the approximate Riesz map
  B = block_diag_mat([LU(B0), LU(B1), LU(B2), LU(B3), B4])
\end{lstlisting}

\begin{figure}
  \centering
  \includegraphics[width=0.8\textwidth]{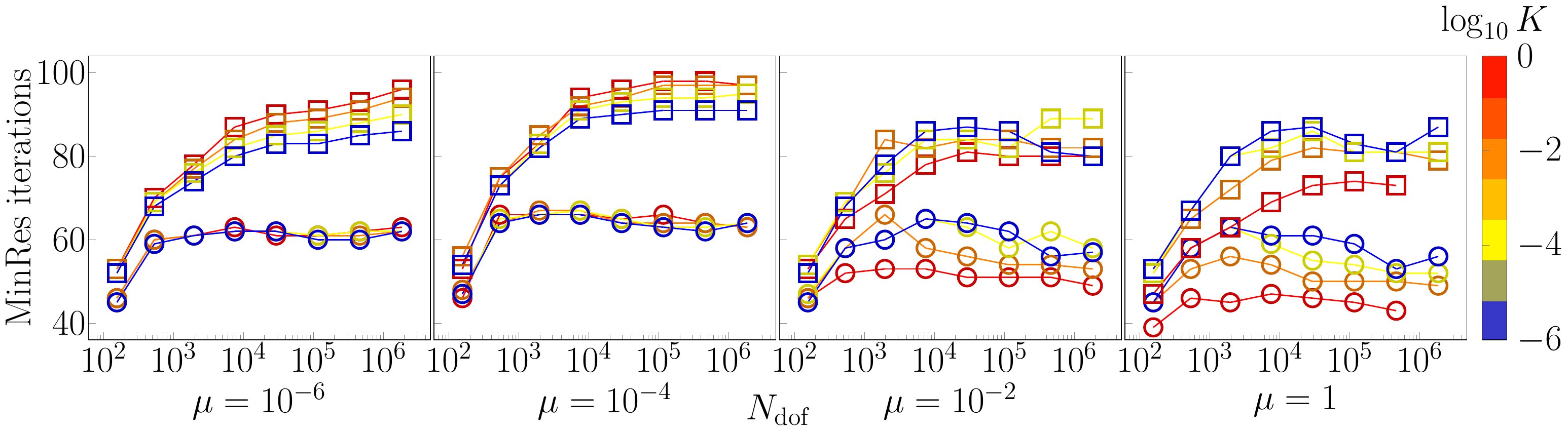}
  \includegraphics[width=0.8\textwidth]{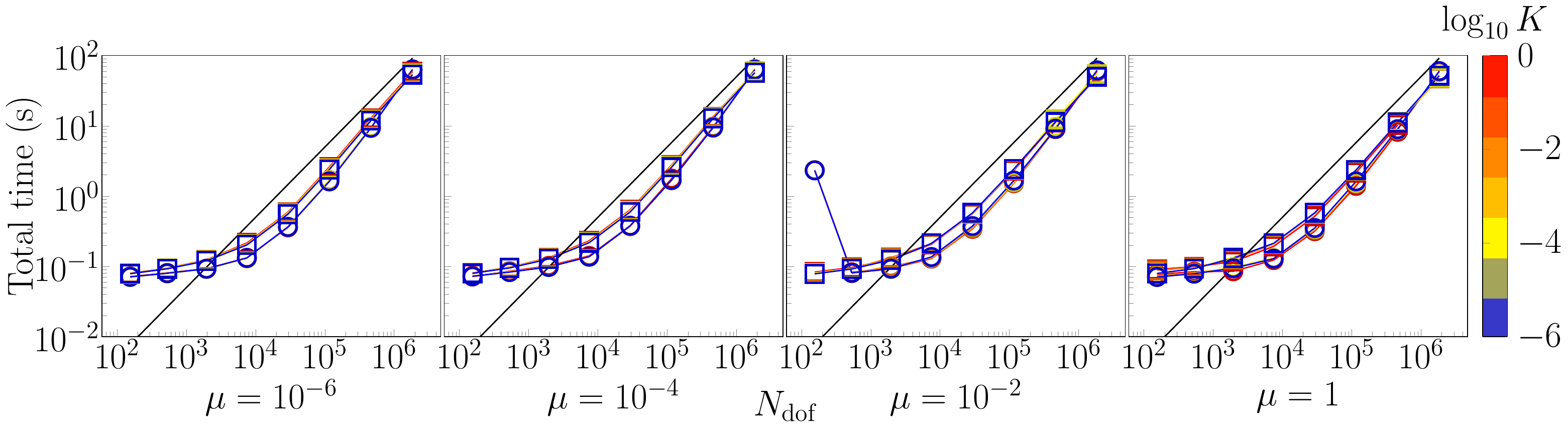}
  \caption{
          Performance of Darcy-Stokes preconditioner \eqref{eq:ds_precond} in case $\Omega_S$, $\Omega_D\subset\mathbb{R}^2$.
Discretization by $(\mathbb{P}_2-\mathbb{P}_1)-(\mathbb{R}\mathbb{T}_0-\mathbb{P}_0)-\mathbb{P}_0$
          elements with $D=0.1$.           
          (Top) Number of MinRes iterations until convergence in relative preconditioned residual norm and tolerance $10^{-12}$
          for different values of $\mu, K$ and mesh sizes. (Bottom) Total solution
          time for solving \eqref{eq:ds_system} including the setup time of the
          preconditioner. Black line indicates linear scaling. In both plots
          data points marked with circles correspond to realization of the preconditioner
          using LU for the leading blocks, see \Cref{lst:ds_prec_lu}, while square markers
          are due to the multilevel approximation in \Cref{lst:ds_prec}.
          Results are obtained by running 
          \texttt{HAZniCS-examples/demo\_darcy\_stokes\_2d\_flat.py} with
          command line switch \texttt{-elm\_family TH}.
  }
  \label{fig:ds_2d_TH}
\end{figure}

\begin{figure}
  \centering
  \includegraphics[width=0.8\textwidth]{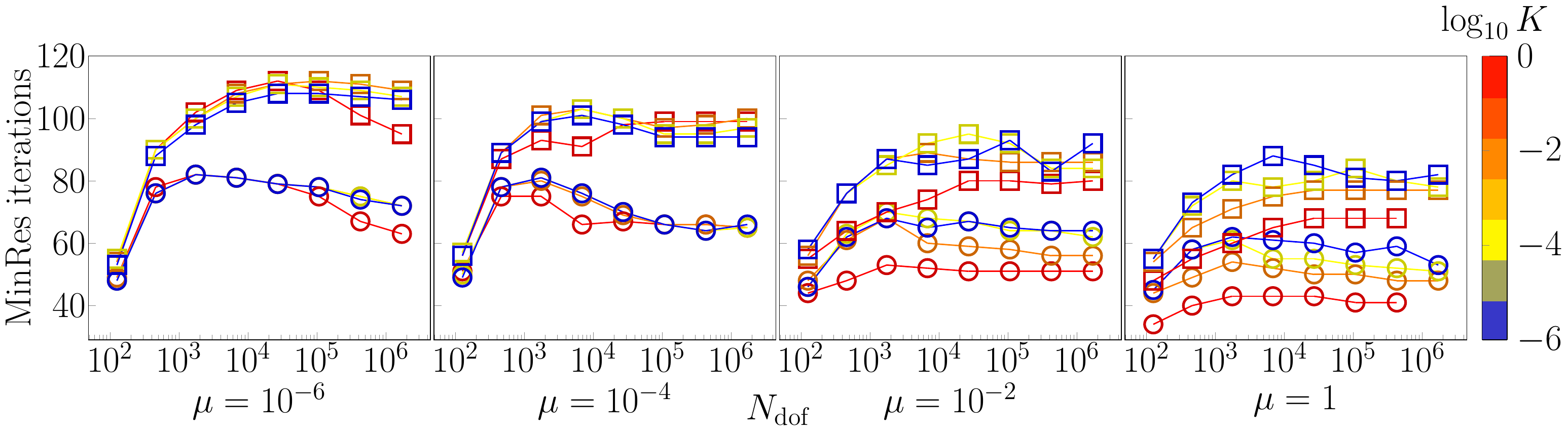}
  \includegraphics[width=0.8\textwidth]{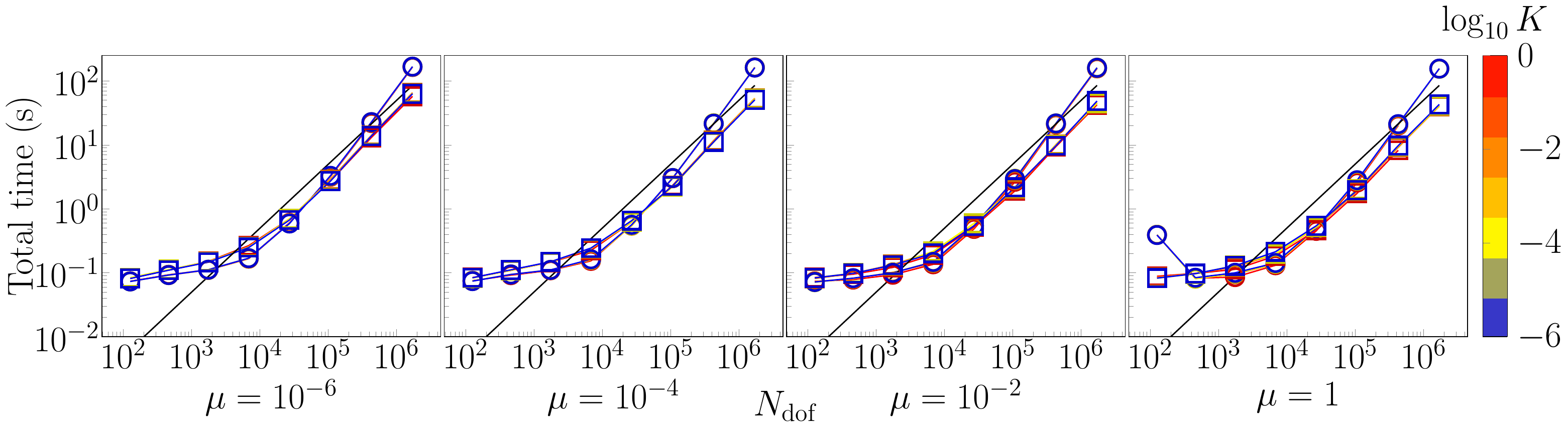}
  \caption{
          Performance of Darcy-Stokes preconditioner \eqref{eq:ds_precond} in case $\Omega_S$, $\Omega_D\subset\mathbb{R}^2$.
Discretization by $(\mathbb{C}\mathbb{R}_1-\mathbb{P}_0)-(\mathbb{R}\mathbb{T}_0-\mathbb{P}_0)-\mathbb{P}_0$
          elements with $D=0.1$.           
          (Top) Number of MinRes iterations until convergence in relative preconditioned residual norm and tolerance $10^{-12}$
          for different values of $\mu, K$ and mesh sizes. (Bottom) Total solution
          time for solving \eqref{eq:ds_system} including the setup time of the
          preconditioner. Black line indicates linear scaling. In both plots
          data points marked with circles correspond to realization of the preconditioner
          using LU for the leading blocks, see \Cref{lst:ds_prec_lu}, while square markers
          are due to the multilevel approximation in \Cref{lst:ds_prec}.
          Results are obtained by running 
          \texttt{HAZniCS-examples/demo\_darcy\_stokes\_2d\_flat.py} with
          command line switch \texttt{-elm\_family CR}.
  }
  \label{fig:ds_2d_CR}
\end{figure}

Performance of the two Darcy-Stokes preconditioners using Taylor-Hood-
and Crouzeix-Raviart based discretizations is shown \Cref{fig:ds_2d_TH} and
\Cref{fig:ds_2d_CR} respectively. In all cases we observe that the number of
MinRes iterations is bounded in mesh size and parameters $\mu$ and $K$. The exact
preconditioners lead to convergence in fewer iterations, e.g. for $K=1$, $\mu=10^{-6}$
the difference is 30 iterations. However, the total solution time is smaller with the
approximate preconditioners using multilevel methods for $\bm{V}_S$ and $\bm{V}_D$
blocks. Moreover, it can be seen that multigrid leads to (close to) optimal scalability of
the preconditioner while the scaling of the exact preconditioner becomes suboptimal. This is
especially the case for the finest meshes and $\mathbb{C}\mathbb{R}$ elements where
$\text{dim}\bm{V}_S=787968$, $\text{dim}{Q}_S=262144$, $\text{dim}\bm{V}_D=393984$,
$\text{dim}{Q}_D=262144$, $\text{dim}\Lambda=512$.

\end{document}